\documentclass{article}

\usepackage{amsmath,amssymb,amsthm,mathrsfs,epsf,a4,graphicx,eucal}
\newtheorem{theorem}{Theorem}[section]

\newtheorem{definition}[theorem]{Definition}

\newtheorem{remark}[theorem]{Remark}

\numberwithin{equation}{section}

\def\version{30.12.2020}\def\users{}  %
\usepackage{psfrag} % it allows for latex fonts in figures
\usepackage{graphicx}

\usepackage{fullpage,bm,amsmath,hyperref,amsfonts,amssymb,color}
\usepackage{mathrsfs} % THIS PACKAGE ALLOWS FOR ``scr'' FONTS
\usepackage{dutchcal}
\DeclareMathAlphabet{\mathdutchcal}{U}{dutchcal}{m}{n}
\DeclareMathAlphabet{\mathpzc}{OT1}{pzc}{m}{it} % small case caligraphic letters %            ~~~~~~~~
\def\users{final-layout}  % when activated, ``our'' debugging is suppressed
\usepackage{ifthen}
\ifthenelse{\equal{\users}{final-layout}}{}{
\usepackage{fancyhdr}
\pagestyle{fancy}
\headheight=28pt\headwidth=17cm
\definecolor{gray}{gray}{0.5}
\rhead{\color{gray}Stefan problem thermomechanically\\
T.Roub\'\i\v cek}
\chead{}
\lhead{Version \version, file: \jobname.tex, %\underline{\Large\color{red}\WHO's working on}
  \\
compiled:
\number\day.\number\month.\number\year\ at
\the\hour:\ifnum\minute<10 0\fi\the\minute\ h }
}

\definecolor{labelkey}{rgb}{1.,.2,0.}

\def\vv{{\bm v}}
\def\uu{{\bm u}}
\def\nn{{\bm n}}
\def\ff{{\bm f}}

\def\vvk{\vv_\etau^k}
  \def\vvkk{\vv_\etau^{k-1}}
\def\alphak{\alpha_\etau^k}
\def\alphal{\alpha_\etau^l}
  \def\alphakk{\alpha_\etau^{k-1}}
  \def\Eek{\Ee_\etau^k}
  \def\Eekk{\Ee_\etau^{k-1}}
 \def\Eetau{\Ee_\etau^{}}
\def\overlineEetau{\hspace*{.2em}\overline{\hspace*{-.2em}\bm E}_\etau^{}}
\def\overlineStau{\hspace*{.2em}\overline{\hspace*{-.2em}\bm S}_\etau^{}}
\def\overlineSstrtau{\hspace*{.2em}{\overline{\hspace*{-.2em}\bm K\hspace*{-.1em}}\hspace*{.1em}}_\etau^{}}
\def\overlineTtau{\hspace*{.2em}{\overline{\hspace*{-.2em}\bm T\hspace*{-.1em}}\hspace*{.1em}}_\etau^{}}

\def\overlineEptau{\hspace*{.3em}\overline{\hspace*{-.3em}\Ep\hspace*{-.1em}}_{\,\etau^{}}}

\def\overlineRRtau{\hspace*{.2em}\overline{\hspace*{-.2em}\RR}_{\etau}^{}}
\def\overlinealphatau{\hspace*{.15em}\overline{\hspace*{-.15em}\alpha}_{\etau}^{}}
\def\overlinechitau{\hspace*{.1em}\overline{\hspace*{-.1em}\chi}_{\etau}^{}}
\def\overlinevvtau{\hspace*{.15em}\overline{\hspace*{-.15em}\vv}_{\etau}^{}}
\def\Epk{{\Ep}_\etau^k}

\def\NU{\omega}
\newcommand{\INT}[2]{\text{$\displaystyle{\int}$}_{\hspace*{-.6em}#1}^{\hspace*{-.0em}#2}}

\def\Vdots{\!\mbox{\setlength{\unitlength}{1em}
\begin{picture}(0,0)
\put(-.07,0){.}
\put(-.07,.3){.}
\put(-.07,.6){.}
\end{picture}
}
}

\newcommand\DELETEDELETE[1]{}
\newcommand\DT[1]{\mathchoice
  {{\buildrel{\hspace*{.1em}\text{\Large\bf.}}\over{#1}}}
    {{\buildrel{\hspace*{.1em}\text{\large\bf.}}\over{#1}}}
                 {{\buildrel{\hspace*{.1em}\text{\Large.}}\over{#1}}}
                 {{\buildrel{\hspace*{.1em}\text{\large.}}\over{#1}}}}

\newcommand\pdt[1]{\frac{\partial{#1}}{\partial t}} %Partial Derivative with respect to T
\newcommand\Ee{{\bm E}}              %Eulerian elastic strain
\newcommand\Ep{{\bm\varPi}}              %Eulerian plastic strain

\def\widetildeEp{\hspace*{.1em}\widetilde{\hspace*{-.1em}\Ep}}
\newcommand\R{\mathbb{R}}

\renewcommand\d{{\rm d}}
\newcommand\EE{{\bm e}}
\newcommand{\lineunder}[2]{\LU{\begin{array}[t]{c}\underbrace{#1}\vspace*{.5em}\end{array}}{\mbox{\footnotesize\rm #2}}}
\newcommand{\LU}[2]{\begin{array}[t]{c}#1\vspace*{-1em}\\_{#2}\end{array}}
\newcommand{\linesunder}[3]{\LSU{\begin{array}[t]{c}\underbrace{#1}\vspace*{.5em}\end{array}}{\mbox{\footnotesize\rm #2}}{\mbox{\footnotesize\rm #3}}}
\newcommand{\LSU}[3]{\begin{array}[t]{c}#1\vspace*{-1em}\\_{#2}\vspace*{-.5em}\\_{#3}\end{array}}

\newcommand{\divS}{\mathrm{div}_{\scriptscriptstyle\textrm{\hspace*{-.1em}S}}^{}}

\newcommand{\nablaS}{\nabla_{\scriptscriptstyle\textrm{\hspace*{-.3em}S}}^{}}
\newcommand\ZJ[1]{\mathchoice
                 {{\buildrel{\hspace*{.1em}{_{\,\Large\boldsymbol\circ}}}\over{#1}}}
                 {{\buildrel{\hspace*{.1em}{_{\,\large\boldsymbol\circ}}}\over{#1}}}
                 {{\buildrel{\hspace*{.1em}{\boldsymbol\circ}}\over{#1}}}
                 {{\buildrel{\hspace*{.1em}{\boldsymbol\circ}}\over{#1}}}}

\newcommand{\bbC}{\mathbb{C}}
\newcommand{\bbK}{\mbox{\footnotesize$\mathscr{K}$}}
\newcommand{\bbI}{\mathbb{I}}

\newcommand{\bbD}{\mathbb D}

\newcommand{\bmE}{\bm{\varTheta}}

\newcommand{\W}{w}
\newcommand{\GM}{G_\text{\sc m}^{}}
\newcommand{\RR}{\ZJEp}
\newcommand{\eps}{\varepsilon}
\newcommand{\etau}{{\eps\tau}}
\newcommand{\barOmega}{\,\overline{\!\varOmega}}

\newcounter{myfigure}
\newenvironment{my-picture}[3]{\refstepcounter{myfigure}\label{#3}\setlength{\unitlength}{1em}\begin{picture}(#1,#2)}{\end{picture}}

\begin{document}

\def\TTT{\color{black}}
\def\EEE{\color{black}}

\noindent
%\title
{\LARGE\bf The Stefan problem in a thermomechanical context\\[.3em]with fracture and
  fluid flow}
{\Large\footnote{This research has been partially supported also from the
CSF (Czech Science Foundation) project 19-04956S, the M\v SMT \v CR
(Ministry of Education of the Czech Republic) project
CZ.02.1.01/0.0/0.0/15-003/0000493, and the institutional support RVO:
61388998 (\v CR). And, according the publisher's demand, the author declares
that ``this work does not have any conflicts of interest'' and that, not having
any juristical education, he is not completely sure what does such sentence
mean in all details.}}

\bigskip\bigskip

\noindent
       {\large\bf Tom\'a\v s Roub\'\i\v cek}

       \bigskip\medskip

 \noindent 
       {Mathematical Institute, Charles University,\\Sokolovsk\'a 83,
CZ--186~75~Praha~8, %\\\country
{Czech Republic},\footnote{
Email: {\tt tomas.roubicek@mff.cuni.cz},\ \ ORCID: 0000-0002-0651-5959.}
\\and\\
 \noindent      %\address[2]
     {Institute of Thermomechanics, Czech Academy of
       Sciences,\\ Dolej\v skova~5, CZ--182 00 Praha 8, %\\\country
       {Czech Republic}
               }

     \bigskip \bigskip

     %\abstract%[Summary]
\begin{minipage}[c]{40em}
  {\baselineskip=11pt
    {\bf Summary}\\
 The classical Stefan problem, concerning mere heat-transfer
during solid-liquid phase transition, is here enhanced towards
mechanical effects. The Eulerian description  at large displacements 
is used with convective and Zaremba-Jaumann corotational time derivatives,
linearized by using the additive Green-Naghdi's decomposition in 
(objective) rates.  In particular, the liquid phase  is a viscoelastic
fluid  while 
creep and rupture of the solid phase is considered in the Jeffreys
viscoelastic rheology exploiting the phase-field model and a concept
of slightly (so-called ``semi'') compressible materials. The $L^1$-theory
for the heat equation is adopted for the Stefan problem relaxed
by allowing for kinetic superheating/supercooling effects during the
solid-liquid phase transition. A rigorous proof of existence of
weak solutions is provided for an incomplete melting, employing
a time-discretisation approximation.}

  \smallskip
  %\keywords
      {\baselineskip=10pt
  {\bf Keywords:} Solid-liquid phase transition, 
  creep, 
  Jeffreys rheology, 
  semi-compressible fluids, 
  Eulerian formulation, 
  Stefan problem,
  melting, 
  solidification, 
  enthalpy formulation,
  phase-field fracture, 
  fully convective model.
  \newline\vspace*{-.5em}
  \newline
  AMS Classification:
35Q74,  % PDEs in connection with mechanics of deformable solids
35R37,  % Moving boundary problems for PDEs
74A30,  % Nonsimple materials
%74C10,  % Plastic materials - internal variable - small-strain, rate-dependent theories
74R20,  % Anelastic fracture and damage
76A10,  % Viscoelastic fluids
80A22. %Stefan problems, phase changes, etc.
}

\end{minipage}

%\maketitle

%\begin{document}

\def\TTT{\color{black}}
\def\EEE{\color{black}}

\allowdisplaybreaks

%\begin{abstract}
%\end{abstract}

%\end{frontmatter}

\baselineskip=12.5pt

\section{Introduction}
%        ~~~~~~~~~~~~

The formulation of {\it solid-liquid phase transitions} has a long history
and had been widely scrutinized under the name of the Stefan problem
\cite{Stef89PTW},
being a prominent free- or moving-boundary problem in the second half of the
20th century, cf.\ the survey \cite{Visi08ISTP}. The original motivation of
Josef Stefan was to model melting of ice in polar caps of the Earth into
merely static water (considering the heat-transfer problem in one dimension
only). Actually, the {\it melting/freezing} ({\it solidification}) phenomena
can be found in other occasions in our planet Earth: melting/solidification
of Iron with Nickel (in solid inner core vs fluidic outer core) or, without
a sharp phase-transition temperature, of rocks vs magma (volcanism in the
crust and magma migration in the mantle). Moreover, this phenomena is
important in other parts of our Solar system, in particular for  icy 
moons of Jupiter (Europa and Ganymede) and of Saturn (Enceladus  and Titan).
Beside these geophysical applications, other motivations can certainly
be found in engineering, e.g.\ recrystallization and fusing of steel  and of
other metals.

In contrast to the original models of the Stefan problem dealing with mere
temperature (or enthalpy) evolution, there are some mechanical phenomena
with which this thermal problem is ultimately coupled.
In particular, the liquid phase is really fluidic and can easily flow.
In turn, the solid phase may exhibit {\it creep} and damage or {\it fracture},
as often seen on the Earth's ice sheet \cite{SchDuv09CFI} or in the lithospheric
rocks (as so-called aseismic slip and tectonic earthquakes, respectively).
Fracture can also be clearly seen on the surface of 
the mentioned Jupiter's and Saturn's moons. 
Involving creep by a permanent inelastic strain is important
for a longer-time scales and for healing (reversible damage).
This allows us to model {\it Maxwellian or Jeffreys visco-elastic
rheologies} in the solid part.

Rather surprisingly, so far the Stefan problem was scrutinized in mechanical
context only in very simplified situations, like solidification with the
incompressible Navier-Stokes and completely rigid solid with only simplified
thermomechanical coupling e.g.\ in
\cite{FukKen05SPCG,HinZie07OCFB,RodUrb02TDCS,XuShi97SPCJ} or a Darcy flow
independent of a possible phase transformation in \cite{RodUrb99DSPA}. Thus
realistic models for e.g.\ floating and melting ice floats in  oceans, which
was the original motivation of Josef Stefan for his 1-dimensional heat-transfer
problem, was not yet successfully formulated and analysed.

Such a thermo-mechanical model will be formulated here in
Section~\ref{sec-fuller}, after preparatory calculations for casting first
a merely visco-elasto-dynamical model in a
{\it convective stress/velocity formulation} in Section~\ref{sec-preliminary}
 to explain how the corotational time derivative for strains in isotropic
materials arises. 
A specialization of general model towards solid-fluid phase transformation
and {\it phase-field fracture} in solid will then be stated in
Section~\ref{sec-special} and involvement of {\it superheating/cooling
effects} into the model in Section~\ref{sec-super}.
The model from Section~\ref{sec-fuller} in the modification
from Section~\ref{sec-super} will be then analyzed by a converging
time-discretisation combined with a certain regularization
in Section~\ref{sec-anal}, but only for an incomplete melting. 
Eventually, various comments and outlined
enhancements of this model will be outlined in Section~\ref{sec-rem}.
 For geophysical modelling of rock-magma phase transition
see \cite{Roub21TCMP}.

\section{A preliminary model: 
  visco-elastodynamics\\in Eulerian stress-velocity formulation}\label{sec-preliminary}
%        ~~~~~~~~~~~~~~~~~~~~~~~~~~~~~~~~~~~~~~~~~~~~~~~~~~~~

There are several aspects which should be reflected in the model. First,
the so-called {\it fluid-structure interaction} is itself a difficult
problem because fluids dictate Eulerian description while solids
are (believed to be) better treated in Lagrangian coordinates.
We however
devise the model in the {\it fully Eulerian description}. This needs 
a careful concept of the {\it objective time derivatives}. 
In fluidic part, it is also more natural to formulate the model
in terms of stresses and velocities  rather 
than in strains and displacements. 
It then suggests to use the stress/velocity formulation in the solid
part, too.

As this ``monolithic'' Eulerian formulation is not standard, 
let us first introduce it for the mere
visco-elastodynamics without considering any internal variable, which is
perhaps of its own importance itself and which justifies the
strain/velocity formulation used later. As to the mentioned objective time
derivatives, for the transport of scalars there is no controversy and the
{\it convective} (material) {\it time derivative} is used when scalars
represent an {\it intensive variable} (like pressure, 
damage phase-field variable, or temperature). For tensors (or vectors), there
are however many options used in literature, like Oldroyd's or
Truesdell's or Lie's derivatives.
Objectivity of rates here means invariance with respect to moving frame of
reference. It also is reasonable to require the tensor time derivatives
(in particular stress rate), to be so-called {\it identical} and
{\it corotational}, meaning that the stress rate vanishes for all
rigid body motions and commutes with index raising and index lowering,
respectively, as articulated in \cite{Prag61EDDS}.

For stress rates (in contrast to general strain rates),
the {\it Zaremba-Jaumann}'s \cite{Jaum11GSPC,Zare03FPTR}
{\it derivative} (which is  corotional and the simplest variant of
the objective derivative, as articulated in \cite{Haup02CMTM})
is considered as the proper one, as also derived by M.\,Biot
\cite[p.494]{Biot65MID} and advocated e.g.\ in
\cite{Bruh09EEBI,Fial11GSSM,Fial20OTDR}. It is widely used both
in theoretical and numerical modelling in particular in geophysics, see e.g.\
\cite{BabSob08HRNM,BuMaMi17SDRN,Gery19INGM,MoDuMu03LIPF,PGYD20SHMM}, although
some spurious ratchetting effects under cyclic loading may occur, as
pointed out in \cite{JiaFis17ADRD,MeXiBrMe03ESRC}, and although some other
corotational derivatives (as Green-Naghdi \cite{GreNag65GTEP} or
Xiao-Bruhns-Meyer \cite{XiBrMe00CORF,XiBrMe06EPSD}) are sometimes
considered, too. Cf.\ \cite{LiuHon99NCHM} for a survey and comparison.
Perhaps, in the context of the solid-liquid phase transformation, the simple
convincing argument is that, when the stress tensor is just
a pressure (i.e.\ in the form of the unit matrix multiplied by a scalar),
then the corotational derivatives (and in particular the Zaremba-Jaumann
derivative as used in \eqref{ED-2} below) degenerates to the convective
derivative, which is indeed expected in fluids,
cf.\ Remark~\ref{rem-semi-compress} below.
The other important attribute of this derivative is that symmetric (resp.\
trace-free) tensors stay symmetric (resp.\ trace-free) when transported
by this derivative.

We use the notation $\R_{\rm sym}^{d\times d}=\{A\in \R^{d\times d};\ A^\top=A\}$ and
$\R_{\rm dev}^{d\times d}=\{A\in \R_{\rm sym}^{d\times d};\ {\rm tr}\,A=0\}$
where ${\rm tr}\,A=\sum_{i=1}^dA_{ii}$ denotes the trace of the matrix $A$. We 
will use the decomposition of square matrices on symmetric and skew-symmetric parts:
\begin{align}
  {\rm sym}\,A:=\frac12A^\top\!+\frac12A\ \ \ \text{ and }\ \ \
  {\rm skew}\,A:=\frac12A^\top\!-\frac12A=A-{\rm sym}\,A\,,
\end{align}
We will also use another decomposition
to the spherical (volumetric) and the deviatoric parts, defined by
\begin{align}
  A={\rm sph}\,A+{\rm dev}\,A\ \ \ \ \text{ with }\ \ \
  {\rm sph}\,A:=\frac1d({\rm tr}\,A)\bbI\ \ \ \text{ and }\ \ \
  {\rm dev}\,A:=A-\frac1d({\rm tr}\,A)\bbI\,,
\end{align}
where ``tr'' denotes for a trace of a matrix. Note that it is an orthogonal
decomposition with respect to the Frobenius norm and that ${\rm dev}\,A$
is trace-free, i.e.\ ${\rm sph}\,A:{\rm dev}\,A=0$ and
${\rm tr}({\rm dev}\,A)=0$. In this notation,
${\rm sym}\,A\in\R_{\rm sym}^{d\times d}$ and
${\rm dev}({\rm sym}\,A)={\rm sym}({\rm dev}\,A)\in\R_{\rm dev}^{d\times d}$.

Also we will use the standard notation $AB$ for the product of
two matrices while $A{:}B$ will denote the scalar product of
these matrices. We recall the algebra
\begin{align}\label{algebra}
  A{:}(BC)=(B^\top\!A){:}C=(AC^\top){:}B
\end{align}
for any square matrices $A$, $B$, and $C$. Moreover, ``$\,\cdot\,$'' and
``$\,\Vdots\,$'' will be used for the scalar products of the vectors and for
the 3rd-order tensors, respectively.

We consider a fixed bounded domain $\varOmega\subset\R^d$ with
a Lipschitz boundary $\varGamma$ and a time interval $I=[0,T]$.
The basic variables and data for this preliminary model are summarized
in the following table:

\begin{center}
\fbox{
\  \begin{minipage}[t]{20em}

$\vv$ velocity

$\EE(\vv)={\rm sym}(\nabla\vv)$

 $\bm{T}$ the Cauchy stress  

  $\bm{S}$ the Piola-Kirchhoff stress tensor

 $\bm{D}$ the dissipative part of the Cauchy stress 

   $\varrho$ mass density  (considered constant) 

\end{minipage}
\begin{minipage}[t]{24em}
    $\varphi^*$ stored energy (in terms of the stress $\bm{S}$)

$\bbC$ elasticity tensor (hence $\bbC^{-1}$ the compliance tensor)

$\bbD$  Stokes viscosity tensor 

  $\ff$ bulk force

$\bm{\varSigma}$ the conservative part of the Cauchy stress 

$\bm g$ tangential traction force

\end{minipage}
} % end of fbox
\end{center}

  We will assume that $\bbC$ is fully symmetric in the sense
  $\bbC_{ijkl}=\bbC_{klij}=\bbC_{jikl}$ while $\bbD$ suffices to
  satisfy only a minor symmetry $\bbD_{ijkl}=\bbD_{jilk}$, both
  these tensors being positive definite, and $\bbC$ as well
  as the mass density $\varrho>0$ independent
  of $x\in\varOmega$. These assumptions are important to ensure
  energetic consistency of the following system considered
  on $I{\times}\varOmega$:
  \begin{subequations}\label{ED}\begin{align}\nonumber
&\varrho\DT\vv={\rm div}\,\bm{T}-\frac\varrho2({\rm div}\,\vv)\,\vv+\ff\,,
      \ \ \ \text{where }\ \ \
      \DT\vv=\pdt\vv+(\vv{\cdot}\nabla)\vv\ \ \text{ and }\ \ \bm{T}=
\bm{\varSigma}+\bm{D}\ 
\\[-.2em]
\label{ED-1}&\hspace{12.8em}
\text{ with }\ \ \ \bm{\varSigma}=
      \bm{S}+\varphi^*(\bm{S})\bbI\,,\ \ \ 
\varphi^*(\bm{S})=\frac12\bm{S}{:}\bbC^{-1}\bm{S}\ \ \text{ and }\ \
\bm{D}=\bbD\EE(\vv)\,,
\\[.1em]&
\ZJ{\bm{S}}=\bbC\EE(\vv)
\ \ \ \ \text{where}\ \ \ \ 
 \ZJ{\bm{S}}= \!\!\linesunder{\DT{\bm{S}}
   -\frac{\nabla\vv{-}\nabla\vv^\top\!\!}2{\bm{S}}
   -{\bm{S}}\frac{\nabla\vv^\top\!{-}\nabla\vv}2}{Zaremba-Jaumann}{
co-rotational derivative}
 \ \text{with }\ \ \DT{{\bm{S}}}=\!\!\!
 \linesunder{\pdt{\bm{S}}
   +(\vv{\cdot}\nabla){\bm{S}}}{convective}{derivative}\!\!\!.
\label{ED-2}\end{align}\end{subequations}
We can equally write $\ZJ{\bm{S}}=\DT{\bm{S}}
-{\rm skew}(\nabla\vv){\bm{S}}+{\bm{S}}{\rm skew}(\nabla\vv)$.
The penultimate  term $\frac\varrho2({\rm div}\,\vv)\,\vv$ 
in the first equation \eqref{ED-1} is a force introduced by Temam \cite{Tema69ASEN}
to balance the energy in flows of fluids which are so-called
quasi-incompressible, cf.\ also \cite{Toma21ITST} for a
certain justification. The  pressure contribution $\varphi^*(\bm{S})$
to the Cauchy stress $\bm{\varSigma}$ is due to the concept of considering the stored
energy per unit actual (not referential) volume and it matches the energy
balance, cf.\ \eqref{energy} below.
The (choice of the) Zaremba-Jaumann derivative in \eqref{ED-2} respects,
in particular, the proper evolution of the eigenvalues of $\bm{S}$,
as articulated in \cite[Thm.8.1]{Haup02CMTM}.

Of course, this system of nonlinear parabolic equations should be
accompanied by some boundary conditions, e.g.\
\begin{align}\nonumber\\[-2.2em]\label{ED-BC}
  \vv{\cdot}\nn=0\ \ \ \ \text{ and }\ \ \ \ \big[\bm{T}\nn\big]_\text{\sc t}^{}={\bm g}\,
\end{align}
with $\nn$ denoting the unit outward normal to $\varGamma$, $[\cdot]_\text{\sc t}^{}$
the tangential component of a vector, i.e.\
$[\bm\sigma]_\text{\sc t}^{}=\bm\sigma-(\bm\sigma{\cdot}\nn)\nn$
for a vector $\bm\sigma$. Naturally, ${\bm g}{\cdot}\nn=0$ is to be assumed.
Further, we prescribe the initial conditions 
\begin{align}\label{ED-IC}
\vv|_{t=0}^{}=\vv_0 \ \ \ \text{ and }\ \ \ \bm{S}|_{t=0}^{}=\bm{S}_0\,.
\end{align}

The energy balance behind this model is revealed by adding the momentum
equation \eqref{ED-1} tested by $\vv$ and the stress equation \eqref{ED-2}
tested by $\bbC^{-1}\bm{S}$. We use several calculations
exploiting integration over $\varOmega$ and Green's formula. In particular:

Exploiting the Green formula, the calculus for the convective term
\begin{align}\nonumber
  \INT{\varOmega}{}\varrho(\widetilde\vv{\cdot}\nabla)\vv\cdot \vv\,\d x
 & =\INT{\varGamma}{}\varrho|\vv|^2(\widetilde\vv{\cdot}\nn)\,\d S
  - \INT{\varOmega}{}\varrho\vv\,{\rm div}(\widetilde\vv{\otimes}\vv)\,\d x
  \\&\nonumber=\INT{\varGamma}{}\varrho|\vv|^2(\widetilde\vv{\cdot}\nn)\,\d S
  - \INT{\varOmega}{}\varrho|\vv|^2{\rm div}\,\widetilde\vv
  +\varrho\widetilde\vv{\cdot}\nabla\vv\cdot\vv
+|\vv|^2(\nabla\varrho{\cdot}\widetilde\vv)\,\d x
  \\&=\INT{\varGamma}{}\frac\varrho2|\vv|^2(\widetilde\vv{\cdot}\nn)\,\d S
  -\INT{\varOmega}{}\frac\varrho2|\vv|^2{\rm div}\,\widetilde\vv
  +\frac12|\vv|^2(\nabla\varrho{\cdot}\widetilde\vv)\,\d x
  \label{convective-tested}\end{align}
is used for $\widetilde\vv=\vv$ and for $\varrho$ constant so that
$\nabla\varrho=0$ and the last term in \eqref{convective-tested} vanishes;
here the role of the bulk force $\frac12\varrho({\rm div}\,\vv)\vv$ in
\eqref{ED-1} is revealed.

Further, by the assumed symmetry of $\bbC$ and $\bbD$,
\begin{align}
  \!\!\INT{\varOmega}{}{\rm div}\,(\bm{\varSigma}{+}\bbD\EE(\vv))
          {\cdot}\vv\,\d x=
          \!\INT{\varGamma}{}(
          \bm{\varSigma}{+}\bbD\EE(\vv)){:}(\vv\otimes\nn)\,\d S-
          \!\!\INT{\varOmega}{}(
          \bm{\varSigma}{+}\bbD\EE(\vv)){:}\EE(\vv)\,\d x.
\label{calculus-3-}\end{align}
For the conservative-stress terms, we futher use
\begin{align}\nonumber
  &\INT{\varOmega}{}\bm{\varSigma}{:}\EE(\vv)\,\d x
=\INT{\varOmega}{}\bm{S}{:}\EE(\vv)+\varphi^*(\bm{S}){\cdot}{\rm div}\,\vv\,\d x
=
\INT{\varOmega}{}\bm{S}{:}\bbC^{-1}\ZJ{\bm{S}}+\varphi^*(\bm{S}){\cdot}{\rm div}\,\vv\,\d x
\\\nonumber&\qquad=\frac{\d}{\d t}
\INT{\varOmega}{}\varphi^*(\bm{S})\,\d x
+\INT{\varOmega}{}\varphi^*(\bm{S}){\cdot}{\rm div}\,\vv
+\bbC^{-1}\bm{S}{:}\big((\vv{\cdot}\nabla){\bm{S}}
({\rm skew}\nabla\vv)
-{\rm skew}(\nabla\bm v)\bm{S}-{\bm{S}}({\rm skew}\nabla\vv)^\top\big)\,\d x
\\
&\qquad=\frac{\d}{\d t}
\INT{\varOmega}{}\varphi^*(\bm{S})\,\d x
+\INT{\varOmega}{}{\rm skew}\big({\bm{S}}(\bbC^{-1}\bm{S})^{\top}
  \!-(\bbC^{-1}\bm{S}){\bm{S}}^{\top}\big){:}\nabla\vv\,\d x=\frac{\d}{\d t}
\INT{\varOmega}{}\varphi^*(\bm{S})\,\d x\,,
\label{calculus-3}\end{align}
where we used also \eqref{ED-2}.  In particular, relying on that $\bbC$ is
constant in space, we used also the Green formula for
\begin{align}
  \INT{\varOmega}{}\varphi^*(\bm{S}){\cdot}{\rm div}\,\vv+
\bbC^{-1}\bm{S}{:}(\vv{\cdot}\nabla){\bm{S}}
\,\d x
=\INT{\varOmega}{}\varphi^*(\bm{S}){\cdot}{\rm div}\,\vv+
\nabla\varphi^*(\bm{S}){\cdot}\vv\,\d x
=\INT{\varGamma}{}\varphi^*(\bm{S})(\vv{\cdot}\nn)\,\d S=0
\label{calculus-4}\end{align}
and the matrix algebra \eqref{algebra} for 
\begin{align}\nonumber
  &\Ee{:}\big(({\rm skew}\nabla\vv)\bm{S}+\bm{S}({\rm skew}\nabla\vv)^\top\big)
  =\frac12\Ee{:}\big((\nabla\vv)\bm{S}-(\nabla\vv)^\top\bm{S}
-\bm{S}(\nabla\vv)+\bm{S}(\nabla\vv)^\top\big)
\\&\qquad\qquad\nonumber
=\frac12\big(\Ee\bm{S}^\top\!-\bm{S}^\top\Ee\big){:}\nabla\vv
-\frac12\bm{S}{:}\big((\nabla\vv)\Ee-\Ee(\nabla\vv)^\top\big)
\\&\qquad\qquad
=\frac12\big(\Ee\bm{S}^\top\!-\bm{S}\Ee^\top\!-\bm{S}^\top\Ee+\Ee^\top\bm{S}\big){:}\nabla\vv
={\rm skew}(\hspace{-.7em}\lineunder{\Ee\bm{S}^\top\!-\bm{S}\Ee^\top\!\!}{$={\bm0}$}\hspace{-.7em}){:}\nabla\vv
={\bm 0}
\label{algebra+}\end{align}
with $\Ee=\bbC^{-1}\bm{S}$;
here we have to assume that  $\bm{S}$ commutes with $\Ee$ and the initial
condition for $\bm{S}$ is valued in $\R_{\rm sym}^{d\times d}$ and exploit that
the Zaremba-Jaumann corotational derivative in \eqref{ED-2} keeps symmetry
of $\bm{S}$ during the whole evolution.

Then we should sum up \eqref{calculus-3-} and \eqref{calculus-3} and
use the boundary conditions \eqref{ED-BC}, so that
$\int_{\varGamma}{}\bm{T}{:}(\vv\otimes\nn)\,\d S
=\int_{\varGamma}{}{\bm g}\cdot\vv\,\d S$ and thus
\begin{align}
  \INT{\varOmega}{}{\rm div}\,\bm{T}{\cdot}\vv\,\d x&=
\INT{\varGamma}{}{\bm g}\cdot\vv\,\d S
-\INT{\varOmega}{}\bm{T}{:}\nabla\vv\,\d x
=\INT{\varGamma}{}{\bm g}\cdot\vv\,\d S
-\INT{\varOmega}{}
\Big(\varrho\DT\vv+\frac\varrho2({\rm div}\,\vv)\vv-\ff\Big)\cdot\vv\,\d x\,.
\label{calculus-6}\end{align}
Thus we eventually obtain the energy balance
\begin{align}
\frac{\d}{\d t}
\INT{\varOmega}{}\!\linesunder{\frac\varrho2|\vv|^2}{kinetic}{$^{^{}}$energy}\hspace{-.7em}+\hspace{-.7em}\linesunder{\frac12\bbC^{-1}\bm{S}{:}\bm{S}}{stored}{$^{^{}}$energy}\hspace{-.7em}\,\d x+
\INT{\varOmega}{}\!\!\lineunder{\bbD\EE(\vv){:}\EE(\vv)_{_{_{_{_{_{_{}}}}}}}\!\!}{dissipation rate}\hspace{-.7em}\d x
=\INT{\varOmega}{}\hspace{-1em}\linesunder{\,\ff\cdot\vv\,}{power of}{external load}\hspace{-1em}\d x
+\INT{\varGamma}{}\hspace{-1em}\linesunder{\,{\bm g}\cdot\vv\,}{power of}{traction load}\hspace{-1em}\d S
\,.
\label{energy}\end{align}

  \begin{remark}[{\sl Isotropic materials}]\label{rem-isotropic}\upshape
    The prominent application of this convective model is for isotropic
    materials which are invariant for rotations. The elasticity tensor
    $\bbC$ is then of the form
  \begin{align}\nonumber\\[-2.3em]\label{C}
    \bbC_{ijkl}:=K_\text{\sc e}^{}\delta_{ij}\delta_{kl}
    +G_\text{\sc e}^{}\Big(\delta_{ik}\delta_{jl}
    +\delta_{il}\delta_{jk}-\frac2d\delta_{ij}\delta_{kl}\Big)\,,
  \end{align}
  where $K_\text{\sc e}^{}$ is the bulk elastic modulus and
  $G_\text{\sc e}^{}$ is the shear elastic modulus with the physical
  dimension Pa=J/m$^3$, while $\delta$ is the Kronecker symbol.
  Actually, $G_\text{\sc e}^{}$ is also called the second Lam\'e parameter
  (denoted as $\mu$) while the first Lam\'e parameter (denoted as $\lambda$)
  is $K_\text{\sc e}^{}-2G_\text{\sc e}^{}/d$, so that $K_\text{\sc e}^{}=\lambda
  +2\mu/d$. The viscosity tensor $\bbD$ for isotropic materials is then in
  an analogical form
  \begin{align}
    \label{D}
    \bbD_{ijkl}:=K_\text{\sc v}^{}\delta_{ij}\delta_{kl}
    +G_\text{\sc v}^{}\Big(\delta_{ik}\delta_{jl}
    +\delta_{il}\delta_{jk}-\frac2d\delta_{ij}\delta_{kl}\Big)\,
 \end{align}
 with $K_\text{\sc v}^{}$ and $G_\text{\sc v}^{}$ the viscosity
 moduli with the physical dimension Pa\,s. In this  isotropic case, 
  the stored energy in \eqref{energy} reads as
\begin{align}\nonumber\\[-2.em]\label{phi+}
  \varphi^*(\bm{S})=\frac12\bbC^{-1}\bm{S}{:}\bm{S}
  =\frac1{2dK_\text{\sc e}^{}}|{\rm sph}\,\bm{S}|^2
  +\frac1{4G_\text{\sc e}^{}}|{\rm dev}\,\bm{S}|^2=
\frac1{2d^2K_\text{\sc e}^{}}|{\rm tr}\,\bm{S}|^2
  +\frac1{4G_\text{\sc e}^{}}|{\rm dev}\,\bm{S}|^2\,,
\end{align}
cf.\ \cite[Formulas (6.7.9)--(6.7.10)]{KruRou19MMCM}. In terms of a strain
tensor $\Ee$, the stress $\bbC\Ee$ in this isotropic case looks as
$dK_\text{\sc e}^{}{\rm sph}\,\Ee+2G_\text{\sc e}^{}{\rm dev}\,\Ee$.
It is important that, taking into account \eqref{C} and the symmetry of $\Ee$,
$\bm{S}=\bbC\Ee$ now means
$\bm{S}_{ij}=K_\text{\sc e}^{}\delta_{ij}\Ee_{kk}+2G_\text{\sc e}^{}(\Ee_{ij}-
\delta_{ij}\Ee_{kk}/d)$ with the summation convention and also
$\bm{S}\Ee=\Ee\bm{S}$. Then it is
easy to check that $\nabla\vv\bbC\Ee=\bbC(\nabla\vv)\Ee$ and also
$(\nabla\vv)^\top\bbC\Ee=\bbC(\nabla\vv)^\top\Ee$.
Altogether, with some $\Ee$ symmetric tensor and $\bbC$ homogeneous
(i.e.\ independent of $x$) isotropic in the sense \eqref{C}, we can see that 
the strain $\Ee$ is also transported by the Zaremba-Jaumann derivative, i.e.
\begin{align}\nonumber\\[-2.1em]\label{commute}
  \ZJ{\bm{S}}=\bbC\ZJ{\Ee}\,.
\end{align}
This will justify the strain/velocity formulation used later in
Sections~\ref{sec-fuller}--\ref{sec-anal}  where the convex-conjugate
energy, i.e.\ here $\varphi(\Ee)=[\varphi^*]^*(\Ee)=
\frac d2K_\text{\sc e}^{}|{\rm sph}\,\Ee|^2+G_\text{\sc e}^{}|{\rm dev}\,\Ee|^2$,
will be used, cf.\ \eqref{psi} below. 
  \end{remark}
  
\begin{remark}[{\sl Semi-compressible fluids}]\label{rem-semi-compress}\upshape
If the elastic shear modulus $G_\text{\sc e}^{}$ in \eqref{C} vanishes
(i.e.\ the elastic response in the solid model on the deviatoric stress
vanishes), the medium becomes fluidic. Yet, the elastic response
on the spherical stress remains and thus the model describes 
so-called viscoelastic fluids. Notably, such fluids facilitate propagation
of pressure waves, in contrast to the merely viscous fluids.
Then $\bbC\EE(\vv)=dK_\text{\sc e}^{}{\rm sph}\,\EE(\vv)=
K_\text{\sc e}^{}({\rm div}\,\vv)\bbI$. 
From the last term in \eqref{phi+}, we can see that ${\rm dev}\,\bm{S}$ must
vanish. Then $\bm{S}(t,x)=-p(t,x)\bbI$ for some variable $p$ meaning a pressure
\eqref{ED-2} turns into
\begin{subequations}\label{semi-compres}\begin{align}\label{semi-compres-1}
 \DT{p}+K_\text{\sc e}^{}{\rm div}\,\vv=0\,.
\end{align}
Here it is important that, for $\bm{S}=-p\bbI$, it holds
$\ZJ{\bm{S}}=-\DT{p}\,\bbI$.
The momentum equation \eqref{ED-1} degenerates to 
\begin{align}
  \label{semi-compres-2}
  \varrho\DT\vv=
{\rm div}\big(\bbD\EE(\vv)\big)-\nabla\Big(p+\frac{p^2}{2K_\text{\sc e}^{}}\Big)
  -\frac\varrho2({\rm div}\,\vv)\,\vv+\ff\,.
\end{align}\end{subequations}
It allows for propagation of
{\it longitudinal} (also called pressure) {\it waves} (so-called P-waves).
In a one-dimensional situation, one can see that 
longitudinal waves propagate with the velocity
$v=v(\lambda)=\sqrt{K_\text{\sc e}^{}/\varrho-K_\text{\sc v}^2/(4\varrho^2\lambda^2)}$
with $\lambda$ denoting the wave length and $K_\text{\sc v}^{}$ is from \eqref{D}.
This shows that, due to the viscosity, this model exhibits a so-called
{\it normal dispersion}, i.e.\ the speed of waves increases
with their wavelength, cf.\ \cite{Roub??QISC}.
\end{remark}

\begin{remark}[{\sl Gradient theories towards analysis}]\label{rem-gradients}\upshape
  The existence of weak solutions of the models \eqref{ED} or
  \eqref{semi-compres} is unfortunately not ensured. To facilitate rigorous
  analysis, it is seems inevitable to involve some gradient theories into the
  model. There are several options.  The first option it to augment the energy $\varphi^*$ by
  $\frac\eps2\nabla{\bm S}{\Vdots}\nabla\bbC^{-1}{\bm S}$ and then the Cauchy stress 
  stress in \eqref{ED-1} by
  $-\eps\Delta{\bm S}+\eps\nabla{\bm S}{\otimes}\bbC^{-1}\nabla{\bm S}$ 
  with a small $\eps>0$, which leads to a conservative expansion of the
  energetics \eqref{energy}. The second option is to augment \eqref{ED-2} by
  $-\eps\Delta{\bm S}$, which leads to a dissipative expansion of the
energetics \eqref{energy} by $\eps\nabla{\bm S}\Vdots\nabla\bbC^{-1}{\bm S}$.
  This is a concept of so-called (small) stress diffusion, which 
  is related with a (large) P\'eclet (in fluids also called Brenner) number.
  The third option is to augment the stress in \eqref{ED-1} by
  $-\nu\Delta\EE(\vv)$, which leads to a dissipative expansion of the
  energetics \eqref{energy} by $\nu|\nabla\EE(\vv)|^2$.
This is a concept of so-called 2nd-grade nonsimple materials.
  In what follows, we will employ a nonlinear variant of this third option.
  These three options can be reflected also in the semi-compressible fluid
  model \eqref{semi-compres} and have been addressed in \cite{Roub??QISC}.
  The first and the third options leads to anomalous dispersion while
  the second one contributes to the normal dispersion of the velocity
  of the P-waves. For the 
  limit for $K_\text{\sc e}^{}\to\infty$ (possibly together with $\eps\to0$ for
  the second option), one can prove that the semi-compressible fluid
turns into the incompressible Navier-Stokes system, which however does not
allow for propagation of P-waves.
\end{remark}

\begin{remark}[{\sl Varying mass density}]\upshape
  The above model  has  used the simplification based on the
  assumption of a constant mass density $\varrho$  even though the model
  is not fully incompressible, cf.\ the calculus \eqref{convective-tested}.
   This is legitimate when, in particular,  variations of $\varrho$
  during volumetric deformation (which is typically indeed small in
  liquids and solids, in contrast to gases)  can be neglected;
  e.g.\ even extremely compressed water in the deepest spot in Earth oceans
  (i.e.\ Mariana trench) shrinks
  only by 5\% in volume and rocks varies their volume less than 0.1\% by longitudianal
  waves emitted during even big earthquakes. (Remarkable density variations can
  occur rather during phase transformations but we neglect it, too.)
   To keep energy balance  while  omitting the continuity equation
   $\DT\varrho=-\varrho{\rm div}\,\vv$
   is then compensated by the force $-\varrho({\rm div}\,\vv)\vv/2$
   in \eqref{ED-1} and in \eqref{ED-1+}, which is pressumably very small
   but not Galilean invariant, as pointed out in \cite{Toma21ITST}. This modeling
   shortcut simplifies a lot of calculations  and allows us to focus on the
   thermomechanical context of the Stefan problem itself. If the medium would
   be semi-compressible but inhomogeneous as far as mass density, one
   should complete the system by a simplified continuity equation
   $\DT\varrho=0$ while still keep the compensating force in the momentum
   equation. Like for the fully compressible model without the compensating
   force but with the full continuity equation 
   $\DT\varrho=-\varrho{\rm div}\,\vv$, e.g.\
   \eqref{test-of-convective}--\eqref{test-of-convective+} and 
   \eqref{strong-hyper} would become considerably more technical
   although, because the concept of nonlinear nonsimple material ensures
    regular velocity fields $\nabla\vv$
   in $L_{\rm w*}^p(I;L^\infty(\varOmega;\R^{d\times d}))$, the analysis
   would be doable. 
 \end{remark}

\def\ZJEp{\Ep}

\section{A thermomechanical model with  inelasticity and damage}\label{sec-fuller}
%        ~~~~~~~~~~~~~~~~~~~~~~~~~~~~~~~~~~~~~~~~~~~~~~~~~~~~~~

Now we will enhance the model 
by some internal variables, namely a scalar-valued {\it damage} $\alpha$ and
an {\it inelastic} (or plastic) {\it strain} $\bm{P}$ which is
considered standardly symmetric and isochoric, i.e.\
${\rm tr}\,\bm{P}=0$.
Moreover, we enhance it towards full thermodynamics by
considering a heat-transfer equation for {\it temperature} $\theta$
coupled with the mechanical part in a thermodynamically consistent way.

Involving damage and creep (or plasticity) into the stress/velocity
formulation is quite tricky. Actually, the creep model (arising from a
dissipation potential quadratic in terms of the rate $\ZJEp$ of
$\bm{P}$) can be understood also as the Jeffreys' rheological model.

However, in the convective damage model, this does not seem working properly
because the constancy of $\bbC$ used in \eqref{calculus-4} cannot hold
any more because $\bbC$ must now depend on damage, otherwise there would not
be any driving force for damage processes.

In order to eliminate the
tensor $\bbC$ from the formulation of the model, we switch into a
{\it strain/velocity formulation}. We confine ourselves on the {\it additive},
also called {\it Green-Naghdi's decomposition} which
 reads at small strains as   $\EE(\uu)=\Ee+\bm{P}$,
where $\Ee$ denotes the elastic strain and $\uu$ the displacement.
Yet, as we do not want to stay on small displacements, we formulate
the additive decomposition rather in terms of rates as
\begin{align}\nonumber\\[-3.2em]\label{additive-rate}
  \EE(\vv)=\ZJ{\overline{\Ee{+}\bm{P}}}=
\ZJ\Ee+\ZJEp\ \ \ \ \text{ if denoting}\ \ \ZJEp=\ZJ{\bm{P}}\,.
\end{align}
The new variable $\ZJEp$ is a rate of plastic strain. 
Thus both $\uu$ and $\bm{P}$ will not be eventually used, and the model
allows for {\it large displacements} even though still using
linearized elastic strains. For example, after solid has being melted, the
resulting liquid may flow on large distances where it can again solidify
etc. Usage of the corotational derivative in \eqref{additive-rate}
is consistent with Section~\ref{sec-preliminary} where, instead of
\eqref{ED-2}, we may use $\ZJ{\Ee}=\EE(\vv)$ provided the material is 
isotropic homogeneous, cf.\ \eqref{commute}. 

As we want to model also solid-liquid phase transition like the
classical Stefan problem in the so-called enthalpy formulation, the
deviatoric elastic response will be influenced by temperature and vanishes
completely above a melting temperature when the material becomes
viscoelastic fluid, cf.\ Remark~\ref{rem-semi-compress}.
However, direct involvement of the temperature in the elastic shear 
modulus would have to be nonlinear and would cause complicated
adiabatic effects. We thus devise the model differently by allowing
temperature to influence only the dissipative part, cf.\ Figure~\ref{fig1}
below.

The additional variables and data employed in the enhancement
and the strain/velocity (instead of stress/velocity) formulation 
of mechanical merely visco-elastodynamic system \eqref{ED}
are summarized in the following table:
\begin{center}
\fbox{
 \begin{minipage}[t]{19em}

$\bm E$ the elastic strain tensor (symmetric)

$\ZJEp$ the inelastic strain rate tensor 

$\alpha$ damage variable

$\theta$ temperature  (in K)

$\W$ enthalpy (internal heat energy in J/m$^3$)

$\eta$ entropy (in J/m$^3K$)

$\chi$ water/ice volume fraction

 ${\bm T}$ total Cauchy stress (symmetric) 

 ${\bm S}$ Piola-Kirchhoff stress 

${\bm K}$ Korteweg-like stress 

\end{minipage}
\begin{minipage}[t]{25.5em}

\vspace*{-.5em}

${\bm D}$ dissipative stress 
  
$\mbox{$\bbK$}=\bbK(\alpha,\W)$ thermal conductivity 
  
$\zeta=\zeta(\alpha,\W;\cdot)$ dissipation potential for damage rate

   $\psi=\psi(\Ee,\alpha,\theta)$ free energy (in Pa=J/m$^3$) 

   $\varphi=\varphi(\Ee,\alpha)$ the stored energy (in Pa=J/m$^3$) 

  $\phi=\phi(\theta)$ thermal part of the free energy (in Pa=J/m$^3$)

  $\beta=\gamma^{-1}$ with $\gamma(\theta)=\phi(\theta)-\theta\phi'(\theta)$

  $\GM\!=\GM(\W)$ Maxwell-viscosity modulus (in Pa\,s=Nm$^{-2}$s)

 $\kappa$, $\varkappa$ length-scale coefficients 

  $\nu$ a hyper-viscosity coefficient.

\end{minipage}
}  % end of fbox
\end{center}

The main ingredient will be the free energy $\psi$, which expands the
stored energy in terms of strains (i.e.\ the dual to the
stored energy $\varphi^*$ in terms of stresses as used in
Section~\ref{sec-preliminary}). We will consider a gradient theory
for damage and, for simplicity, and additive split of mechanical
and thermal parts, namely 
\begin{align}\label{psi-}
\psi(\bm E,\alpha,\nabla\alpha,\theta)&=\varphi(\bm E,
\alpha)+\frac{\kappa}2|\nabla\alpha|^2+\phi(\theta)\,.
\end{align}
We may consider a general stored energy $\varphi=\varphi(\bm E,\alpha)$
possibly non-quadratic in terms of $\Ee$ while, as usual, we are using
a linear gradient theory for the damage variable $\alpha$.

Beside the free energy, another important ingredient is the (pseudo)potential
of dissipative forces acting on rates of internal variables, i.e.\ 
here of the inelastic strain and the damage. As to the inelastic strain rate
$\ZJEp$, we consider a materially linear Maxwellian creep together with
a gradient theory, so that the corresponding
potential is quadratic in terms of rate,
i.e.\ $\frac12\GM|\ZJEp|^2+\frac12\varkappa|\nabla\ZJEp|^2$
for some (presumably small) coefficient $\varkappa>0$;
in what follows, $\GM$ will depend on temperature (or rather on enthalpy
$\W$).
A nonlinear generalization will be outlined in Remark~\ref{rem-nonlin-creep}
below. On the other hand, evolution of damage is surely always a very
nonlinear phenomenon, in particular that damaging of solid is a fast process
while healing is rather slow (except when damaged solid starts melting
and healing is to be fast so that $\alpha$ goes fast to 1 in created liquid
phase and possible later solidification leads again to an undamaged solid).
Thus for the damage rate $\DT{\alpha}$, we will consider a general convex
potential depending also on state variables $(\alpha,\W)$, let us
denote it by $\zeta=\zeta(\alpha,\theta;\cdot)$. 

We used the gradient theory for the rate $\ZJEp$, as suggested in
\cite{DaRoSt??NHFV},  which does not cause any  spurious hardening-like
effects during large slips particularly along cracks or fluid-solid interface.
Simultaneously, as outlined already in
Remark~\ref{rem-gradients}, we will use a higher-order velocity gradient.
Here, the goal is to ensure
$\nabla\vv\in L_{\rm w*}^1(I;L^\infty(\varOmega;\R^{d\times d}))$,
which further facilitates regularity of the transport, cf.\ 
in particular the estimates 
\eqref{Gronwall} and \eqref{transport+++} below. Therefore, another
enhancement will be by dissipative gradient terms, exploiting the concept of
the so-called {\it nonsimple fluids}, devised by E.\,Fried and M.\,Gurtin
\cite{FriGur06TBBC} and earlier, even more generally and nonlinearly as
multipolar fluids, by J.\,Ne\v cas at al.\
\cite{BeBlNe92PBMV,NeNoSi89GSIC,NecRuz92GSIV}. More specifically, we use 
nonlinear 2nd-grade nonsimple fluids, also called {\it bipolar fluids},
monolithically also in the solid part although making the coefficient
$\nu$ in \eqref{ED-1+} dependent on enthalpy could distinguish
``intensity'' of this concept in the liquid and in the solid areas.

The above introduced model now reads specifically as
\begin{subequations}\label{ED+}\begin{align}
 \nonumber &\varrho\DT\vv={\rm div}\,{\bm T}
-\frac\varrho2({\rm div}\,\vv)\,\vv+\ff\,
\ \ \ \text{where }\ \ \
           {\bm T}=\bm{\varSigma}+{\bm K}+\bm{D}
\\\nonumber
&\hspace*{3em}
\ \ \ \ \ \text{ with }\
\bm{\varSigma}=\varphi_\Ee'(\Ee,\alpha)+\psi(\Ee,\alpha,\theta)\bbI\,,
\ \ \ {\bm K}=
\kappa\nabla\alpha{\otimes}\nabla\alpha-\frac{\kappa}2|\nabla\alpha|^2\bbI\,,
\\&\hspace*{3em}
\ \ \ \ \ \text{ and }\ \ 
\bm{D}=\bbD\EE(\vv)
-{\rm div}\big(\nu|\nabla\EE(\vv)|^{p-2}\nabla\EE(\vv)\big)\,,\ 
\label{ED-1+}
\\[-.3em]&
\ZJ\Ee=\EE(\vv)-\ZJEp\,,
\label{ED-2+}
\intertext{and  the equation for creep rate and the flow-rules for
 damage}
 \label{ED-3+}
 &\GM(\W)\ZJEp={\rm dev}\,\bm{\varSigma}+\varkappa\Delta\ZJEp\,,
        \\\label{ED-4+}
        &\partial_{\DT\alpha}\zeta\big(\alpha,\W;
   \DT{\alpha}\big)+\varphi_\alpha'(\Ee,\alpha)\ni \kappa\Delta\alpha\,.
   \intertext{Moreover, the system should be thermodynamically closed by the
heat-transfer equation, cf.\ \eqref{heat-eq} below, as}
\nonumber
&\pdt\W+{\rm div}\big(\vv\,\W-\bbK(\alpha,\W)\nabla\theta\big)
=\GM(\W)|\ZJEp|^2\!
+\xi\big(\alpha,\W;\DT\alpha\big)
+\bbD\EE(\vv){:}\EE(\vv)\!
\\&\label{ED-6+thermo}
\hspace{2.3em}+\nu|\nabla\EE(\vv)|^p\!+
\varkappa|\nabla\ZJEp|^2+
\phi(\theta){\rm div}\,\vv
\ \ \ \ \ \text{ with }\ \ \
\xi\big(\alpha,\W;\DT\alpha\big)
=\partial_{\DT\alpha}\zeta\big(\alpha,\W;\DT\alpha\big)\DT\alpha
\\[.2em]&\hspace*{20em}
\text{and\ \ with }\ \theta=\beta(\W)\,,
\label{beta}\end{align}\end{subequations}
where  $\beta$ is related to $\phi$ by 
$\beta^{-1}(\theta)=\phi(\theta)-\theta\phi'(\theta)$.

Note that, in \eqref{ED-3+}, ${\rm dev}\bm{\varSigma}={\rm dev}\bm{S}$
with the Piola-Kirchhoff stress $\bm{S}=\varphi_\Ee'(\Ee,\alpha)$
independent of temperature.
Also note that \eqref{ED-3+} is a quasistatic equation from which,
when prescribing the boundary conditions \eqref{BC3}, $\ZJEp$ can be explicitly
isolated as a function of $\W$ and $\bm{S}$. Moreover, $\ZJEp$
is always trace-free and, if $\bm{S}$ is symmetric,
$\ZJEp$ is also symmetric; this last attribute is indeed granted  
if the initial elastic strain is symmetric, cf.\ also
Remark~\ref{rem-plast-strain}. 

We complete the system \eqref{ED+} by suitable boundary conditions.
This is however little delicate due to the
2nd-grade nonsimple fluid concept and naturally needs 
2 boundary conditions for the momentum equation \eqref{ED-1+}. 
One of variationally natural option is to prescribe
boundary traction while putting higher derivatives zero, i.e.\
Neumann-type boundary condition. Here we take it in the tangential
direction while combining it with Dirichlet/Neumann condition in  
the normal direction. As it concerns $d$-dimensional vector field $\vv$,
we actually prescribe $2d$ conditions on the boundary. More specifically,
for the heat flux $h_{\rm ext}$, we consider
\begin{subequations}\label{BC}
\begin{align}\label{BC1}
  &
  \vv{\cdot}\nn=0\,,\ \ \ \ \big[{\bm T}
  \nn{+}\divS\big(\nu|\nabla\EE(\vv)|^{p-2}\nabla\EE(\vv){\cdot}\nn\big)
  \big]_\text{\sc t}^{}={\bm0}\,,
 \ \ \ \ \nabla\EE(\vv){:}(\nn{\otimes}\nn)={\bm0}\,,
\\\label{BC3}
  &(\nn{\cdot}\nabla)\ZJEp=0,
\ \ \ \ \ \ \ 
\nabla\alpha{\cdot}\nn=0,\ \ \ \ \ \ \
\nn{\cdot}\bbK(\alpha,\W)\nabla\theta=h_{\rm ext}\,,
\end{align}\end{subequations}
where $\divS={\rm tr}(\nablaS)$ with ${\rm tr}(\cdot)$ being the trace of a
$(d{-}1){\times}(d{-}1)$-matrix, denotes the $(d{-}1)$-dimensional
surface divergence and $\nablaS v=\nabla v-(\nabla v{\cdot}\nn)\nn$ 
being the surface gradient of $v$.
In contrast to \eqref{ED-BC}, we consider homogeneous mechanical condition
rather from the analytical reasons, cf.\ Remark~\ref{rem-g}.

We will consider an initial-value problem for the evolution boundary-value
problem \eqref{ED+}--\eqref{BC} and prescribe the corresponding initial
conditions enhancing \eqref{ED-IC}, namely
\begin{align}\label{IC}
  \vv|_{t=0}^{}=\vv_0,\ \ \ \ \ \bm{S}|_{t=0}^{}=\bm{S}_0,\ \ \ \ \ 
\alpha|_{t=0}^{}=\alpha_0,\ \ \text{ and }\ \ \W|_{t=0}^{}=\W_0\,.
\end{align}
Let us note that we thus prescribe also the initial
condition for $\bm E$, namely
\begin{align}
\Ee|_{t=0}^{}={\Ee}_0=\big[\varphi_\Ee'(\cdot,\alpha_0)\big]^{-1}(\bm{S}_0)\,.
\end{align}

The mechanical-energy balance behind the model (\ref{ED+}a-d) is revealed
by adding the momentum equation \eqref{ED-1+} tested by $\vv$ as before
in Section~\ref{sec-preliminary} and now the strain equation \eqref{ED-2+}
tested by the Piola-Kirchhoff stress $\bm{S}=\varphi_\Ee'(\Ee,\alpha)$
and further \eqref{ED-3+} tested by $\ZJEp$ and
\eqref{ED-4+} tested by $\DT\alpha$.
At least formally, we use the following calculations:

From \eqref{ED-1+} tested by $\vv$ and using again
\eqref{convective-tested}--\eqref{calculus-3}, we obtain
\begin{align}
  &
\frac{\d}{\d t}\INT{\varOmega}{}\,\frac\varrho2|\vv|^2\,\d x
+\!\INT{\varOmega}{}\big(\bm{S}{+}\varphi(\Ee,\alpha)\bbI{+}\bm K
{+}\bbD\EE(\vv)\big){:}\EE(\vv)+\nu|\nabla\EE(\vv)|^p
  +\phi(\theta){\rm div}\vv\,\d x=\!\INT{\varOmega}{}\ff{\cdot}\vv\,\d x\,,
\label{test-velocity}\end{align}
where we used that $\bm{S}=\varphi_\Ee'(\bm E,\alpha)$ is symmetric
and also the boundary conditions \eqref{BC} have been used;
actually, treatment of the boundary conditions for 2nd-grade nonsimple material
is indeed nonsimple and, beside usage of Green formula twice, needs also
usage of a surface Green formula, cf.\ \cite[Sect.\,2.4.4]{Roub13NPDE}. 
The conservative part of the Cauchy stress is to be treated
by using \eqref{ED-2+} tested by $\bm{S}$
and \eqref{ED-3+}  tested by $\ZJEp$.
Similarly as in \eqref{calculus-3} we obtain
\begin{align}\nonumber
  &
\INT{\varOmega}{}\big(\bm{S}{+}\varphi(\Ee,\alpha)\bbI{+}\bm K\big){:}\EE(\vv)\,\d x
  =\INT{\varOmega}{}
  \bm{S}{:}\ZJ\Ee{+}{\rm dev}\bm{S}{:}\ZJEp{+}\bm K{:}\EE(\vv)
   {+}\varphi(\Ee,\alpha){\rm div}\vv\,\d x
\\[-.2em]&=\INT{\varOmega}{}\varphi_\Ee'(\Ee,\alpha){:}\Big(\pdt{\Ee}
+(\vv{\cdot}\nabla)\Ee\Big)+\GM(\W)|\ZJEp|^2+\varkappa|\nabla\ZJEp|^2+\bm K{:}\EE(\vv)
  +\varphi(\Ee,\alpha)\,{\rm div}\vv\,\d x\,,
 \label{test-1}\end{align}
where we used also
$\bm{S}{:}({\rm skew}(\nabla\vv)\Ee-\Ee{\rm skew}(\nabla\vv))=0$ similarly
as \eqref{algebra+}, exploiting again the symmetry of $\Ee$ which is now
based on the symmetry not only of $\EE(\vv)$ but also of $\ZJEp$.

Furthermore, from \eqref{ED-4+} tested by $\DT\alpha$, we obtain
\begin{align}\nonumber
  \INT{\varOmega}{}
  \partial_{\DT\alpha}&\zeta\big(\alpha,\W;
  \DT{\alpha}\big)\DT{\alpha}\,\d x
=\INT{\varOmega}{}
\big(\kappa\Delta\alpha-\varphi_\alpha'(\Ee,\alpha)\big)\DT{\alpha}\,\d x
\\&\nonumber
=
\INT{\varGamma}{} \kappa(\nn{\cdot}\nabla\alpha)\DT{\alpha}\,\d S
-\frac{\d}{\d t}\INT{\varOmega}{}\frac{\kappa}2|\nabla\alpha|^2\,\d x
-\INT{\varOmega}{}\nabla\alpha{\cdot}\nabla(\vv{\cdot}\nabla\alpha)
+\varphi_\alpha'(\Ee,\alpha)\Big(\frac{\partial\alpha}{\partial t}
+\vv{\cdot}\nabla\alpha\Big)\,\d x
\\&\nonumber
=\INT{\varGamma}{} \kappa(\nn{\cdot}\nabla\alpha)
\DT\alpha
-\frac{\kappa}2|\nabla\alpha|^2(\vv{\cdot}\nn)\,\d S
-\frac{\d}{\d t}\INT{\varOmega}{}\frac{\kappa}2|\nabla\alpha|^2\,\d x
\\[-.3em]&\qquad\qquad\qquad\quad\
-\INT{\varOmega}{}\hspace{-.5em}
\lineunder{\kappa(\nabla\alpha{\otimes}\nabla\alpha){:}\EE(\vv)-
\frac{\kappa}2|\nabla\alpha|^2{\rm div}\,\vv}{$=\bm{K}{:}\EE(\vv)$}\hspace{-.9em}
+\varphi_\alpha'(\Ee,\alpha)\Big(\frac{\partial\alpha}{\partial t}
+\vv{\cdot}\nabla\alpha\Big)\,\d x\,,
\label{test-damage}\end{align}
where the boundary integral vanishes due to the boundary conditions
$\nn{\cdot}\nabla\alpha=0$ and $\vv{\cdot}\nn=0$. The two convective terms
arising in \eqref{test-1} and \eqref{test-damage} are to be handled
jointly, which gives 
\begin{align}&
\INT{\varOmega}{}\varphi_\Ee'(\Ee,\alpha){:}(\vv{\cdot}\nabla)\Ee
+\varphi_\alpha'(\Ee,\alpha)(\vv{\cdot}\nabla\alpha)\,\d x
\label{calculus-toward-pressure}
=\INT{\varOmega}{}\nabla\varphi(\Ee,\alpha){\cdot}\vv\,\d x
=\INT{\varGamma}{}\varphi(\Ee,\alpha)\vv{\cdot}\bm{n}\,\d S
-\INT{\varOmega}{}\varphi(\Ee,\alpha){\rm div}\,\vv\,\d x\,.
\end{align}
This then uses the boundary condition $\vv{\cdot}\bm{n}=0$ while
the bulk term is balanced with the corresponding pressure-type term  in
\eqref{test-1}. 

 Substituting \eqref{test-1} into \eqref{test-velocity} and summing it with
\eqref{test-damage} gives the mechanical energy balance
\begin{align}\nonumber
  &\frac{\d}{\d t}
  \INT{\varOmega}{}\!\!\linesunder{\frac\varrho2|\vv|^2}{kinetic}{$^{^{}}$energy}\!\!\!\!
  +\hspace{-.7em}\lineunder{\varphi(\Ee,\alpha)
    +\frac{\kappa}2|\nabla\alpha|^2}{stored energy}\hspace{-.7em}\,\d x
  \\[-.1em]&  +
  \INT{\varOmega}{}\!\lineunder{\!\!
\GM(\W)|\ZJEp|^2\!+\xi\big(\alpha,\W;\DT\alpha\big)
    +\bbD\EE(\vv){:}\EE(\vv)+\nu|\nabla\EE(\vv)|^p\!
   +\varkappa|\nabla\ZJEp|^2}{dissipation rate}\hspace{-.7em}\d x
  =\INT{\varOmega}{}\hspace{-.7em}\linesunder{\ff{\cdot}\vv-\phi(\theta){\rm div}\,\vv_{_{_{_{}}}}\!\!}{power of external load}{and of adiabatic effects}\hspace{-.7em}\d x
   \,.
  \label{energy+}\end{align}
Adding \eqref{ED-6+thermo} tested by 1 and using also the last boundary
condition in \eqref{BC3}, we obtain the total energy balance 
\begin{align}
  &\frac{\d}{\d t}
  \INT{\varOmega}{}\!\!\linesunder{\frac\varrho2|\vv|^2}{kinetic}{$^{^{}}$energy}\!\!\!
  +\hspace*{-.7em}\linesunder{\varphi(\Ee,\alpha)
    +\frac{\kappa}2|\nabla\alpha|^2}{stored}{$^{^{}}$energy}\hspace*{-.7em}
  +\!\!\!\!\!\!\!\!\!\!\!\linesunder{\W_{_{_{_{}}}}\!\!\!}{internal heat}{$^{^{}}$energy (enthalpy)}
  \!\!\!\!\!\!\!\!\!\!\!\d x
  =\INT{\varOmega}{}\!\!\linesunder{\ff\cdot\vv_{_{_{_{}}}}\!\!\!}{external}{$^{^{^{}}}$bulk load}\!\!\!\d x
  +\INT{\varGamma}{}\!\!\!\!\!\!\!\linesunder{
    {h_{\rm ext}}_{_{_{_{_{}}}}}\!\!\!}{power of}{external heating}\!\!\!\!\!\!\!\d S\,.
  \label{energy+++}\end{align}

The thermodynamical context of this model relies on an additive splitting
of the specific free energy $\psi$ into the purely mechanical part and the
thermal part $\phi$, cf.\ an example \eqref{psi} below.
An important attribute of the model, beside keeping the energetics
\eqref{energy+}--\eqref{energy+++}, is the entropy balance, and in
particular the {\it Clausius-Duhem inequality}  due to non-negativity
of the entropy production in \eqref{entropy-balance}. The specific 
entropy $\eta=-\psi_\theta'$ is an extensive variable (in JK$^{-1}$m$^{-3}$) and  its transport and
production is governed  by the entropy equation 
\begin{align}
  \theta\Big(\pdt\eta{+}{\rm div}\big(\vv\,\eta\big)\Big)
=\xi-{\rm div}\,{\bm j}\qquad\qquad\qquad
\label{entropy-eq}\end{align}
with $\xi$ denoting the heat production rate (here by mechanical dissipation
processes) and ${\bm j}$ the heat flux (here governed by the Fourier law
${\bm j}=-\bbK\nabla\theta$). The ultimate assumptions $\xi\ge0$ and $\bbK>0$
then ensure the entropy balance
\begin{align}
  \frac{\d}{\d t}\INT{\varOmega}{}\,\eta\,\d x
  =\INT{\varOmega}{}\,\frac{\xi-{\rm div}\,{\bm j}}\theta\,\d x
  +\INT{\varGamma}{}\,\eta\vv{\cdot}\nn\,\d S=
  \INT{\varOmega}{}\!\!\!\!\!\!\!\lineunder{\frac\xi\theta+\bbK\frac{|\nabla\theta|^2}{\theta^2}}{entropy production $\ge0$}\!\!\!\!\!\!\!\!\!\d x
+\INT{\varGamma}{}\!\!\!\!\linesunder{(\vv\eta-{\bm j}/\theta)}{entropy}{flux}\!\!\!\!{\cdot}\nn\,\d S
\label{entropy-balance}\end{align}
relying on positivity of temperature.
Substituting $\eta=-\phi'(\theta)$ into \eqref{entropy-eq}, we obtain
the heat-transfer equation 
\begin{align}
c(\theta)\DT\theta=\xi-{\rm div}\,{\bm j}\,-\,\theta
\eta\,{\rm div}\,\vv\ \ \ \ \text{ with the heat capacity }\ \ \
c(\theta)=-\theta\phi''(\theta)\,;
\label{heat-eq-}\end{align}
note that temperature (in K) is an intensive variable and is transported
by the material derivative while the adiabatic heat source/sink 
term $\theta\eta\,{\rm div}\,\vv$ occurs on the right-hand side due to
the compressibility of the solid or fluid. Furthermore,
the internal energy is given by the Gibbs relation
$\psi+\theta\eta$, and splits here into the purely mechanical part and the
thermal part $\W=\phi(\theta)-\theta\phi'(\theta)=:\gamma(\theta)$.
The thermal internal energy $\W$ in Jm$^{-3}$ is again an extensive
variable and is transported like \eqref{entropy-eq}, resulting here to the
equation
\begin{align}
  \pdt\W+{\rm div}\big(\vv\,\W\big)=\xi-{\rm div}\,{\bm j}
  +\phi(\theta){\rm div}\vv\,,
\label{heat-eq}\end{align}
which reveals the structure of \eqref{ED-6+thermo}.

\begin{remark}[{\sl Inelastic strain $\bm{P}$}]\label{rem-plast-strain}\upshape
   As already noted, we did not explicitly involve
  the inelastic strain $\bm{P}$ in \eqref{additive-rate}, neither
  in the list of initial conditions \eqref{IC}. Yet, if an initial
  condition $\bm{P}|_{t=0}^{}=\bm{P}_0^{}$ would be prescribed, we could
  substitute $\ZJEp=\ZJ{\bm{P}}$ into (\ref{ED+}b,c,e).
  It is important that, formulating the model in terms of $\bm{P}$ complies
  with the generally accepted concept of volume preserving (so-called isochoric)
inelastic deformations,
which is reflected by trace-free inelastic strain $\bm{P}$. This would be
ensured if $\bm{P}_0^{}$ is valued in $\R_{\rm dev}^{d\times d}$. Then $\bm{P}(t)$ remains 
 deviatoric  (=\,symmetric trace-free) for all $t>0$ because
\begin{align*}
  {\rm tr}\,\ZJ{\bm{P}}=({\rm tr}\,\bm{P})\!\DT{^{}}\ \ \ \ \text{ and }\ \ \ \ 
  \ZJ{\bm{P}}^\top=(\bm{P}^\top)\!\ZJ{^{}}
\end{align*}
and because ${\rm dev}\,\bm{S}$ which drives the evolution of $\bm{P}$ in
\eqref{ED-3+} with $\Ep=\ZJ{\bm{P}}$, is symmetric trace-free.
 This is a general property of the Zaremba-Jaumann derivative,
 while others, non-corotational (like Oldroyd's or Truesdell's) derivatives
 would not guarantee this property. Similarly, if $\Ee_0$
is valued in $\R_{\rm sym}^{d\times d}$, then remains
symmetric for all $t>0$ because $\ZJ\Ee^\top=(\Ee^\top)\!\ZJ{^{}}$ and 
because $\EE(\vv)$ driving the evolution of $\Ee$ in
\eqref{ED-2+} is symmetric.
\end{remark}

\section{A specialization towards phase transformation\\and phase-field fracture}
%        ~~~~~~~~~~~~~~~~~~~~~~~~~~~~~~~~~~~~~~~~~~~~~~~~~~~~~~~~~~~~~~~~~~~~~~
\label{sec-special}
For isotropic material as in Remark~\ref{rem-isotropic} combined with damage,
the above model together with the analysis in Section~\ref{sec-anal} below
covers the particular case of the free energy \eqref{psi-} as
\begin{align}
  \psi(\bm E,\alpha,\nabla\alpha,\theta)&
  =\frac d2K_\text{\sc e}^{}|{\rm sph}\,\Ee|^2+G_\text{\sc e}^{}(\alpha)
  |{\rm dev}\,\Ee|^2
  +\frac{1}{2\kappa}G_\text{\sc d}^{}(1{-}\alpha)^2
+\frac{\kappa}2|\nabla\alpha|^2+\phi(\theta)\,.
\label{psi}\end{align}
The parameters are bulk modulus $K_\text{\sc e}^{}$ as before,
a ``damageable'' shear modulus $G_\text{\sc e}^{}=G_\text{\sc e}^{}(\alpha)$
with $G_\text{\sc e}^{}(0)=G_\text{\sc e}'(0)=0$ and $G_\text{\sc e}^{}(1)\ge0$,
which is important to keep $\alpha$ valued in $[0,1]$. Further ingredients
in \eqref{psi} are  a fracture-toughness-like modulus $G_\text{\sc d}^{}>0$, 
and the thermal part $\phi$ of the free energy.  Moreover, the  
parameter $\kappa>0$ in \eqref{psi} occurring already in \eqref{ED+} is in
N=Pa\,m$^2$. The $G_\text{\sc d}^{}$-term gives rise
to a driving for healing processes, i.e.\ $\varphi_\alpha'(0,\alpha)>0$
if $\alpha<1$, cf.\ e.g.\ \cite[Sect.\,7.5]{KruRou19MMCM}
or \cite{Roub19MDDP}. The most popular choice is 
$G_\text{\sc e}^{}(\alpha)=G_\text{\sc e,0}^{}(\alpha^2{+}\epsilon^2)$
with some constant $G_\text{\sc e,0}^{}$,
which is known as an Ambrosio-Tortorelli functional \cite{AmbTor90AFDJ}.
While $G_\text{\sc e,0}^{}$ as well as $K_\text{\sc e}^{}$
are in Pa=J/m$^3$, the physical dimension of $G_\text{\sc d}^{}$ is N/m$^2$.
For $\kappa$ small, the damage tends to be localized along thin regions,
called cracks. Various modifications fitted to both initiation and propagation
of cracks are devised in literature and known under the name of
{\it phase-field} fracture, cf.\ \cite[Sect.\,3]{Roub19MDDP}.

Here the bulk modulus $K_\text{\sc e}^{}$ is considered independent of damage,
which can be interpreted as that the material is surely so compressed
that no tension (i.e.\ no opening also called 
fracture in More I) can occur. This means that the fracture might develop only
in Mode II (shearing) or Mode III (tearing). 
The situations when the material is well compressed are actually what
we have implicitly in mind not only for solid part but also for the
liquids (modelled here as semicompressible fluids) which can withstand only 
non-negative (or very small negative) pressures.

The rheology behind the equations (\ref{ED+}a--c) in this special
case can schematically be depicted by Figure~\ref{fig1}:
\begin{center}
\begin{my-picture}{25}{12}{fig1}
\psfrag{e}{\small $\EE(\uu)$}
\psfrag{s}{\small $\sigma$}
\psfrag{f}{\small $\bm f$}
\psfrag{s1}{\small $\sigma_{\rm sph}$}
\psfrag{s2}{\small $\sigma_{\rm dev}$}
\psfrag{e1}{\small ${\rm dev}\,\bmE$}
\psfrag{e2}{\small ${\rm sph}\,\Ee$}
\psfrag{e3}{\small ${\rm dev}\,\bm E$}
\psfrag{C1}{\small $K_\text{\sc e}^{}$}
\psfrag{C2}{\small $G_\text{\sc e}^{}$}
\psfrag{D1}{\small $G_\text{\sc v}^{}$}
\psfrag{D2}{\small $K_\text{\sc v}^{}$}
\psfrag{D3}{\small $\GM$}
\psfrag{q}{\small $\theta$}
\psfrag{a}{\small $\alpha$}
\psfrag{P}{\small $\bm{P}$}
\psfrag{spherical part}{\!\!\!\!\footnotesize\sf\begin{minipage}[t]{10em}\hspace*{0em}spherical part
\\\hspace*{0em}(Kelvin-Voigt)
\end{minipage}}
\psfrag{deviatoric part}{\!\!\!\!\footnotesize\sf\begin{minipage}[t]{10em}\hspace*{0em}deviatoric part
  \\[-.2em]\hspace*{0em}(Jeffreys or,
  \\[-.2em]\hspace*{0em}if $\GM=0$,
 \\[-.2em]\hspace*{0em}mere Stokes)
\end{minipage}}
\hspace*{-5em}\includegraphics[width=31em]{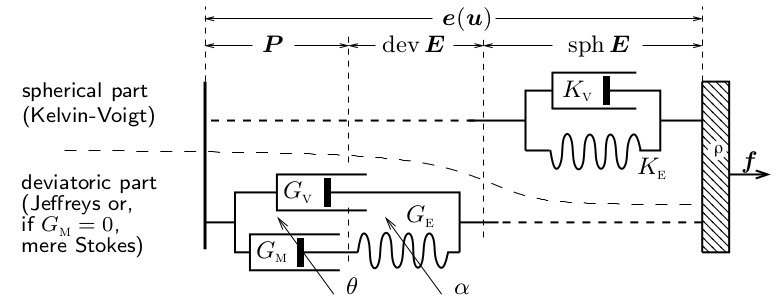}
\end{my-picture}
\nopagebreak
    {\small\sl\hspace*{-5em}Fig.~\ref{fig1}:~\begin{minipage}[t]{34em}
 A schematic 0-dimensional diagramme of the mixed rheology acting
differently on spherical (volumetric) and the deviatoric parts.
The Jeffreys rheology in the deviatoric part combines
Stokes' and Maxwell rheologies in parallel and may degenerate
to mere Stokes-type fluid if $\GM$ vanishes within melting.
\end{minipage}
} 
\end{center}

The ansatz \eqref{psi} is thermodynamically simple situation because the
mechanical and the thermal parts are coupled additively. This particularly
means that the entropy $\eta=-\psi_\theta'=-\phi'(\theta)$ does not depend on
mechanical variables and the internal
energy $\psi+\theta\eta=
\varphi(\Ee,\alpha)+\frac12\kappa|\nabla\alpha|^2+\phi(\theta)-\theta\phi'(\theta)$
also couples additively the mechanical and the heat parts, cf.\  
\eqref{energy+++} for $\W=\phi(\theta)-\theta\phi'(\theta)$.
This reveals that $\gamma(\theta)=\phi(\theta)-\theta\phi'(\theta)$ and that
the heat capacity $c=c(\theta)=\gamma'(\theta)=-\theta\phi''(\theta)$
in \eqref{heat-eq}
also does not depend on mechanical variables. In our Stefan problem, the heat
capacity is a distribution containing a Dirac measure at the phase-transition
temperature $\theta=\theta_\text{\sc pt}$, cf.\ Figure~\ref{fig2}:

\begin{center}
\begin{my-picture}{25}{13}{fig2}
\psfrag{q}{\small $\theta=\beta(\W)$}
\psfrag{qsf}{\small $\theta_\text{\sc pt}$}
\psfrag{Gviscous}{\small $[\GM\circ\gamma](\theta)$}
\psfrag{T}{\small $\W\in\gamma(\theta)$}
\psfrag{TS}{\small $\W_\text{\sc s}$}
\psfrag{TF}{\small $\W_\text{\sc l}$}
\psfrag{solid}{\footnotesize\sf solid}
\psfrag{solid elastic}{\!\footnotesize\sf\begin{minipage}[t]{10em}\hspace*{0em}solid with\\[-.3em]\hspace*{0em}dominating\\[-.3em]\hspace*{0em}elasticity and\\[-.3em]\hspace*{0em}possible fracture
\end{minipage}}
\psfrag{creep}{\!\footnotesize\sf\begin{minipage}[t]{10em}\hspace*{0em}creep
  in solid\\[-.3em]\hspace*{0em}in temperature close\\[-.3em]\hspace*{0em}to melting
\end{minipage}}
\psfrag{PT}{\!\footnotesize\sf\begin{minipage}[t]{10em}\hspace*{0em}solid-liquid
  phase\\[-.3em]\hspace*{0em}transformation\\[-.3em]\hspace*{0em}(melting/solidification)
\end{minipage}}
\psfrag{temperature}{\footnotesize temperature}
\psfrag{melting}{\footnotesize melting/freezing} 
\psfrag{liquid}{\footnotesize\sf liquid} 
\psfrag{latent}{\footnotesize\sf latent}
\psfrag{heat}{\footnotesize\sf heat $\ell$}
\psfrag{fluid}{\footnotesize\sf fluid (liquid)}
  \psfrag{t1}{$\theta$}
\vspace*{5em}
\hspace*{-8em}{\includegraphics[width=37em]{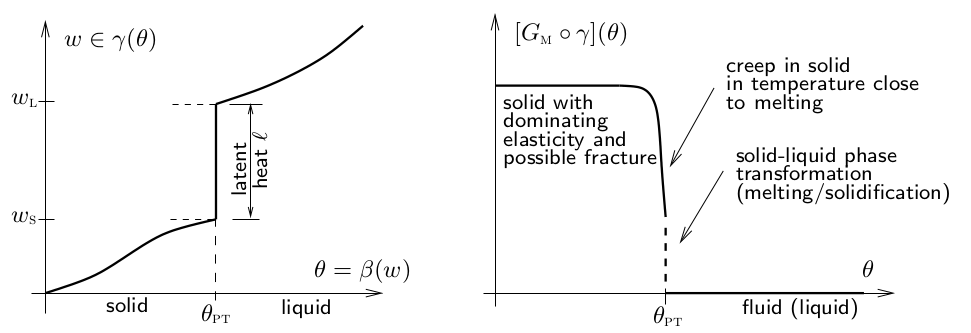}}
\end{my-picture}
\nopagebreak
    {\small\sl\hspace*{-5em}Fig.~\ref{fig2}:~\begin{minipage}[t]{38em}
    A basic philosophy of the model: dependence of
 the enthalpy $\W$ on temperature $\theta$ in the classical Stefan
    problem (left) and the dependence of
    the Maxwellian viscosity shear modulus $\GM$ (right). As $\gamma$ is
    multivalued at the solid-liquid phase-transition temperature
    $\theta=\theta_\text{\sc pt}$, the heat capacity
    $\gamma'(\cdot)$ contains a Dirac measure supported at
    $\theta=\theta_\text{\sc pt}$.
\end{minipage}
} 
\end{center}

The dependence on the Maxwellian modulus $\GM$ (controlling creep effects)
on temperature in the solid state as schematically depicted on
Figure~\ref{fig2}-right is in the case of ice is supported by experiments
showing that when the ice temperature is close to melting, the creep effects
accelerate by several orders, cf.\ 
\cite[Fig.3]{MelTes69ETCI}. Across the melting temperature, $\GM$ falls to
zero and solid ice becomes fluidic water. This discontinuity in terms of
temperature is in fact possible to turn into continuity when expressed
in terms of the enthalpy $\W$, cf.\ Figure~\ref{fig3}.

An example for the damage dissipation potential
  $\zeta(\W;\DT{\alpha})$ as another important ingredient of the model
not explicitly depicted in Figure~\ref{fig1} can be 
\begin{align}\label{dissip-exa}
  \zeta\big(\W;\DT{\alpha}\big)
=-\sigma_\text{\sc f}^{}\min(0,\DT\alpha)+A(\W)\max(0,\DT\alpha)^2
+\epsilon\DT\alpha^2
\end{align}
with $A:\R\to\R$ continuous non-increasing very large for
solid regions (i.e.\ for $\W\le\W_\text{\sc s}$) while zero for
liquid regions (i.e.\ for $\W\ge\W_\text{\sc l}$). The meaning of
$\sigma_\text{\sc f}^{}>0$ is a so-called fracture toughness, i.e.\ the energy
per volume needed to damage the solid material completely, while
$A=A(\W)\ge0$ influences the healing rate and, if $A\to0$, then $\alpha$
tends to come fast to 1 (i.e.\ to undamaged material) due to the term
$G_\text{\sc d}^{}(1{-}\alpha)^2/(2\kappa)$ in the stored energy in
\eqref{psi}. This holds if $\epsilon=0$, while later we will assume
$\epsilon>0$ small for the analytical reasons, cf.\ \eqref{ass:2} below.
Thus, in liquid phase, $\alpha\sim1$ and therefore the role of
damage is naturally eliminated and solidification leads to
the originally undamaged solid.
In general, $\zeta$ may depend also on $\alpha$ itself through an
$\alpha$-dependence of both $\sigma_\text{\sc f}^{}$ and $A$. 

 \begin{center}
\begin{my-picture}{25}{19}{fig3}
\psfrag{q}{\small $\theta=\beta(\W)$}
\psfrag{qsf}{\small $\theta_\text{\sc pt}$}
\psfrag{Gviscous}{\small $\GM=\GM(\W)$}
\psfrag{TS}{\small $\W_\text{\sc s}$}
\psfrag{TF}{\small $\W_\text{\sc l}$}
\psfrag{solid}{\footnotesize\sf solid}
\psfrag{temperature}{\footnotesize\sf temperature}
\psfrag{melting}{\footnotesize\sf melting/freezing} 
\psfrag{qsf}{\small $\theta_\text{\sc pt}$}
\psfrag{liquid}{\footnotesize\sf liquid} 
\psfrag{latent}{\footnotesize\sf latent}
\psfrag{heat}{\footnotesize\sf heat $\ell$}
\psfrag{fluid}{\footnotesize\sf fluid}
\psfrag{w}{$\W$}
\vspace*{5em}
\hspace*{0em}{\includegraphics[width=21em]{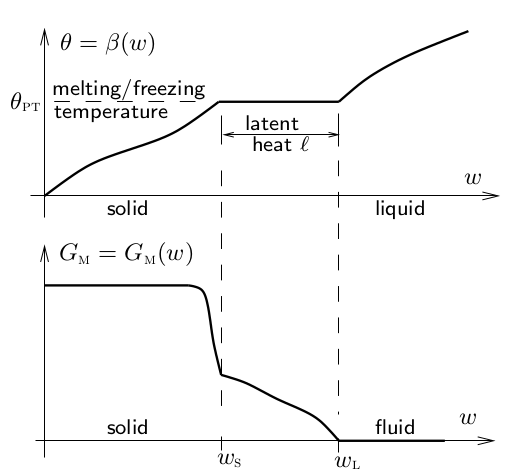}}
\end{my-picture}
\nopagebreak
\\
  {\small\sl\hspace*{-5em}Fig.~\ref{fig3}:~\begin{minipage}[t]{24em}
Temperature $\theta$ and the Maxwellian viscous creep modulus $\GM$ depending
      continuously on the enthalpy $\W$.
\end{minipage}
} 
\end{center}
  
\begin{remark}[{\sl Semi-compressible fluids II}]\label{rem-fluids-II}\upshape
  When the deviatoric elastic response vanishes so that
  $\varphi(\bm E)=dK_\text{\sc e}|{\rm sph}\,\bm E|^2/2$, we again obtain the
  model for semi-compressible fluids as in Remark~\ref{rem-semi-compress} except
  that the momentum balance \eqref{semi-compres-2} is augmented by the
higher-gradient term ${\rm div}^2(\nu|\nabla\EE(\vv)|^{p-2}\nabla\EE(\vv))$,
related with the mentioned 2nd-grade nonsimple fluid concept. Indeed, if
$\GM(0)=0$, then ${\rm dev}\,\Ee=0$ and the elastic stress
simplifies to $\bm{S}=\varphi_\Ee'(\Ee,\alpha)=dK_\text{\sc e}{\rm sph}\,\Ee
  =K_\text{\sc e}({\rm tr}\,\Ee)\bbI
  =p\bbI$ when denoting $p=K_\text{\sc e}{\rm tr}\,\Ee$, and then
  \eqref{ED-2+} yields $\DT{p}/K_\text{\sc e}={\rm div}\,\vv$, which is just
  \eqref{semi-compres-1}.
  Also, \eqref{ED-1+} with omitted the internal variable 
  $\alpha$ (as indeed natural in fluids) gives \eqref{semi-compres-2}.
 Here, the fluidic region can be, at least formally, combined with the solid regions
  in one ``monolithic'' model.
  Pressure (longitudial) waves can propagate thorough solid and liquid regions,
  being possibly refracted and partly reflected on the solid-liquid interfaces,
  while shear waves can propagate only in solid regions, being fully reflected
  on the solid-liquid interfaces or possibly partly transformed to longitudial
  waves in fluidic regions. Of course, these solid-liquid interfaces may
  evolve during solidification or melting.  Moreover, the
  mentioned higher-gradient term likely causes an {\it anomalous dispersion}
  of these waves 
  (like if the 2nd-grade nonsimple term would act in the stored energy, cf.\
  \cite[Remark~6.3.6]{KruRou19MMCM}).
 \end{remark}

\section{Two-phase Stefan problem involving 
  supercooling/superheating by phase relaxation}\label{sec-super}
%        ~~~~~~~~~~~~~~~~~~~~~~~~~~~~~~~~~~~~~~~~~~

The basic Stefan problem was enriched by various fine effects, in
particular towards  a kinetic  supercooling and superheating
(sometimes referred as undercooling and overheating) situations
\cite{ColGra93HPCP,ColRec02CSPP,ColSpr95PFMZ,Schi00SCRC,Visi85SPPH,Visi85SSEP},
relying on that melting or solidification proceed at a rate that increases
with the difference between the actual temperature and the equilibrium
  (=\,phase transformation) 
  temperature \cite{Chal64PS}. This enrichment (also called a relaxed
  Stefan problem) will also bring
  analytical benefits especially in the context of mechanical
  coupling with $L^1$-heat sources. In particular, it allows for strong
  convergence of enthalpy $\W$ so that the mechanical properties
  may depend directly on $\W$ and may thus exhibit jumps when temperature
  crosses the phase-transformation threshold $\theta_\text{\sc pt}$.
  Even more important, the usual limit passage in the relation
  $\theta=\beta(\W)$ relying on estimation of $\nabla\W
  =\gamma'(\theta)\nabla\theta$ from an estimate of $\nabla\theta$ does not
  work for Stefan problem where $\gamma'(\theta)$ is unbounded while
  the estimate of $\frac{\partial}{\partial t}\theta$ as in
  \cite[Remark~8.2.4]{KruRou19MMCM} does not work.

  Here, of course, we present the supercooling/superheating in the convective
  variant and in the context of the $L^1$-theory for heat-transfer problem.

  Having in mind the situation from Fig.~\ref{fig2}, i.e.\ the so-called
  two-phase Stefan problem with only one phase-transformation temperature 
  $\theta_\text{\sc pt}$, we split the set-valued function $\gamma$
  into a continuous part and the multi-valued Heaviside function $H$, i.e.
\begin{subequations}\label{split}\begin{align}\label{split-H}
 &\gamma(\theta)=\widetilde\gamma(\theta)+\ell H\Big(\,\frac\theta{\theta_\text{\sc pt}\!\!}-1\Big)\ \ \
    \text{ with }\ \ H(\theta)=\begin{cases}\ \ 0&\text{if }\ \theta<0\,,\quad
    \text{\sf (complete solid)}
      \\[-.2em]
        [0,1]\!\!&\text{if }\ \theta=0\,,\quad\text{\sf (in between solid and fluid)}
        \\[-.2em]\ \ 1&\text{if }\ \theta>0\,,\quad
        \text{\sf (complete fluid)}
    \end{cases}
\intertext{and correspondingly we split also the heat part of the free energy
      and the internal energy (enthalpy) as}
  &\phi(\theta)=\widetilde\phi(\theta)
    -\ell\Big(\,\frac\theta{\theta_\text{\sc pt}\!\!}-1\Big)^+\ \ \text{ and }\ \ 
    \label{split-phi}
    \\  &
    \W=\vartheta+\ell\chi\ \ \text{ with }\ \
    \vartheta=\widetilde\gamma(\theta)=\widetilde\phi(\theta)-\theta\widetilde\phi'(\theta)\ \
     \label{split-w}\\  &\hspace*{5.7em}\text{ and with}\ \
    \chi\in H\Big(\,\frac\theta{\theta_\text{\sc pt}\!\!}-1\Big),\ \
    \text{ or equivalently }\
  H^{-1}(\chi)\ni\frac\theta{\theta_\text{\sc pt}\!\!}-1\,,
  \label{split-3}\end{align}\end{subequations}
where $(\cdot)^+$ denotes the positive part of its argument. Now
$\widetilde\gamma(\cdot)$ is a continuous single-value function, in contrast
to $\gamma$. The new variable $\chi$ has a meaning of a volume fraction. 

The model from Sections~\ref{sec-fuller}--\ref{sec-special} will now be slightly
modified by  ``relaxing'' the last inclusion in \eqref{split-3} as
\begin{align}
  \NU\DT\chi+H^{-1}(\chi)\ni
  \Upsilon\Big(\,
  \frac\theta{\theta_\text{\sc pt}\!\!}-1
  \Big)
\label{relaxed}\end{align}
with $\NU>0$ a (typically small) relaxation time
and with $\Upsilon$ some continuous increasing function with $\Upsilon(0)=0$.
For $\NU=0$, we obtain exactly the Stefan model \eqref{split}.
In literature, usually $\Upsilon(\theta)=\theta$ but here we consider a general
$\Upsilon$ with some sublinear growth, cf.\ \eqref{Y-ass} below, which will be needed for the test \eqref{test-dchi-dt} which is then needed for the
strong convergence of $\chi$'s. 

Let us note that we used the convective time derivative in \eqref{relaxed},
which is related with the attribute of the ``volume-fraction'' variable $\chi$
as  an intensive variable taking values in $[0,1]$. In contrast to it,
$\ell\chi$ in \eqref{split-3} in the dimension J/m$^3$ is an extensive
variable and can be summed up with $\vartheta$, giving $\W$ 
which should be transported as $\W$ in \eqref{heat-eq}.

Altogether, we now consider the system (\ref{ED+}a-e) with $\phi$ from
\eqref{split-phi} and, instead of \eqref{beta}, now with $\W$ from
\eqref{split-w} completed with the volume-fraction equation \eqref{relaxed}.

The mechanical-energy balance \eqref{energy+} is to be obtained as before by
testing the four equations/inclusion (\ref{ED+}a--d) respectively by $\vv$,
${\bm S}$, $\ZJEp$, and $\DT\alpha$. The total-energy conservation
\eqref{energy+++} is again by summing it with \eqref{ED-6+thermo} tested by 1.
The volume-fraction equation \eqref{relaxed} is not directly included in the
energetics. It also needs its own initial condition $\chi|_{t=0}^{}=\chi_0^{}$.
Then we can prescribe the initial temperature  $\theta|_{t=0}^{}=\theta_0^{}$
because the initial enthalpy is the uniquely determined as
$\W|_{t=0}^{}=\widetilde\gamma(\theta_0^{})+\ell\chi_0^{}$.

\begin{remark}[{\sl Multi-phase Stefan problem}]\upshape
  Actually, a generalization to more than  two phases in materials 
undergoing  several phase transformations at several
phase-transformation temperatures might be routinely possible by using more
than one volume fraction $\chi$'. This occurs in 
metals when a recrystallization with a certain latent heat (and with changing
mechanical properties) occurs still in the solid phase 
at temperature below an actual solid-fluid transformation
(i.e.\ melting/solidification), e.g.\ iron recrystalizes at $\sim450^\circ$C
 (i.e.\ a Stefan problem in solid phase) 
while melts at $1538^\circ$C  (i.e.\ a Stefan problem accompained by
solid-fluid transition).
Rather for notational simplicity, we are considering only the most
interesting scenario, i.e. an only two-phase problem with a
solid-fluid transition.
\end{remark}

\section{Analysis for incomplete melting by time discretisation}\label{sec-anal}
%        ~~~~~~~~~~~~~~~~~~~~~~~~~~~~~~~~~~~~~~~~~~~~~~~~~~~~~~

Although the semi-compressible fluids resulting after the complete
melting are well amenable to analysis, the transition
from solid to completely melted fluid (or reversely from the
semi-compressible fluid through freezing to solid) is analytically difficult,
cf.\ Remark~\ref{rem-fluids-II}. These troubles are not surprising
and occur in other similar situations where some coercivity continuously
degenerates like in contact mechanics (Coulomb friction).
We thus perform the analysis only for
incomplete melting, cf.\ the assumption \eqref{GM-ass} below, which can be
understood as a certain regularization of the desired model.
 It may be particularly relevant in situations
when $\GM$ drops by many orders within solid-fluid transition
but anyhow stays positive; e.g. in rocks $\GM\sim10^{18-25}$Pa\,s while
in magma $\GM\sim10^{1-7}$Pa\,s.
Note that, if $\GM=0$, $\ZJEp$ is not controlled (except its gradient but
the coefficient $\varkappa>0$ is expected small) and
$\Ee$ is not linked with $\EE(\vv)$ through \eqref{ED-2+} so that, naturally,
$\Ee\sim{\bm0}$ as expected in fluids and $\ZJEp\sim\EE(\vv)$. 
In other words,
the incomplete melting is expectedly an acceptable modelling
 approximation if $\GM$ drops much below $G_\text{\sc v}^{}$, cf.~Figure~\ref{fig1}, or if the flow
is not too fast (i.e.\ not exhibiting too big shear rates). 
On the other
hand, the nonsimple-material concept allows us to analyse
even a complete damage as far as the elastic response, which otherwise is
possible only in rather special rate-independent situations,
cf.\ \cite{BoMiRo09CDPS} or \cite[Prop.~4.3.22]{MieRou15RIST}.
In fact, this is possible due to the mentioned viscous response which (even
if generalized for being damage dependent) is to stay always incomplete.

We will use the standard notation concerning the Lebesgue and the Sobolev
spaces, namely $L^p(\varOmega;\R^n)$ for Lebesgue measurable functions
$\varOmega\to\R^n$ whose Euclidean norm is integrable with $p$-power, and
$W^{k,p}(\varOmega;\R^n)$ for functions from $L^p(\varOmega;\R^n)$ whose
all derivative up to the order $k$ have their Euclidean norm integrable with
$p$-power. We also write briefly $H^k=W^{k,2}$.
 The notation $2^*$ will denote the exponent from the embedding
$H^1(\varOmega)\subset L^{2^*}(\varOmega)$, i.e.\ $2^*=2d/(d{-}2)$. 
Moreover, for a Banach space
$X$ and for $I=[0,T]$, we will use the notation $L^p(I;X)$ for the Bochner
space of Bochner measurable functions $I\to X$ whose norm is in $L^p(I)$, 
and $H^1(I;X)$ for functions $I\to X$ whose distributional derivative
is in $L^2(I;X)$. Furthermore, $C_{\rm w}(I;X)$ will denote
the Banach space of weakly continuous functions $I\to X$, and $C_{\rm w*}(I;X)$
of weakly* continuous if $X$ has a predual, i.e.\ there is $X'$ such that
$X=(X')^*$ where $(\cdot)^*$ denotes the dual space. Occasionally, we will use
$L_{\rm w*}^p(I;X)$ the space of weakly* measurable mappings $I\to X$;
recall that $L_{\rm w*}^p(I;X)=L^p(I;X)$ if $X$ is separable reflexive.
Also, ${\rm Meas}(\barOmega)=C(\barOmega)^*$ is a space of Borel measures
on the closure $\barOmega$ of $\varOmega$. 

We will assume, with some $\epsilon>0$ arbitrarily small, that
\begin{subequations}\label{ass}\begin{align}\nonumber
    &\varphi:\R^{d\times d}\times[0,1]\to\R\ \text{ twice continuously
      differentiable, bounded from below with}
    \\&\nonumber\qquad\quad
    \forall(\Ee,\alpha)\in\R_{\rm sym}^{d\times d}\times[0,1]:
    \ \ |\varphi_\Ee'(\Ee,\alpha)|^2\le\varphi(\Ee,\alpha)/\epsilon,
\\&\label{ass:1-}\qquad\qquad\hspace*{10.5em}
|\varphi_\alpha'(\Ee,\alpha)|\le(1{+}|\Ee|^2)/\epsilon,\ \ \ \ 
\varphi_\alpha'(\Ee,0)\le0,\ \ \ \varphi_\alpha'(\Ee,1)\ge0,
   \\&\nonumber\widetilde\phi:\R\to\R\ \text{ continuously differentiable, }\
    \widetilde\phi(0)=0\,,
     \\&\qquad\quad\label{ass:1}
 \text{$\widetilde\gamma:\theta\mapsto\widetilde\phi(\theta){-}\theta\widetilde\phi'(\theta)$ 
 and $\widetilde\gamma^{-1}(\cdot)$ have at most linear growth,}
   \\&\label{GM-ass}
   \GM:\R\to\R\ \text{ continuous and bounded with }\ \inf\GM(\cdot)>0\,,
       \\&\label{Y-ass}
\Upsilon:\R\to\R\ \text{ continuous, increasing},\ \Upsilon(0)=0,\
\sup_{\R}\frac{|\Upsilon(\cdot)|}{1{+}|\cdot|^{1/2}}<\infty\,,
  \\&\nonumber
\zeta:[0,1]\times\R\times\R\to\R\
\text{ continuous with}
\ \ \zeta(\alpha,\W;\cdot):\R\to\R
\ \text{ convex and smooth on $\R\setminus\{0\}$,}
\\&\qquad\quad\epsilon
\big|\DT\alpha\big|^2\le\zeta\big(\alpha,\W;\DT\alpha\big)\le
\big(1+\big|\DT\alpha\big|^2\big)/\epsilon\,,
\label{ass:2}
\\
&\beta:\R\to\R\ \text{ continuous, nondecreasing, with at most linear
  growth,}
\\
&\bbK:\R^2\to\R\ \text{ continuous, bounded, }\ \  \inf\bbK>0\,,
\\\label{ass:4}
&p>d,\ \ \varrho,\ \kappa,\ \varkappa,\
\nu,\ \theta_\text{\sc pt}>0,
\ \bbD\in\R^{d^4}\ \text{ symmetric positive definite},
\\\nonumber
&\vv_0\!\in\! L^2(\varOmega;\R^d),\ \
         {\Ee}_0\!\in\! H^1(\varOmega;\R_{\rm sym}^{d\times d}),\ \
        \alpha_0\!\in\! H^1(\varOmega),\ \ 0\le\alpha_0\le1,
 \\\label{ass:5}
 &\hspace{7.5em}
 \theta_0\!\in\! L^1(\varOmega),\ \,\theta_0\ge0,\ \,
\chi_0\!\in\! W^{1,q}(\varOmega),\,
\ 0\le\chi_0\le1\ \text{ with }1\le q<(d{+}2)/(d{+}1),
\\\label{ass:6}&\ff\in L^1(I;L^2(\varOmega;\R^d))\,,\ \ \ \ \
h_{\rm ext}^{}\in L^1(I{\times}\varGamma)\,,\ \ \ \ h_{\rm ext}^{}\ge0\,.
\end{align}\end{subequations}

We have formulated our growth assumption 
\eqref{ass:1-} to be valid for $d=3$ and $d=2$ too, but especially for the
latter case they can be weakened.
Let us emphasize that we do not assume $\varphi$ convex, which
allows for realistic damage models where $\varphi$ is intentionally always
nonconvex.

For our time-discretisation method used in below, we will still add a
certain mild qualification of $\varphi$, cf.\ the Remark~\ref{rem-example}
below, namely that
   \begin{subequations}\label{phi-eps}
    \begin{align}\nonumber
      &\forall\eps>0\
      \exists\varphi_\eps:\R_{\rm sym}^{d\times d}\times[0,1]\to\R\
     \text{ continuously differentiable, bounded from below with}
    \\&\label{ass:1-eps}\quad
    \forall(\Ee,\alpha)\in\R_{\rm sym}^{d\times d}\times[0,1]:
    \ \ |[\varphi_\eps]_\Ee'(\Ee,\alpha)|^2\le\varphi_\eps(\Ee,\alpha)/\epsilon,
\\&\hspace*{12.5em}
|[\varphi_\eps]_\alpha'(\Ee,\alpha)|^2\le(1{+}|\Ee|^2)/\epsilon,\ \ \
[\varphi_\eps]_\alpha'(\Ee,0)\le0,\ \ \ [\varphi_\eps]_\alpha'(\Ee,1)\ge0,
    \\&\label{semi-convex}\qquad
    \exists K=K(\eps):\qquad
    (\Ee,\alpha)\mapsto\varphi_\eps(\Ee,\alpha)+\frac12K|\alpha|^2
    \ \text{ is convex,}
   \\[-.3em]&\label{phi-eps0}\quad
    \forall\,\Ee\in\R_{\rm sym}^{d\times d},\ \alpha\in[0,1]:\ \ \
\big|\varphi_\eps(\Ee,\alpha)-\varphi(\Ee,\alpha)\big|\le o(\eps)(1+|\Ee|^2)\,,
    \\&\hspace*{12.5em}\label{phi-eps1}
    \big|[\varphi_\eps]_\Ee'(\Ee,\alpha)-\varphi_\Ee'(\Ee,\alpha)\big|\le
    o(\eps)(1+|\Ee|)\,,\ \text{ and}
    \\&\hspace*{12.5em}\big|[\varphi_\eps]_\alpha'(\Ee,\alpha)-\varphi_\alpha'(\Ee,\alpha)\big|\le o(\eps)(1+|\Ee|^2)\,.
  \label{phi-eps2}\end{align}\end{subequations}

  We should realize that $H^{-1}$ in \eqref{relaxed} is the normal cone to
$[0,1]$, i.e.\ the subdifferential to the indicator function to the interval
$[0,1]$, and thus \eqref{relaxed} means $\chi$ is valued in $[0,1]$ and
satisfies the variational inequality
\begin{align}
  \Big(
  \Upsilon\big(\,
  \frac\theta{\theta_\text{\sc pt}\!\!}-1
  \big)
  -\NU\pdt\chi-\NU\vv{\cdot}\nabla\chi\Big)(\widetilde\chi-\chi)\ge0
\label{relaxed+}\end{align}
to be valid a.e.\ on $I{\times}\varOmega$ for all $\widetilde\chi$ valued in
$[0,1]$. As $\chi$ is however not smooth in space/time, we have to understand
\eqref{relaxed+} only weakly by using Green's formula
$\int_\varOmega\vv{\cdot}\nabla\chi(\widetilde\chi-\chi)\,\d x=
\int_\varOmega\frac12({\rm div}\,\vv)\chi^2
-\chi{\rm div}(\vv\widetilde\chi)\,\d x$ and by-part integration in time for
the term $\frac{\partial}{\partial t}\chi(\widetilde\chi{-}\chi)$, resulting
thus to \eqref{VI} in the following:

\begin{definition}[Weak solutions to the system
      (\ref{ED+}a-e)--(\ref{BC})--(\ref{IC}) with (\ref{split}b,c)--(\ref{relaxed}).]\label{def}
  A six-tuple\linebreak$(\vv,\Ee,\Ep,\alpha,\theta,\chi)$ with
\begin{subequations}\begin{align}
    &\vv\!\in\! C_{\rm w}(I;L^2(\varOmega;\R^d))\cap 
    L^p(I;W^{2,p}(\varOmega;\R^d))
    \ \mbox{ with \ $\frac{\partial}{\partial t}$}\vv\!\in\! L^{p'}\!(I;W^{2,p}(\varOmega;\R^d)^*)+L^1(I;L^2(\varOmega;\R^d)),\!
    \\&\Ee\in C_{\rm w}(I;H^1(\varOmega;\R_{\rm sym}^{d\times d}))\,\cap\,H^1(I,L^2(\varOmega;\R_{\rm sym}^{d\times d}))
    \,,
    \\& \ZJEp\in L^2(I;H^1(\varOmega;\R_{\rm dev}^{d\times d}))
    \,,
    \\&\alpha\in C_{\rm w}(I;H^1(\varOmega))\,\cap\,H^1(I;L^2(\varOmega))
  \ \ \text{ and }\ 0\le\alpha\le1\ \text{ a.e. on }I{\times}\varOmega\,,
   \\&\theta\in C_{\rm w}(I;L^1(\varOmega))\,\cap\,L^r(I;W^{1,r}(\varOmega))\ \ \
    \mbox{ with $1\le r<\frac{d{+}2}{d{+}1}$},
      \ \ \text{and}\ \theta\ge0\ \text{ a.e. on }I{\times}\varOmega\,,
    \ \ \text{ and }\label{temp-w-cont}
    \\&\chi\in L^\infty(I{\times}\varOmega)\,\cap\,H^1(I;H^1(\varOmega)^*)\ \
    \mbox{ with \ \ 
      $0\le\chi\le1$ a.e.\ on $I{\times}\varOmega$}
\end{align}\end{subequations}
 will be called a weak solution to the boundary-value problem
 (\ref{ED+}a-e)--\eqref{BC} with (\ref{split}b,c)--\eqref{relaxed}
   and with the initial conditions \eqref{IC}  and $\chi|_{t=0}^{}=\chi_0$  if
 ${\bm S}=\varphi_\Ee'(\Ee,\alpha)\in L^2(I;H^1(\varOmega;\R_{\rm sym}^{d\times d}))$, $\varphi_\alpha'(\Ee,\alpha)\in L^2(I{\times}\varOmega)$,
   $\Delta\alpha\in L^2(I{\times}\varOmega)$,
    $\nn{\cdot}\vv=0$ on $I{\times}\varGamma$,
   $\W=\widetilde\gamma(\theta){+}\ell\chi\in C_{\rm w}(I;L^1(\varOmega))$
    with $\frac{\partial}{\partial t}\W\in
    L^1(I;H^{d+1}(\varOmega)^*)$, and if 
  \begin{subequations}\label{def+}\begin{align}\nonumber
      &\INT{0}{T}\!\!\INT{\varOmega}{}\bigg(\varrho\Big(\!(\vv{\cdot}\nabla)\vv
      {+}\frac12({\rm div}\,\vv)\vv\Big){\cdot}\widetilde\vv
      +\big({\bm S}
       {+}\bbD\EE(\vv){+}{\bm K}\big){:}\EE(\widetilde\vv)
  +\psi(\Ee,\alpha,\theta){\rm div}\widetilde\vv 
        \\[-.4em]&\hspace{6em}
    {+}\nu|\nabla\EE(\vv)|^{p-2}\nabla\EE(\vv)\Vdots\nabla\EE(\widetilde\vv)
   -\varrho\vv{\cdot}\pdt{\widetilde\vv}\bigg)\,\d x\d t=\INT{\varOmega}{}\varrho\vv_0{\cdot}\widetilde\vv(0)\,\d x
      +\INT{0}{T}\!\!\INT{\varOmega}{}\ff{\cdot}\widetilde\vv\,\d x\d t
\intertext{with ${\bm K}$ from \eqref{ED-1+} for all $\widetilde\vv\in 
        H^1(I;L^2(\varOmega;\R^d))\cap L^p(I;W^{2,p}(\varOmega;\R^d))
       $ with $\nn{\cdot}\widetilde\vv=0$ on $I{\times}\varGamma$ and
        $\widetilde\vv(T)={\bm0}$, and if}
&\INT{0}{T}\!\!\INT{\varOmega}{}
\big(\GM(\W)\ZJEp{+}{\rm dev}\,\bm{S}\big){:}\widetildeEp
+\varkappa\nabla\ZJEp\Vdots\nabla\widetildeEp
\,\d x\d t=0
\intertext{holds for all $\widetildeEp\in
  H^1(I{\times}\varOmega;\R_{\rm dev}^{d\times d})$, and if}
     &\INT{0}{T}\!\!\INT{\varOmega}{}\zeta(\alpha,\W;\widetilde{\alpha})
      +\varphi_\alpha'(\Ee,\alpha)
      \big(\widetilde{\alpha}{-}\DT\alpha\big)
     +\kappa\nabla\alpha{\cdot}\nabla\widetilde{\alpha}
      +\kappa\Delta\alpha(\vv{\cdot}\nabla\alpha)\,\d x\d t
      \nonumber\\[-.3em]&\qquad\qquad\qquad\qquad\qquad\qquad
      +\INT{\varOmega}{}
      \frac{\kappa}2|\nabla\alpha_0|^2\,\d x\ge
    \INT{\varOmega}{}\frac{\kappa}2|\nabla\alpha(T)|^2\,\d x
    +\INT{0}{T}\!\!\INT{\varOmega}{}\zeta\big(\alpha,\W;
    \DT\alpha\big)\,\d x\d t
   \label{VI}
\intertext{holds for all
  $\widetilde{\alpha}\in L^2(I;H^1(\varOmega))$, and}\nonumber
&\INT{0}{T}\!\!\INT{\varOmega}{}
\big(\bbK(\alpha,\W)\nabla\theta{-}\vv\,\W\big){\cdot}
\nabla\widetilde\W-\Big(\GM(\W)|\ZJEp|^2
{+}\bbD\EE(\vv){:}\EE(\vv){+}\nu|\nabla\EE(\vv)|^p
+\phi(\theta){\rm div}\,\vv\Big)\widetilde\W
-\W\pdt{\widetilde\W}\,\d x\d t
\\[-.4em]&\label{ED-6+weak}
\hspace*{18em}
=\INT{\varOmega}{}\big(\widetilde\gamma(\theta_0){+}\ell\chi_0\big)
\widetilde\W(0)\,\d x
+\INT{0}{T}\!\!\INT{\varGamma}{} h_{\rm ext}\widetilde\W\,\d S\d t
\intertext{for any $\widetilde\W\in W^{1,\infty}(I{\times}\varOmega)$ with
$\widetilde\W(T)=0$, and furthermore}\nonumber
&\INT{0}{T}\!\!\INT{\varOmega}{}\!\Upsilon
\Big(\,\frac\theta{\theta_\text{\sc pt}\!\!}-1\Big)
(\widetilde\chi-\chi)+\NU\chi\pdt{\widetilde\chi}
-\frac\NU2({\rm div}\,\vv)\chi^2+\NU\chi{\rm div}(\vv\widetilde\chi)\,\d x\d t
\\[-.5em]&\label{ED-7+weak}
\hspace*{15.3em}
+\INT{\varOmega}{}\frac\NU2\chi_0^2-\NU\chi_0^{}\chi(0)\,\d x
\ge\INT{\varOmega}{}\frac\NU2\chi(T)^2-\NU\chi(T)\widetilde\chi(T)\,\d x
  \end{align}\end{subequations}
  for any $\widetilde\chi\in W^{1,\infty}(I{\times}\varOmega)$ with
  $0\le\widetilde\chi\le1$ a.e.\ on $I{\times}\varOmega$, and eventually also
$\ZJ\Ee+\ZJEp=\EE(\vv)$ holds a.e.\ on $I{\times}\varOmega$
and  $\Ee(0)=\Ee_0$ a.e.\ on $\varOmega$.
\end{definition}

\begin{theorem}
Let \eqref{ass} and \eqref{phi-eps} be valid. Then the
initial-boundary-value problem {\rm(\ref{ED+}a-e)}--\eqref{BC}--\eqref{IC} with
{\rm(\ref{split}b,c)}--\eqref{relaxed} has a weak solution
$(\vv,\Ee,\Ep,\alpha,\theta,\chi)$ according
Definition~\ref{def} and every such solution also satisfies
the mechanical-energy balance \eqref{energy+} and
conserves the total energy in the sense \eqref{energy+++}.
\end{theorem}

%\begin{proof}
\noindent{\it Proof.}  For clarity, we divide the proof into five steps.

  \medskip{\it Step 1: Approximation by time discretisation.}
  As we need testing by convective (but not mere partial) time
  derivatives, using of Galerkin method would be very technical (if not
  just impossible). Therefore, 
we use the Rothe method, i.e.\ the fully implicit time discretisation
with an equidistant partition of the time interval $I$ with the time step
$\tau>0$. This approximation is also rather technical because
  the damage-dependent stored energy $\varphi$ is necessarily
  nonconvex. In addition to time discretisation, we thus use also an
  approximation of the stored energy $\varphi$ by a semi-convex
  $\varphi_\eps$ satisfying \eqref{phi-eps}. 

  We denote by $\vvk$, $\Eek$, ${\bm S}_\etau^k$ ...
  the approximate values of $\vv$, $\Ee$, ${\bm S}$ ...
  at time $k\tau$ with $k=1,2,...,T/\tau$.
  We introduce a shorthand notation for the bi-linear operator
    \begin{align}
      \bm B_\text{\sc zj}^{}(\vv,\bm E)=
  (\vv{\cdot}\nabla)\bm E
      -{\rm skew}(\nabla\vv)\Ee+\Ee\,{\rm skew}(\nabla\vv)\,.
    \end{align}
We will then use the following recursive regularized time-discrete scheme
  \begin{subequations}\label{ED+disc}\begin{align}
 \nonumber
 &\varrho \Big(\frac{\vvk{-}\vvkk\!\!\!}\tau+(\vvk{\cdot}\nabla)\vvk
 \Big)
 ={\rm div}{\bm T}_\etau^k-\frac\varrho2({\rm div}\,\vvk)\,\vvk+\ff_\tau^k\,
 \\\nonumber
 &\hspace*{3em}\ \ \ \ \text{ with }\ \
 {\bm T}_\etau^k={\bm S}_\etau^k+\big(\varphi_\eps(\Ee_\etau^k,\alphak){+}\phi(\theta_\etau^k)\big)\bbI
       +{\bm K}_\etau^k\!+\bbD\EE(\vvk)
       -{\rm div}\big(
       \nu|\nabla\EE(\vvk)|^{p-2}\nabla\EE(\vvk)\big)\,,
\\&\hspace*{3em}\ \ \ \ \text{ and }\ \ \ {\bm S}_\etau^k=[\varphi_\eps]_\Ee'(\Ee_\etau^k,\alphak)
\ \ \text{ and }\ \ {\bm K}_\etau^k\!=\kappa\nabla\alphak{\otimes}\nabla\alphak
-\frac{\kappa}2|\nabla\alphak|^2\bbI\,,
\label{ED-1+disc}
\\[-.3em]&\frac{\bm E_\etau^k{-}\bm E_\etau^{k-1}\!\!\!}\tau
+\bm B_\text{\sc zj}^{}(\vvk,\bm E_\etau^k)
=\EE(\vvk)-\RR_\etau^k\,,
 \label{ED-2+disc}
 \\\label{ED-3+disc}
 &
 \GM(\W_\etau^{k-1})\RR_\etau^k={\rm dev}\,\bm{S}_\etau^k+\varkappa\Delta\RR_\etau^k\,,
        \\
        &\partial_{\DT\alpha}\zeta\Big(\alpha_\etau^{k-1},\W_\etau^{k-1};
       \frac{\alphak{-}\alpha_\etau^{k-1}\!\!\!}\tau+
        \vvk{\cdot}\nabla\alphak\Big)
 \label{ED-4+disc}
  +[\varphi_\eps]_\alpha'(\Ee_\etau^k,\alphak)
 +\frac{\alphak{-}\alphakk\!\!\!}{\sqrt\tau}\ni \kappa\Delta\alphak\,,\\
\nonumber
&\frac{\W_\etau^k{-}\W_\etau^{k-1}\!\!\!}\tau
+{\rm div}\big(\vvk\,\W_\etau^k-\bbK(\alpha_\etau^{k-1},\W_\etau^{k-1})\nabla\theta_\etau^k\big)=
\GM(\W_\etau^{k-1})|\RR_\etau^k|^2+\bbD\EE(\vvk){:}\EE(\vvk)
\\[-.6em]\nonumber
&\hspace*{8em}+(1{-}\sqrt[4]{\tau})\xi\Big(\alpha_\etau^{k-1},\W_\etau^{k-1};
\frac{\alphak{-}\alpha_\etau^{k-1}\!\!\!}\tau+
\vvk{\cdot}\nabla\alphak\Big) +\nu|\nabla\EE(\vvk)|^p
\\[-.0em]&\hspace*{8em}
+\varkappa|\nabla\RR_\etau^k|^2
 +\phi(\theta_\etau^k){\rm div}\,\vvk
 \quad\text{ with }\ \quad
\W_\etau^k=\widetilde\gamma(\theta_\etau^k)+\ell\chi_\etau^k\,,\ \ \text{ and}
 \label{ED-6+disc}
 \\[-.3em]
 &\NU\frac{\chi_\etau^k{-}\chi_\etau^{k-1}\!\!\!}\tau
 +\NU\vvk{\cdot}\nabla\chi_\etau^k+H^{-1}(\chi_\etau^k)\ni
\Upsilon\Big(\frac{\theta_\etau^k}{\theta_\text{\sc pt}}-1\Big)
   \label{ED-7+disc}\end{align}\end{subequations}
with $\xi=\xi(\alpha,\W;\DT\alpha)$ from \eqref{ED-6+thermo} and with
$\widetilde\gamma$ from \eqref{split-w}. We complete the system
\eqref{ED+disc} by the corresponding boundary conditions, i.e.
\begin{subequations}\label{BC-disc}
  \begin{align}\label{BC1-disc}
    &\vvk{\cdot}\nn=0\,,\ \ \ \ \big[{\bm T}_\etau^k
\nn+\divS(\nu|\nabla\EE(\vvk)|^{p-2}\nabla\EE(\vvk)\nn)
\big]_\text{\sc t}^{}={\bm0}\,,
  \\\label{BC2}
   &\nabla\EE(\vvk){:}(\nn{\otimes}\nn)={\bm0}\,,\ \ \ 
 (\nn{\cdot}\nabla)\RR_\etau^k={\bm0}\,,\ \ \ 
 \nabla\alphak{\cdot}\nn=0,\ \ \text{ and }\ \
 \\[-.0em]
  &
 \bbK(\alpha_\etau^{k-1},\W_\etau^{k-1}){:}(\nabla\theta_\etau^k\otimes\nn)
 = h_{\rm ext,\tau}^k \,.\label{BCe}
\end{align}
\end{subequations}
Here we used $\ff_\tau^k:=\int_{(k-1)\tau}^{k\tau}\ff(t)\,\d t$ and
similarly also for $\bm{g}_\tau^k$ and $h_{\rm ext,\tau}^k$.
Let us point out that the term $(\alphak{-}\alphakk)/\sqrt\tau$ in
\eqref{ED-4+disc} is devised to convexify $\varphi_\eps$ due to
\eqref{semi-convex} for small $\tau>0$ and vanishes in the limit.
This system of boundary-value problems is to be solved
recursively for $k=1,2,...,T/\tau$, starting 
with the initial conditions for $k=1$: 
\begin{align}\label{IC-disc}
\vv_\etau^0=\vv_0,\ \ \ \ \ {\Ee}_\etau^0={\Ee}_0,\ \ \ \ \ 
\alpha_\etau^0=\alpha_0,\ \ \ \ \ \W_\etau^0=
    \W_0\,,\ \ \text{ and }\ \ \chi_\etau^0=\chi_0
    \,.
\end{align}

The existence of a weak solution $(\vvk,\Eek,\Epk,\alphak,
\chi_\etau^k,\W_\etau^k)\in W^{2,p}(\varOmega;\R^d)\times
L^2(\varOmega;\R_{\rm sym}^{d\times d})
\times H^1(\varOmega;\R_{\rm sym}^{d\times d})\times H^1(\varOmega)\times L^\infty(\varOmega)\times W^{1,1}(\varOmega)$ of the coupled quasi-linear boundary-value problem
(\ref{ED+disc})--(\ref{BC-disc}) can thus be seen by a combination
of the quasilinear technique for \eqref{ED-1+disc} involving the quasilinear
term ${\rm div}^2(\nu|\nabla e(\vvk)|^{p-2}\nabla e(\vvk))$, with the usual
semi-linear technique for (\ref{ED+disc}c--d), with the $L^1$-technique for the
heat equation \eqref{ED-6+disc}, and with the set-valued inclusion (in fact a
variational inequality) \eqref{ED-7+disc} provided $\vvkk$, $\Eekk$,
$\alphakk$, $\W_\etau^{k-1}$, and
$\chi_\etau^{k-1}$ are known from the previous time step. Actually,
$\W_\etau^k\in W^{1,q}(\varOmega)$ for any $q<d'$. Let us note that the system
\eqref{ED+disc} is indeed fully coupled due to the convective derivatives
and due to the adiabatic effects, cf.\ the term  
$\phi(\theta_\etau^k)$ in \eqref{ED-1+disc}, and it seems
not possible to devise some decoupled (staggered) discrete scheme which
would allow for some reasonable estimation strategy of the recursive scheme.
The mentioned coercivity is a particular consequence of the a-priori estimates
derived below. Thus also $\RR_\etau^k\in H^1(\varOmega;\R^{d\times d})$ is obtained.
Let us note that, due to the convective terms,
this system does not have any potential so the rather nonconstructive Brouwer
fixed-point arguments combined with the Galerkin approximation  are to
be used. Here  we rely on  strict
monotonicity of the main parts of  (\ref{ED+disc}a--d) so that the approximated
right-hand side of the semilinear equation \eqref{ED-6+disc}  can be shown
to converge  strongly in $L^1(\varOmega)$.
In general, one cannot expect any uniqueness of this solution.

Since we assume $[\varphi_\eps]_\alpha'(\Ee,0)\le0$ and
$[\varphi_\eps]_\alpha'(\Ee,1)\ge0$, the damage $\alphak$ remains valued in
$[0,1]$ on a.e.\ $\varOmega$ provided $\alpha_0$ is so. Similarly, also
$\chi_\etau^k$ remains valued in $[0,1]$ on a.e.\ $\varOmega$ due to the
constraints involved in $H^{-1}(\cdot)$ in \eqref{ED-7+disc}. Moreover,
$\W_\etau^k\ge0$ a.e.\ on $\varOmega$ for at least one weak solution.
To this goal, we test \eqref{ED-6+disc} by the negative part of $\theta_\etau^k$.
Here we exploit also $\DT\chi\min(0,\theta)=0$ on the discrete level
\eqref{ED-7+disc} when $\theta_\text{\sc pt}>0$ is taken into account,
we prove that $\theta_\etau^k\ge0$. Using also $\chi_\etau^k\ge0$ and
$\widetilde\gamma([0,+\infty))\ge0$,
then also $\W_\etau^k=\widetilde\gamma(\theta_\etau^k)+\ell\chi_\etau^k\ge0$.

Using the values $(\vvk)_{k=0}^{T/\tau}$, we define the piecewise constant and
the piecewise affine interpolants respectively as
\begin{align}\label{def-of-interpolants}
&\overlinevvtau(t):=\vvk,\ \ \ \underline\vv_\etau(t):=\vvkk,
\ \text{ and }\ \vv_\etau(t):=\Big(\frac t\tau{-}k{+}1\Big)\vvk
\!+\Big(k{-}\frac t\tau\Big)\vvkk\ \ \text{ for }\ \
(k{-}1)\tau<t\le k\tau
\end{align}
for $k=0,1,...,T/\tau$. Analogously, we define also 
$\Eetau$, $\overlineEetau$, $\overlineEptau$,  $\underline{\W}_\etau$, etc. 
Thus, \eqref{ED+disc} holding a.e.\ on $\varOmega$ for $k=1,...,T/\tau$
can be written ``compactly'' as 
 \begin{subequations}\label{ED+d}\begin{align}
 \nonumber
 &\varrho
 \Big(\frac{\partial\vv_\etau}{\partial t}
 +(\overlinevvtau{\cdot}\nabla)\overlinevvtau\Big)
 ={\rm div}\overlineTtau
 -\frac\varrho2({\rm div}\,\overlinevvtau)\,\overlinevvtau+\ff_\tau^k\,
 \\\nonumber
  &\hspace*{3em}\ \ \ \ \text{ with }\ \ \overlineTtau=\overlineStau
 +\big(\varphi_\eps(\overlineEetau,\overlinealphatau)
{+}\phi(\overline\theta_\etau)\big)\bbI
               +\bbD\EE(\overlinevvtau)+\overlineSstrtau\!
               -{\rm div}\big(\nu|\nabla e(\overlinevvtau)|^{p-2}\nabla e(\overlinevvtau)\big)\,,
     \\&\hspace*{3em}\ \ \ \ \text{ and }\ \ \
     \overlineStau=[\varphi_\eps]_\Ee'(\overlineEetau,\overlinealphatau)\ \ \text{ and }\ \ 
\overlineSstrtau\!=\kappa\nabla\overlinealphatau{\otimes}\nabla\overlinealphatau-\frac{\kappa}2|\nabla\overlinealphatau|^2\bbI\,,
\label{ED-1+d}
\\[-.3em]&
\frac{\partial\Ee_\etau}{\partial t}
+\bm B_\text{\sc zj}^{}(\overlinevvtau,\overlineEetau)
=\EE(\overlinevvtau)-\overlineRRtau\,,
 \label{ED-2+d}
 \\\label{R}
 &
\GM(\underline\W_\etau)\overlineRRtau={\rm dev}\,\overlineStau+\varkappa\Delta\overlineRRtau\,,
        \\
        &\partial_{\DT\alpha}\zeta\Big(\underline\alpha_\etau,\underline\W_\etau;
 \frac{\partial\alpha_\etau}{\partial t}+
        \overlinevvtau{\cdot}\nabla\overlinealphatau\Big)
 \label{ED-4+d}
 +[\varphi_\eps]_\alpha'(\overlineEetau,\overlinealphatau)
 +
 \sqrt\tau\frac{\partial\alpha_\etau}{\partial t}
 \ni \kappa\Delta\overlinealphatau\,,
\\
\nonumber
&\frac{\partial\W_\etau}{\partial t}
+{\rm div}\big(\overlinevvtau\,\overline\W_\etau-\bbK(\underline\alpha_\etau,\underline\W_\etau)\nabla\overline\theta_\etau\big)=
\GM(\underline\W_\etau)|\overlineRRtau|^2+
\bbD\EE(\overlinevvtau){:}\EE(\overlinevvtau)
\\[-.3em]\nonumber
&\hspace*{6em}+
(1{-}\sqrt[4]{\tau})\xi\Big(\underline\alpha_\etau,\underline\W_\etau;
\frac{\partial\alpha_\etau}{\partial t}
+\overlinevvtau{\cdot}\nabla\overlinealphatau\Big) \\[-.3em]
&\hspace*{6em}+\varkappa|\nabla\overlineRRtau|^2+
\nu|\nabla\EE(\overlinevvtau)|^p
+\phi(\overline\theta_\etau){\rm div}\,\overlinevvtau
 \quad
\text{ with }\ 
\quad
\overline\W_\etau=\widetilde\gamma(\overline\theta_\etau)+\ell\overlinechitau\,,\ \ \text{ and}
 \label{ED-6+d}
 \\[-.3em]
 &\NU\frac{\partial\chi_\etau}{\partial t}
 +\NU\overlinevvtau{\cdot}\nabla\overlinechitau
 +H^{-1}(\overlinechitau)\ni
\Upsilon\Big(\frac{\overline\theta_\etau}{\theta_\text{\sc pt}}-1\Big)
   \label{ED-7+d}\end{align}\end{subequations}
 holding  on $I{\times}\varOmega$ either a.e.\ or in a weak sense
 involving also the boundary conditions \eqref{BC-disc} which are
 to be written analogously in terms of the above introduced interpolants.

\medskip{\it Step 2: A-priori estimates}.
The a-priori estimation is based on the energy test for the mechanical
part combined with the heat problem. This means here 
the test of \eqref{ED-1+disc} by $\vvk$ while using also \eqref{ED-2+disc}
tested by $\bm{S}_\etau^k=[\varphi_\eps]_\Ee'(\Ee_\etau^k,\alphak)$, then the
test \eqref{ED-3+disc} by $\RR_\etau^k$ and the inclusion \eqref{ED-4+disc} by
$(\alphak{-}\alphakk)/\tau+\vvk{\cdot}\nabla\alphak$, and \eqref{ED-6+disc}
by 1. We thus obtain an energy-like inequality for the time-discrete
approximation, but in contrast to \eqref{energy+++}, we still keep part of
the dissipation rate like in \eqref{energy+}.

More specifically, the terms related to inertia in \eqref{ED-1+disc} uses
the calculus 
\begin{align}\nonumber
  \Big(\varrho\frac{\vvk{-}\vvkk\!\!}\tau
  +\varrho(\vvk{\cdot}\nabla)\vvk+\frac\varrho2(\operatorname{div}\vvk)\vvk\Big){\cdot}\vvk&=
\frac{\varrho}2\frac{|\vvk|^2-|\vvkk|^2\!\!}\tau
  +\varrho(\vvk{\cdot}\nabla)\vvk\cdot\vvk
  \\[-.4em]&\quad\ +\frac{\varrho}2({\rm div}\,\vvk)|\vvk|^2
  +\tau\frac{\varrho}2\Big|\frac{\vvk{-}\vvkk\!\!}\tau\ \Big|^2\,.
\label{test-of-convective}\end{align}
This holds pointwise and, when integrated over $\varOmega$, we further use
also \begin{align}\nonumber
 \INT{\varOmega}{}\varrho(\vvk{\cdot}\nabla)\vvk{\cdot}\vvk\,\d x&=
  -\INT{\varOmega}{}\frac{\varrho}2|\vvk|^2({\rm div}\,\vvk)\,\d x
  +\INT{\varGamma}{}\frac{\varrho}2|\vvk|^2(\vvk{\cdot}\nn)\,\d S
\\&=
-\INT{\varOmega}{}\Big(\frac{\varrho}2|\vvk|^2\bbI\Big):\EE(\vvk)\,\d x
  +\INT{\varGamma}{}\frac{\varrho}2|\vvk|^2(\vvk{\cdot}\nn)\,\d S\,.
\label{convective-tested+}\end{align}
The last term in \eqref{test-of-convective} is non-negative and will simply
be forgotten, which will give the inequality
\begin{align}\nonumber
  &\INT{\varOmega}{}\!\Big(\varrho\frac{\vvk{-}\vvkk\!\!}\tau
  +\varrho(\vvk{\cdot}\nabla)\vvk
  +\frac\varrho2({\rm div}\,\vvk)\vvk\Big)\cdot\vvk\,\d x
  \stackrel{\eqref{test-of-convective}}{\ge}\INT{\varOmega}{}\!\Big(
\frac{\varrho}2\frac{|\vvk|^2-|\vvkk|^2\!\!}\tau
   \\&\ +
  \varrho(\vvk{\cdot}\nabla)\vvk{\cdot}\vvk
  +\frac\varrho2({\rm div}\,\vvk)|\vvk|^2\Big)\,\d x
  \stackrel{\eqref{convective-tested+}}{=}
  \INT{\varOmega}{}\frac{\varrho}2\frac{|\vvk|^2\!-|\vvkk|^2\!\!}\tau\,\d x
  +\INT{\varGamma}{}\frac\varrho2|\vvk|^2(\vvk{\cdot}\nn)\,\d S\,.
\label{test-of-convective+}\end{align}
The last term vanishes due to the boundary condition \eqref{BC1-disc}.
The further term in \eqref{ED-1+disc} uses the calculus
\begin{align}\nonumber
  &\INT{\varOmega}{}\operatorname{div}{\bm S}_\etau^k
  \cdot\vvk\,\d x=\INT{\varGamma}{}{\bm S}_\etau^k{:}(\vvk\otimes\nn)\,\d S
  -\INT{\varOmega}{}{\bm S}_\etau^k{:}\EE(\vvk)\,\d x
\\&\nonumber\stackrel{\eqref{ED-2+disc}}{=}
\INT{\varGamma}{}{\bm S}_\etau^k{:}(\vvk\otimes\nn)\,\d S
-    \INT{\varOmega}{}
{\bm S}_\etau^k{:}\Big(\frac{\Ee_\etau^k{-}\Ee_\etau^{k-1}}\tau+
\bm B_\text{\sc zj}^{}(\vvk,\Ee_\etau^k)-\RR_\etau^k
\Big)\,\d x
\\&\nonumber\!\!\stackrel{(\ref{ED+disc}a,c)}{=}
\!\!\INT{\varGamma}{}{\bm S}_\etau^k{:}(\vvk\otimes\nn)\,\d S
+\INT{\varOmega}{}\bigg( 
\GM(\W_\etau^{k-1})|\RR_\etau^k|^2+\varkappa|\nabla\RR_\etau^k|^2
\\[-.4em]&
\hspace*{15em}
-[\varphi_\eps]_\Ee'(\Eek,\alphak){:}\Big(\frac{\Eek{-}\Eekk}\tau
    {+}\bm B_\text{\sc zj}^{}(\vvk,\Eek)\Big)\bigg)\,\d x\,,
\label{test-by-v}\end{align}
where we used also \eqref{ED-3+disc} tested by $\RR_\etau^k$.

Testing \eqref{ED-4+disc} by the discrete convective derivative
of $\alpha$, we obtain 
\begin{align}\nonumber
\INT{\varOmega}{}&\bigg(\xi\Big(\alpha_\etau^{k-1},\W_\etau^{k-1};
       \frac{\alphak{-}\alpha_\etau^{k-1}\!\!\!}\tau+
        \vvk{\cdot}\nabla\alphak\Big)
  +\Big([\varphi_\eps]_\alpha'(\Ee_\etau^k,\alphak)+\frac{\alphak{-}\alphakk\!\!\!}{\sqrt\tau}\ \,\Big)
 \,\frac{\alphak{-}\alpha_\etau^{k-1}\!\!\!}\tau
 \\\label{ED-4+disc-tested}
 &\hspace*{7em}
 +\frac{\kappa}{2\tau}|\nabla\alphak|^2\bigg)\,\d x
\le\INT{\varOmega}{}\frac{\kappa}{2\tau}|\nabla\alphakk|^2
 -\Big(\Delta\alphak+\frac{\alphak{-}\alphakk}{\sqrt\tau\!\!\!}\Big)\,
 \vvk{\cdot}\nabla\alphak\,.
\end{align}

For summing the $\varphi$-terms in \eqref{test-by-v} and
\eqref{ED-4+disc-tested}, we use the semi-convexity \eqref{semi-convex}
of $\varphi_\eps$ to estimate
\begin{align}\nonumber
  &{\bm S}_\tau^k{:}\frac{\Eek{-}\Eekk\!\!}\tau
+\Big([\varphi_\eps]_\alpha'(\Eek,\alphak)
+\frac{\alphak{-}\alphakk\!\!}{\sqrt\tau}\ \Big)\,\frac{\alphak{-}\alphakk\!\!}\tau  
\\&\nonumber\hspace*{0em}
=[\varphi_\eps]_\Ee'(\Eek,\alphak){:}\frac{\Eek{-}\Eekk\!\!}\tau+
  \Big([\varphi_\eps]_\alpha'(\Eek,\alphak)
  +\frac{\alphak}{\sqrt\tau}\Big)\,\frac{\alphak{-}\alphakk\!\!}\tau
-\frac{\alphakk}{\sqrt\tau}\,\frac{\alphak{-}\alphakk\!\!}\tau
  \\&\nonumber\hspace*{0em}
  \ge
  \frac{\varphi_\eps(\Eek,\alphak)
  -\varphi_\eps(\Eekk,\alphakk)}\tau+\frac1{2\sqrt\tau}\frac{|\alphak|^2{-}|\alphakk|^2}\tau
-\frac{\alphakk}{\sqrt\tau}\,\frac{\alphak{-}\alphakk\!\!}\tau
  \\&\hspace*{0em}
  =\frac{\varphi_\eps(\Eek,\alphak)-\varphi_\eps(\Eekk,\alphakk)}\tau
-\frac{\sqrt\tau}2\Big|\frac{\alphak{-}\alphakk\!\!}\tau\ \Big|^2
\,,\label{semi-convex-trick}
\end{align}
cf.\ also the calculation in \cite[Remark~8.24]{Roub13NPDE}. This holds
a.e.\ on $\varOmega$ and is to be integrated over $\varOmega$.
For the remaining two  convective terms arising from these tests, we use
the calculus
\begin{align}\nonumber&
\INT{\varOmega}{}\Big([\varphi_\eps]_{\Ee}'(\Eek,\alphak){:}(\vvk{\cdot}\nabla)\Eek
+[\varphi_\eps]_\alpha'(\Eek,\alphak){\cdot}(\vvk{\cdot}\nabla\alphak)\Big)\,\d x
=\INT{\varOmega}{}
\nabla\varphi_\eps(\Eek,\alphak){\cdot}\vvk\,\d x
\\[-.3em]&
\qquad
=\INT{\varGamma}{}\!\varphi_\eps(\Eek,\alphak)\vvk\cdot{\bm n}\,\d S-\INT{\varOmega}{}\!
\varphi_\eps(\Eek,\alphak){\rm div}\,\vvk\,\d x
=-\INT{\varGamma}{}\varphi_\eps(\Eek,\alphak)\bbI{:}\EE(\vvk)\,\d x\,,
\label{calc-energy-pressure}\end{align}
which cancels with the pressure-type stress contribution
$\varphi_\eps(\Eek,\alphak)\bbI$, cf.\ \eqref{calculus-toward-pressure}.
Noteworthy, by the test of the regularizing term
$(\alphak{-}\alphakk)/\sqrt\tau$ by the convective time difference
$(\alphak{-}\alphakk)/\tau+\vvk{\cdot}\nabla\alphak$,
as in \cite{RouTom21CMPE}, we obtain still the term
\begin{align}\nonumber
\INT{\varOmega}{}\frac{\alphak{-}\alphakk\!\!\!}{\sqrt\tau}\ \vvk{\cdot}\nabla\alphak\,\d x
&=
\sqrt\tau\INT{\varOmega}{}\frac{\alphak{-}\alphakk\!\!\!}{\tau}\ \vvk{\cdot}\nabla\alphak\,\d x
\\\nonumber&
=\sqrt\tau\INT{\varOmega}{}\frac{\alphak{-}\alphakk}{\tau}
\Big(\frac{\alphak{-}\alphakk\!\!}\tau+\vvk{\cdot}\nabla\alphak\Big)
-\Big|\frac{\alphak{-}\alphakk\!\!\!}{\tau}\ \Big|^2\,\d x
\\&\le\frac{\sqrt\tau}2\Big\|\frac{\alphak{-}\alphakk\!\!}\tau
+\vvk{\cdot}\nabla\alphak\Big\|_{L^2(\varOmega)}^2\!
-\frac{\sqrt\tau}2\Big\|\frac{\alphak{-}\alphakk\!\!\!}{\tau}\ \Big\|_{L^2(\varOmega)}^2\,.
\label{semi-convex-trick+}\end{align}
The last term advantageously absorbs the last term in
\eqref{semi-convex-trick} while the penultimate term in
\eqref{semi-convex-trick+} can be absorbed for sufficiently small $\tau$'s
in the dissipation term $\xi(\alphak,\W_\etau^k;\cdot)$ later in
\eqref{mech-energy-disc}.

The other stress contributions, i.e.\ ${\bm K}_\etau^k\!+\bbD\EE(\vvk)
-{\rm div}(\nu|\nabla\EE(\vvk)|^{p-2}\nabla\EE(\vvk))$, can be treated by
adapting the calculus \eqref{test-1} and \eqref{test-damage}. 

Altogether, after summation over $k=1,2,....$, we obtain \eqref{energy}
as an upper estimate up
to an error term which is small for $\tau>0$ small, so that it can be used
for a-priori estimates. More specifically, we thus obtain the discrete
mechanical-energy balance as an inequality
\begin{align}\nonumber
  &
  \INT{\varOmega}{}\frac\varrho2|\vv_\etau^k|^2+\varphi_\eps(\Ee_\etau^k,\alphak)
    +\frac{\kappa}2|\nabla\alphak|^2
    \,\d x+\tau\!\INT{\varOmega}{}\!\bigg(
    \GM(\W_\etau^{k-1})|\RR_\etau^k|^2+
\bbD\EE(\vvk){:}\EE(\vvk)
\\[-.2em]\nonumber
&\qquad\qquad
+\nu|\nabla\EE(\vvk)|^p
+\xi\Big(\alpha_\etau^{k-1},\W_\etau^{k-1};
 \frac{\alphak{-}\alpha_\etau^{k-1}\!\!}\tau+\vvk{\cdot}\nabla\alphak\Big)
 +\varkappa|\nabla\RR_\etau^k|^2\bigg)\,\d x
\\&\hspace*{2em}\nonumber
  \le\INT{\varOmega}{}\frac\varrho2|\vv_\etau^{k-1}|^2+\varphi_\eps(\Ee_\etau^{k-1},
    \alpha_\etau^{k-1})+\frac{\kappa}2|\nabla\alpha_\etau^{k-1}|^2
  +\tau\ff_\tau^k{\cdot}\vv_\etau^k
 -\tau\phi(\theta_\etau^k){\rm div}\vv_\etau^k\,\d x
 \\[-.3em]
 &\hspace*{21em}
+\frac{\sqrt\tau}2\Big\|\frac{\alphak{-}\alphakk\!\!}\tau
+\vvk{\cdot}\nabla\alphak\Big\|_{L^2(\varOmega)}^2 \,.
\label{mech-energy-disc}\end{align}

Then we test the discrete heat-transfer equation \eqref{ED-6+disc} by 1 and
add \eqref{mech-energy-disc}, we can cancel the adiabatic terms
$\pm\tau\phi(\theta_\etau^k){\rm div}\vv_\etau^k$. Summing it for $k=1,...,l$
with $l\le T/\tau$, we obtain a discrete energy-like inequality:
\begin{align}\nonumber
\!\!\INT{\varOmega}{}\!\frac\varrho2|\vv_\etau^l|^2\!+\varphi_\eps(\Ee_\etau^l,\alphal)
  +\frac{\kappa}2|\nabla\alphal|^2\!+\W_\etau^l\,\d x
  &\le\INT{\varOmega}{}\frac\varrho2|\vv_0|^2\!+\varphi_\eps(\Ee_0,\alpha_0)
+\frac{\kappa}2|\nabla\alpha_0|^2\!
  +\W_0\,\d x
 \\[-.4em]
 &\hspace*{3.5em}
 +\!\tau\sum_{k=1}^l\bigg(\INT{\varOmega}{}\!\ff_\tau^k{\cdot}\vv_\etau^k
 \,\d x
 +\!\INT{\varGamma}{} h_{\rm ext,\tau}^k\d S\!\bigg)\!
\label{energy-disc}\end{align}
provided $\tau\le2^8\epsilon^4$ with $\epsilon$ from \eqref{ass:2}; here the
smallness of $\tau$ allows for absorption of the last term in
\eqref{mech-energy-disc} when using the uniform convexity of the rest of the
dissipation potential $\sqrt[4]{\tau}\zeta(\alpha,\W;\cdot)$ which also
implies the uniform convexity of the corresponding dissipation rate
$\sqrt[4]{\tau}\xi(\alpha,\W;\cdot)$. This reveals why we used the
factor $1{-}\sqrt[4]{\tau}$ in \eqref{ED-6+disc}.

Then, by the discrete Gronwall inequality, 
we obtain the a-priori estimates
\begin{subequations}\label{est-disc}\begin{align}\label{est-disc1}
    &\|\vv_\etau^{}\|_{L^\infty(I;L^2(\varOmega;\R^d))}^{}\le C,
  \\\label{est-disc5}
&\|{\bm S}_\etau^{}\|_{L^\infty(I;L^2(\varOmega;\R_{\rm sym}^{d\times d}))}^{}\le C\,,
\\& \|\RR_\etau^{}\|_{L^2(I;H^1(\varOmega;\R_{\rm dev}^{d\times d}))}^{}\le C\,.
\label{est-R}
\\\label{est-disc4}
   &\|\alpha_\etau^{}\|_{L^\infty(I;H^1(\varOmega))
     \,\cap\,L^\infty(I{\times}\varOmega)}^{}\le C,
   \\\label{est-disc-w}
   &\|\W_\etau^{}\|_{L^\infty(I;L^1(\varOmega))}^{}\le C\ \ \text{ and }\ \
   \|\theta_\etau^{}\|_{L^\infty(I;L^1(\varOmega))}^{}\le C\,.
\end{align}\end{subequations}
Actually, from \eqref{energy-disc} we can see
$\|\varphi_\eps(\overlineEetau,\overlinealphatau)\|_{L^\infty(I;L^1(\varOmega))}\le C$  and then from \eqref{ass:1-eps} we obtain the estimate $\|[\varphi_\eps]_\Ee'(\overlineEetau,\overlinealphatau)\|_{L^\infty(I;L^2(\varOmega;\R^{d\times d}))}\le C$, i.e.\ \eqref{est-disc5}.

Furthermore, we test the volume-fraction equation \eqref{ED-7+d} by
the convective derivative of $\chi$, i.e.\ in the discrete form by
$\frac{\partial}{\partial t}\chi_\etau
+\overlinevvtau{\cdot}\nabla\overlinechitau$.
We obtain
\begin{align}\nonumber
\INT{\varOmega}{}\!
  \NU\Big|\frac{\partial\chi_\etau}{\partial t}
  +\overlinevvtau{\cdot}\nabla\overlinechitau\Big|^2\!\!
  +H^{-1}(\overlinechitau)\frac{\partial\chi_\etau}{\partial t}
  +H^{-1}(\overlinechitau)\overlinevvtau{\cdot}\nabla\overlinechitau\,\d x\
  \\[-.3em]
  =\INT{\varOmega}{}\!\Upsilon\Big(\,\frac{\overline\theta_\etau}{\theta_\text{\sc pt}}-1\Big)
  \Big(\frac{\partial\chi_\etau}{\partial t}
  +\overlinevvtau{\cdot}\nabla\overlinechitau\Big)\,\d x\,.
  \label{test-dchi-dt}\end{align}
Realizing that
$H^{-1}(\overlinechitau)\frac{\partial}{\partial t}\chi_\etau
=\frac{\partial}{\partial t}\delta_{[0,1]}(\chi_\etau)=0$
and
$H^{-1}(\overlinechitau)\overlinevvtau{\cdot}\nabla\overlinechitau
=\overlinevvtau{\cdot}\nabla\delta_{[0,1]}(\overlinechitau)=0$
because $\overlinechitau$ is valued in $[0,1]$ so that
$\delta_{[0,1]}(\overlinechitau)=0$,
we obtain that
$\frac{\partial}{\partial t}\chi_\etau{+}\overlinevvtau{\cdot}\nabla\overlinechitau$
bounded in $L^2(I{\times}\varOmega)$. Here we also used that
$\Upsilon(\overline\theta_\etau/\theta_\text{\sc pt}{-}1)$ is surely bounded in
$L^2(I{\times}\varOmega)$ due to the already proved bound for
$\overline\theta_\etau\in L^\infty(I;L^1(\varOmega))$ and due to the growth assumption
\eqref{Y-ass}.

As the next step, we use the $L^1$-technique to estimate of temperature
gradient developed by Boccardo and Gallou{\"e}t \cite{BocGal89NEPE}
exploiting sophistically Gagliardo-Nirenberg inequality, cf.\
\cite[Prop.\,8.2.1]{KruRou19MMCM}. The essence is to test the heat-transfer
equation \eqref{ED-6+d} by a smoothened the Heaviside function $H$ from
\eqref{split-H}, say $\omega(\theta):=1-(1{+}\theta)^{-\epsilon}$
for $\epsilon>0$, as suggested in \cite{FeiMal06NSET}. The modification in
comparison with the
usual ``heat operator'' in the form $\frac{\partial}{\partial t}
\vartheta-{\rm div}(\bbK\nabla\theta)$ with $\vartheta=\widetilde\gamma(\theta)$
and with an $L^1$-right-hand side
consists in that $\DT\W$ contains two other terms more, namely
$\DT\W=\frac{\partial}{\partial t}\widetilde\gamma(\theta)
+\ell\DT\chi+{\rm div}(\widetilde\gamma(\theta)\vv)$. In the
discrete form, these additional terms are
$\ell(\frac{\partial}{\partial t}\chi_\etau
{+}\overlinevvtau{\cdot}\nabla\overlinechitau)
+{\rm div}(\widetilde\gamma(\overline\theta_\etau)\overlinevvtau)$. As for
$\ell(\frac{\partial}{\partial t}\chi_\etau{+}
\overlinevvtau{\cdot}\nabla\overlinechitau)$, we have it already
estimated in $L^2(I{\times}\varOmega)$ and it remains so after being tested by
$\omega(\overline\theta_\etau)$. As for the convective term
${\rm div}(\widetilde\gamma(\overline\theta_\etau)\overlinevvtau)$ tested by
$\omega(\overline\theta_\etau)$, we can estimate it ``on the right-hand side''
as
\begin{align}\nonumber
 -\INT{\varOmega}{}{\rm div}(\widetilde\gamma(\overline\theta_\etau)\overlinevvtau)
  \omega(\overline\theta_\etau)\,\d x
  &=\INT{\varOmega}{}\overlinevvtau{\cdot}\widetilde\gamma(\overline\theta_\etau)
  \omega'(\overline\theta_\etau)\nabla\overline\theta_\etau\,\d x
  \\&\le\INT{\varOmega}{}\frac1\epsilon|\overlinevvtau|^2
  \widetilde\gamma(\overline\theta_\etau)^2\omega'(\overline\theta_\etau)\,\d x
  +\epsilon\INT{\varOmega}{}\omega'(\overline\theta_\etau)|\nabla\overline\theta_\etau|^2\,\d x\,.
\end{align}
As $\widetilde\gamma(\theta)=\mathscr{O}(\theta)$ while
$\omega'(\theta)=\mathscr{O}(1/\theta)$, we have
$[\widetilde\gamma^2\omega'](\theta)=\mathscr{O}(\theta)$ so that
$\widetilde\gamma(\overline\theta_\etau)^2\omega'(\overline\theta_\etau)$ is
bounded in $L^\infty(I;L^1(\varOmega))$ while $|\overlinevvtau|^2$
is surely bounded in $L_{\rm w*}^{p/2}(I;L^\infty(\varOmega))$, so that 
the integrand in the penultimate integral is bounded in $L^{p/2}(I;L^1(\varOmega))$.
For the last integral, this is exactly fitted with the estimation in the
$L^1$-theory, and for $\epsilon>0$ sufficiently small can be absorbed in
the respective estimation, cf.\ \cite[Formula (8.2.17)]{KruRou19MMCM}.
Using also the already obtained estimate \eqref{est-disc-w}, the resulted
``prefabricated'' estimate is that $\|\theta_\etau\|_{L^s(I{\times}\varOmega)}^s$
and $\|\nabla\theta_\etau\|_{L^r(I{\times}\varOmega;\R^d)}^r$ is
estimated by the $L^1(I{\times}\varOmega)$-norm of the right-hand side of
\eqref{ED-6+d} provided $s<1+2/d$ and $r<(d{+}2)/(d{+}1)$,
cf.\ \cite[Formula~(8.2.6b)]{KruRou19MMCM}.

Then, to estimate the dissipation rate, we come back to \eqref{mech-energy-disc}
and estimate $|\int_{\varOmega}{}\phi(\theta_\etau^k){\rm div}\vv_\etau^k\,\d x|
\le C\|\theta_\etau^k\|_{L^1(\varOmega)}\|\vv_\etau^k\|_{W^{2,p}(\varOmega;\R^d)}$, exploiting 
$p>d$ and \eqref{est-disc-w}. Actually, we have now even a better estimate of
$\theta_\etau^k$ than $L^1(\varOmega)$ but we do not need it.

To summarize, in addition to \eqref{est-disc}, we have proved
\begin{subequations}\label{est-disc+}\begin{align}
\label{est-disc+1}
    &\|\vv_\etau^{}\|_{L^p(I;W^{2,p}(\varOmega;\R^d))}^{}\le C,
   \\\label{est-disc+4}
   &\|\alpha_\etau^{}\|_{H^1(I;L^2(\varOmega))}^{}\le C,
\\\label{est-disc+theta}
&\|\overline\theta_\etau\|_{L^r(I;W^{1,r}(\varOmega))}^{}\le C_r\ \ \ \ \ \ \ \ 
\text{ with $\ 1\le r<\frac{d+2}{d+1}$}\,,\ \text{ and} 
\\[-.4em]\label{est-disc+vartheta}
&\|\overline\vartheta_\etau\|_{L^r(I;W^{1,r}(\varOmega))
  \,\cap\,L^s(I{\times}\varOmega)}^{}\le C_{r,s}
\ \ \ \ \text{ with $\ \ 1\le s<1+\frac2d$.}
\end{align}\end{subequations}
The $L^r$-estimate \eqref{est-disc+vartheta} of $\nabla\overline\vartheta_\etau$ 
can be read from \eqref{est-disc+theta}
due to $\nabla\overline\vartheta_\etau=\nabla\widetilde\gamma(\overline\theta_\etau)
=\widetilde\gamma'(\overline\theta_\etau)\nabla\overline\theta_\etau$.
The $L^s$-estimate of $\overline\vartheta_\etau$ 
is to be read by the Gagliardo-Nirenberg interpolation
of the first estimate in \eqref{est-disc+vartheta} with  \eqref{est-disc-w}.
For \eqref{est-disc+4}, we use the bound for
$\frac{\partial}{\partial t}\alpha_\etau+\overlinevvtau{\cdot}\nabla\overlinealphatau$ in $L^2(I{\times}\varOmega)$ obtained from \eqref{energy-disc} and that
$\overlinevvtau{\cdot}\nabla\overlinealphatau$ is bounded in
$L^2(I{\times}\varOmega)$ due to \eqref{est-disc1} with the former part of
\eqref{est-disc+4}.
By comparison from $\Delta\overlinealphatau\in
(\partial_{\DT\alpha}\zeta(\underline\alpha_\etau,\underline\W_\etau;
\frac{\partial}{\partial t}\alpha_\etau{+}\overlinevvtau{\cdot}\nabla\overlinealphatau)+[\varphi_\eps]_\alpha'(\overlineEetau,\overlinealphatau)
+\tau\frac{\partial}{\partial t}\alpha_\etau)/\kappa$, we also obtain
\begin{align}\label{est-Delta}
\|\Delta\overlinealphatau\|_{L^2(I{\times}\varOmega)}^{}\le C\,.
\end{align}

\def\rexp{\sigma}
An important attribute of the model is that the convective transport of
variables via the velocity field
$\vv\in L_{\rm w*}^1(I;W^{1,\infty}(\varOmega;\R^d))$ or, in the discrete variant
by $\overline\vv_\tau$ bounded in $L^1(I;W^{1,\infty}(\varOmega;\R^d))$,
qualitatively well copies  regularity  properties of the initial conditions.
We use this phenomenon particularly for the Zaremba-Jaumann time difference
 of the elastic strain $\Ee=[E_{ij}]$ and also for
the convective time-difference of the volume fraction $\chi$.
Let us consider a general tensor-valued source ${\bm F}\in
L^2(I{\times}\varOmega;\R^{d\times d})$ for $\ZJ\Ee={\bm F}$, i.e.\ in the difference variant
\begin{align}
\frac{{\bm E}_\tau^k-{\bm E}_\tau^{k-1}\!\!}\tau
  +B_\text{\sc zj}(\vv_\tau^k,{\bm E}_\tau^k)={\bm F}_\tau^k\,.
  \label{ZJ-evol-abstract}\end{align}
For $\rexp>1$, we use the following calculus exploiting the Green formula with
the boundary condition $\vv{\cdot}\nn=0$: 
\begin{align}
  \INT{\varOmega}{}(\vv{\cdot}\nabla z)|z|^{\rexp-2}z\,\d x
  &=\INT{\varGamma}{}\!|z|^\rexp(\vv{\cdot}\nn)\,\d S
  -\INT{\varOmega}{}\!(\rexp{-}1)|z|^{\rexp-2}z(\vv{\cdot}\nabla z)+({\rm div}\,\vv)|z|^\rexp\d x
=-\frac1\rexp\INT{\varOmega}{}({\rm div}\,\vv)|z|^\rexp\d x\,.
\label{transport}\end{align}
We test \eqref{ZJ-evol-abstract}
by $|{\bm E}_\tau^k|^{\rexp-2}{\bm E}_\tau^k$, which gives
\begin{align}
  \frac1\rexp\INT{\varOmega}{}\!\frac{|{\bm E}_\tau^k|^\rexp\!-|{\bm E}_\tau^{k-1}|^\rexp\!}\tau
  \,\d x\le\INT{\varOmega}{}\!\frac{{\rm div}\,\vv_\tau^k}\rexp|{\bm E}_\tau^k|^\rexp
  +2|{\rm skew}(\nabla\vv_\tau^k)|\,|{\bm E}_\tau^k|^\rexp
  +|{\bm E}_\tau^k|^{\rexp-2}{\bm E}_\tau^k{:}{\bm F}_\tau^k\,\d x\,,
\label{transport+}\end{align}
where we used \eqref{transport} for each component $z=E_{ij}$. From
\eqref{transport+}, by the Young and the discrete Gronwall inequalities,
we obtain the estimate
\begin{align}
  \|{\bm E}_\tau^k\|_{L^\rexp(\varOmega;\R^{d\times d})}^\rexp\le
  C{\rm e}^{1+2k\tau\max_{l=1,...,k}\|\nabla\vv_\tau^k\|_{L^\infty(\varOmega;\R^{d\times d})}}\bigg(\|{\bm E}_\tau^0\|_{L^\rexp(\varOmega;\R^{d\times d})}^\rexp\!+
  \tau\sum_{l=1}^k\|{\bm F}_\tau^l\|_{L^{\rexp'}(\varOmega;\R^{d\times d})}^{\rexp'}\bigg)\,
\label{Gronwall}\end{align}
for some $C$ and for $\tau>0$ sufficiently small, namely for all $\tau\le1/
((2{+}4\rexp)\max_{l=1,...,k}\|\nabla\vv_\tau^k\|_{L^\infty(\varOmega;\R^{d\times d})}+2\rexp)$.

Moreover, we can also test \eqref{ZJ-evol-abstract} by the $q$-Laplacian
$-{\rm div}(|\nabla{\bm E}_\tau^k|^{q-2}\nabla{\bm E}_\tau^k)$. More specifically,
we can apply the $\nabla$-operator to \eqref{ZJ-evol-abstract}
and test it by $|\nabla{\bm E}_\tau^k|^{q-2}\nabla{\bm E}_\tau^k$.
Instead of \eqref{transport}, we use the calculus 
\begin{align}\nonumber
&\INT{\varOmega}{}\nabla(\vv{\cdot}\nabla z){\cdot}|\nabla z|^{q-2}\nabla z\,\d x
=\INT{\varOmega}{}|\nabla z|^{q-2}\nabla\vv{:}(\nabla z{\otimes}\nabla z)
+(\vv{\cdot}\nabla^2z){\cdot}|\nabla z|^{q-2}\nabla z\,\d x
\\[-.4em]&\nonumber\quad=\INT{\varGamma}{}\!|\nabla z|^q(\vv{\cdot}\nn)\,\d S
-\INT{\varOmega}{}\!|\nabla z|^{q-2}\nabla\vv{:}(\nabla z{\otimes}\nabla z)
+(q{-}1)|\nabla z|^{q-2}\nabla z{\cdot}(\vv{\cdot}\nabla^2z)
+({\rm div}\,\vv)|\nabla z|^q\,\d x
\\[-.4em]&\hspace{20.5em}
=\INT{\varOmega}{}|\nabla z|^{q-2}\nabla\vv{:}(\nabla z{\otimes}\nabla z)
-\frac1q({\rm div}\,\vv)|\nabla z|^q\,\d x\,.
\label{transport-}\end{align}
In the tensorial situation \eqref{ZJ-evol-abstract}, we use it again for
$z=E_{ij}$ and then we use also 
\begin{align}
\nabla\Big({\rm skew}(\nabla\vv_\tau^k){\bm E}_\tau^k
-{\bm E}_\tau^k{\rm skew}(\nabla\vv_\tau^k)\Big)
\Vdots|\nabla{\bm E}_\tau^k|^{q-2}\nabla{\bm E}_\tau^k\le
2|\nabla\vv_\tau^k|\,|\nabla{\bm E}_\tau^k|^q
+2|\nabla^2\vv_\tau^k|\,|{\bm E}_\tau^k|\,|\nabla{\bm E}_\tau^k|^{q-1}.
\label{transport--}\end{align}
Instead of \eqref{transport+}, this gives
\begin{align}\nonumber
\frac1q\INT{\varOmega}{}
\frac{|\nabla {\bm E}_\tau^k|^q\!-|\nabla {\bm E}_\tau^{k-1}|^q\!}\tau
\,\d x&\le\INT{\varOmega}{}\bigg(\frac{{\rm div}\,\vv_\tau^k}q|\nabla{\bm E}_\tau^k|^q
+\Big(2+\frac1q\Big)|\nabla\vv_\tau^k|
\,|\nabla{\bm E}_\tau^k|^q
\\[-.7em]\nonumber
&\qquad
+|\nabla{\bm E}_\tau^k|^{q-2}\nabla{\bm E}_\tau^k{\Vdots}\nabla{\bm F}_\tau^k
+2|\nabla^2\vv_\tau^k|\,|{\bm E}_\tau^k|\,|\nabla{\bm E}_\tau^k|^{q-1}
\bigg)\,\d x
\\[-.3em]&\nonumber\le
C+C\big(1{+}\|\nabla\vv_\tau^k\|_{L^\infty(\varOmega;\R^{d\times d})}^{}\big)
\|\nabla{\bm E}_\tau^k\|_{L^q(\varOmega;\R^{d\times d\times d})}^q
\\[-.1em]
&\qquad
+\|\nabla{\bm F}_\tau^k\|_{L^q(\varOmega;\R^{d\times d\times d})}^q+
\|\nabla^2\vv_\tau^k\|_{L^p(\varOmega;\R^{d\times d\times d})}^p
+\|{\bm E}_\tau^k\|_{L^\rexp(\varOmega;\R^{d\times d\times d})}^\rexp
\,.
\label{transport+++}\end{align}
 for some $C$ sufficiently large, provided $1/p+1/\rexp+1/q'\le1$ with $p$
from \eqref{est-disc+1}. By a discrete Gronwall inequality 
like \eqref{Gronwall}, provided $\tau$ is sufficiently small, we obtain
\begin{align}
  \|\nabla{\bm E}_\tau^k\|_{L^q(\varOmega;\R^{d\times d\times d})}^q\le C
\bigg(\|\nabla{\bm E}_\tau^0\|_{L^q(\varOmega;\R^{d\times d\times d})}^q\!+
  \tau\sum_{l=1}^k\|\nabla{\bm F}_\tau^l\|_{L^{q}(\varOmega;\R^{d\times d\times d})}^{q}\bigg)\,
\label{Gronwall+}\end{align}
with some $C$ depending on
$\max_{l=1,...,k}\|\nabla\vv_\tau^k\|_{L^\infty(\varOmega;\R^{d\times d})}$ and on 
$\tau\sum_{l=1}^k\|\nabla^2\vv_\tau^k\|_{L^p(\varOmega)}^p\!+\|{\bm E}_\tau^k\|_{L^\rexp(\varOmega)}^\rexp$.

Considering $\vv_\tau^k=\vvk$ and ${\bm F}_\tau^k=\EE(\vvk)-\RR_\etau^k$
and using the already obtained estimate \eqref{est-R} and
the initial condition $\Ee_0\in H^1(\varOmega;\R_{\rm dev}^{d\times d})$, 
we use the calculus \eqref{transport+}--\eqref{Gronwall}
with  $\rexp=2^*$ and $q=2$ for \eqref{ED-2+disc}. 
Thus we obtain
\begin{align}
\|\nabla\overlineEetau\|_{L^\infty(I;L^2(\varOmega;\R^{d\times d\times d}))}^{}\le C\,.
\end{align}

The arguments \eqref{transport+}--\eqref{Gronwall} and \eqref{transport+++}--\eqref{Gronwall+}
hold even simplified for the convective transport of scalar variables; in particular
\eqref{transport--} and last two terms in \eqref{transport+++} would be omitted. Here,  having
now $\nabla\overline\theta_\etau$ estimated, by using the calculus \eqref{transport-} with $z=\chi$,
from the discrete volume-fraction evolution \eqref{ED-7+disc}
we thus obtain
\begin{align}\nonumber
  &\frac1q\INT{\varOmega}{}\!
  \frac{|\nabla\chi_\etau^k|^q{-}|\nabla\chi_\etau^{k-1}|^q\!\!}\tau\,\d x
  \le \INT{\varOmega}{}
\Big(\frac{\nabla\theta_\etau^k\!\!}{\theta_\text{\sc pt}}-\vvk{\cdot}\nabla\chi_\etau^k\Big)\cdot|\nabla\chi_\etau^k|^{q-2}\nabla\chi_\etau^k\,\d x
\\&\nonumber\qquad=\INT{\varOmega}{}\!
\frac{\nabla\theta_\etau^k\!}{\theta_\text{\sc pt}}\cdot|\nabla\chi_\etau^k|^{q-2}
\nabla\chi_\etau^k+|\nabla\chi_\etau^k|^{q-2}\big(\nabla\chi_\etau^k{\otimes}\nabla\chi_\etau^k\big)
:\EE(\vvk)-
\frac1q|\nabla\chi_\etau^k|^q{\rm div}\,\vvk\,\d x
\\&\qquad\le \theta_\text{\sc pt}^{-q}\|\nabla\theta_\etau^k\|_{L^q(\varOmega;\R^d)}^q\!+
\big(1{+}\|\EE(\vvk)\|_{L^\infty(\varOmega;\R^{d\times d})}^{}\big)\|\nabla\chi_\etau^k\|_{L^q(\varOmega;\R^d)}^q\,,
\label{grad-chi}\end{align}
where we again used $\vvk{\cdot}\nn=0$. The first inequality in
\eqref{grad-chi} follows from the convexity of $|\cdot|^q$ and from that,
written formally,
$\nabla H^{-1}(\chi_\etau^k)\cdot|\nabla\chi_\etau^k|^{q-2}\nabla\chi_\etau^k
=\partial^2\delta_{[0,1]}^{}(\chi_\etau^k)|\nabla\chi_\etau^k|^q\ge0$, where
$\partial^2\delta_{[0,1]}^{}$ denotes the (generalized) Hessian of the convex
indicator function of the interval $[0,1]$. When $q<(d{+}2)/(d{+}1)$, we can
use \eqref{est-disc+1} and the assumption \eqref{ass:5} on $\chi_0$ and, by
the discrete Gronwall inequality, we obtain a bound for
$\nabla\overlinechitau$ in $L^\infty(I;L^{q}(\varOmega;\R^d))$.

Since now we have $\nabla\overlinechitau$ estimated, we can see that  
$\overlinevvtau{\cdot}\nabla\overlinechitau$ is surely bounded in
$L^1(I{\times}\varOmega)$. Since
$\frac{\partial}{\partial t}\chi_\etau{+}\overlinevvtau{\cdot}\nabla\overlinechitau$
has already been proved bounded in $L^2(I{\times}\varOmega)$,
we thus obtain a bound for $\frac{\partial}{\partial t}\chi_\etau=
(\frac{\partial}{\partial t}\chi_\etau{+}\overlinevvtau{\cdot}\nabla\overlinechitau)
-\overlinevvtau{\cdot}\nabla\overlinechitau$ in $L^1(I{\times}\varOmega)$, and
$\frac{\partial}{\partial t}\overlinechitau$ in ${\rm BV}(I;L^1(\varOmega))$.
To summarize, we proved
\begin{align}\label{est-disc-chi}
  &\|\overlinechitau\|_{L^\infty(I{\times}\varOmega)}^{}\le C\ \ \ \text{ and }\ \ \ 
  \|\overlinechitau\|_{L^\infty(I;W^{1,q}(\varOmega))\,\cap\, {\rm BV}(I;L^1(\varOmega))}^{}\le C_q
\ \ \ \text{ with }\ q<\frac{d{+}2}{d{+}1}\,.
\end{align}

\medskip{\it Step 3: convergence in the mechanical part for $\tau\to0$.}
We now consider $\eps>0$ fixed.
By the Banach selection principle, we obtain a subsequence
converging weakly* with respect to the topologies indicated in
\eqref{est-disc} and \eqref{est-disc+} to some limit
$(\vv_\eps,\Ee_\eps,\RR_\eps,\alpha_\eps,\theta_\eps,\chi_\eps)$.

By the uniform monotonicity of the operators $\frac{\partial}{\partial t}$
and ${\rm div}({\rm div}(\nu|\nabla\EE(\cdot)|^{p-2}\nabla\EE(\cdot))-
\bbD\EE(\cdot))$, we obtain the strong convergence
\begin{subequations}\label{velocity-strong}\begin{align}
    &&&\overlinevvtau\to\vv_\eps&&\text{in }\ L^p(I;W^{2,p}(\varOmega;\R^d))\ \ \text{ and }&&&&
    \\&&&\vv_\etau(T)\to\vv_\eps(T)&&\text{in }\ L^2(\varOmega;\R^d)\,.
\end{align}\end{subequations}\def\whvv{\widehat{\vv}_\tau}
More in detail, we use the discrete momentum equation tested by
$\overlinevvtau{-}\whvv$ with some
$\whvv$ piecewise constant on the time-partition of the
times step $\tau$ and converging strongly to $\vv_\eps$ in
$L^p(I;W^{2,p}(\varOmega;\R^d))\cap
L^\infty(I;L^2(\varOmega;\R^d))$.
After integration over time, we obtain
\begin{align}\nonumber
  &c\|\vv_\etau(T){-}\whvv(T)\|_{L^2(\varOmega;\R^d)}^2+
c\|\nabla\EE(\overlinevvtau{-}\whvv)\|_{L^p(I{\times}\varOmega;\R^{d\times d\times d})}^p
 \le \INT{\varOmega}{}\frac\varrho2|\vv_\etau(T){-}\whvv(T)|^2\d x
  \\[-.3em]&\nonumber
  +\INT{0}{T}\!\!\INT{\varOmega}{}\!\bbD\EE(\overlinevvtau{-}\whvv){:}
  \EE(\overlinevvtau{-}\whvv)
  +\nu\big(|\nabla\EE(\overlinevvtau)|^{p-2}\nabla\EE(\overlinevvtau)
  -|\nabla\EE(\whvv)|^{p-2}\nabla\EE(\whvv)\big)\Vdots
  \nabla\EE(\overlinevvtau{-}\whvv)\,\d x\d t
  \\&\le\nonumber
  \INT{0}{T}\!\!\INT{\varOmega}{}\bigg(
  \overline\ff_\tau{\cdot}(\overlinevvtau{-}\whvv)
  -\big(\overlineStau\!+\overlineSstrtau\big):\EE(\overlinevvtau{-}\whvv)
 -\nu|\nabla\EE(\whvv)|^{p-2}\nabla\EE(\whvv)\Vdots
 \nabla\EE(\overlinevvtau{-}\whvv)
 \\[-.6em]&\hspace{14em}
-\bbD\EE(\whvv){:}\EE(\overlinevvtau{-}\whvv)\bigg)\,\d x\d t
-\INT{\varOmega}{}\varrho\whvv(T){\cdot}(\vv_\etau(T){-}\whvv(T))\,\d x
\label{strong-hyper}\end{align}
with some small $c>0$. As the right-hand side of \eqref{strong-hyper}
converges to 0, we obtain \eqref{velocity-strong}.

Moreover, we use the Aubin-Lions compact-embedding theorem to show
the strong convergence
\begin{subequations}\label{est-strong}\begin{align}\label{est-strong-E}
    &&&\overlineEetau\to\Ee_\eps\ \
\text{ and }\ \ \overlineStau\to\bm{S}_\eps=[\varphi_\eps]_\Ee'(\Ee_\eps,\alpha_\eps)
&&\text{ in }\ L^2(I{\times}\varOmega;\R_{\rm sym}^{d\times d})\,,
&&&&
\\&&&\label{est-strong-a}
    \overlinealphatau\to\alpha_\eps&&\text{ in }\
    L^q(I{\times}\varOmega),\ \ q<\infty\,;
\end{align}\end{subequations}
in fact, \eqref{est-strong-E} holds even in better spaces (i.e.\ in stronger
modes). For this, we need some information about time derivatives. It is
directly at disposal for \eqref{est-strong-a} from the a-priori estimate
\eqref{est-disc-chi}.
For the former convergence in \eqref{est-strong-E}, we can see some
information by comparison from \eqref{ED-2+d},
from which we can see that the sequence
$\{\frac{\partial}{\partial t}\Ee_\etau\}_{\tau>0}$ is bounded in
$L^2(I;H^1(\varOmega;\R_{\rm sym}^{d\times d})^*)$. Actually,
usage of the mentioned Aubin-Lions theorem for piece-wise constant
functions in time needs its generalized version for functions whose
time-derivatives are measures, cf.\ \cite[Cor.7.9]{Roub13NPDE}.
The later strong convergence in \eqref{est-strong-E} is then simply by
continuity of the Nemytski\u{\i} operator induced by $[\varphi_\eps]_\Ee'$,
even without having any information about the time derivative.

As $\nabla\alpha$ occurs nonlinearly in the stress $\bm{K}$ and
also in the weak formulation \eqref{VI} multiplied by $\Delta\alpha$,
we need to prove also the strong convergence
  \begin{align}
&&&\nabla\overlinealphatau\to\nabla\alpha_\eps&&\hspace*{-2em}\text{strongly in }\ L^2(I{\times}\varOmega;\R^{d})\,.&&&&
\label{nabla-dam-strong}\end{align}
To prove it, we take a sequence $\{\widetilde\alpha_\tau\}_{\tau>0}^{}$
piecewise constant in time with respect to the partition with the time
step $\tau$  and converging strongly towards 
$\alpha_\eps$ for $\tau\to0$. Then, using the variational inequality
arising from the inclusion \eqref{ED-4+d} tested by
$\nabla(\overlinealphatau{-}\widetilde\alpha_\tau)$ 
we can see that, written in terms of the interpolants,
\begin{align}\nonumber
\nonumber
  &\INT{0}{T}\!\!\INT{\varOmega}{} \kappa|\nabla(\overlinealphatau{-}\widetilde\alpha_\tau)|^2\,\d x\d t
  \le-\INT{0}{T}\!\!\INT{\varOmega}{}
\bigg(\partial_{\DT\alpha}\zeta\Big(\underline\alpha_\etau,\underline\W_\etau;
\pdt{\alpha_\etau}+
\overlinevvtau{\cdot}\nabla\overlinealphatau
\Big)+\sqrt\tau\pdt{\alpha_\etau}
\\[-.3em]&\hspace{15em}+[\varphi_\eps]_\alpha'(\overlineEetau,\overlinealphatau)
\bigg)(\overlinealphatau{-}\widetilde\alpha_\tau)
+\kappa\nabla\widetilde\alpha_\etau{\cdot}\nabla(\overlinealphatau{-}\widetilde\alpha_\tau)\,\d x\d t
  \to0\,.
\label{strong-nabla-alpha}
\end{align}
Here we used \eqref{ass:1} with \eqref{est-strong-E} so that
$[\varphi_\eps]_\alpha'(\overlineEetau,\overlinealphatau)$ is bounded in
$L^{2}(I{\times}\varOmega)$ 
while $\overlinealphatau-\widetilde\alpha_\tau\to0$ strongly in
$L^{q}(I{\times}\varOmega)$ with any $q<+\infty$;
this is due to an interpolated Aubin-Lions theorem, relying on
\eqref{est-disc4} with \eqref{est-disc+4}.
In \eqref{strong-nabla-alpha}, we used that
$\partial_{\DT\alpha}\zeta(\underline\alpha_\etau,\underline\W_\etau;
\pdt{}\alpha_\etau{+}\overlinevvtau{\cdot}\nabla\overlinealphatau)$ is bounded
in $L^2(I{\times}\varOmega)$ and, moreover, that
$\|\sqrt\tau\pdt{}\alpha_\etau\|_{L^2(I\times\varOmega)}^{}
=\mathscr{O}(\sqrt\tau)\to0$ due to \eqref{est-disc+4}.

As $\GM(\cdot)$ and $\zeta(\alpha,\cdot;\DT\alpha)$ depend nonlinearly on
the enthalpy $\W=\vartheta+\ell\chi$, we need a strong convergence of both
$\underline\vartheta_\etau$ and $\underline\chi_\etau$. The former one follows
from the estimates on the gradient \eqref{est-disc+vartheta}
and from some information about time derivative of $\underline\vartheta_\etau$ 
from the discrete equations \eqref{ED-6+d} by the (generalized) Aubin-Lions
compact embedding theorem. More in detail, from
\eqref{ED-6+d} we can see that $\frac{\partial}{\partial t}\underline\W_\etau
\in{\rm BV}(I;H^2(\varOmega)^*)$ so that, taking
$\frac{\partial}{\partial t}\underline\chi_\etau\in{\rm BV}(I;L^1(\varOmega))$
proved already before in \eqref{est-disc-chi}, we obtain
$\frac{\partial}{\partial t}\underline\vartheta_\etau=
\frac{\partial}{\partial t}\underline\W_\etau
-\ell\frac{\partial}{\partial t}\underline\chi_\etau
\in{\rm BV}(I;H^2(\varOmega)^*)$. Thus we proved the strong convergence
\begin{align}
\underline\W_\etau\!=
\widetilde\gamma(\underline\theta_\etau)+\ell\underline\chi_\etau\!\to
\widetilde\gamma(\theta_\eps)+\ell\chi_\eps\!=\W_\eps
\ \text{ strongly in }\ L^s(I{\times}\varOmega)
\label{beta-discrete}\end{align}
with $s$ from \eqref{est-disc+vartheta}.
For the latter one, i.e.\ for $\underline\chi_\etau\to\chi_\eps$ strongly in
$L^q(I{\times}\varOmega)$ for any $q<+\infty$, we use an interpolated (and
generalized) Aubin-Lions compact embedding theorem exploiting the a-priori
estimates \eqref{est-disc-chi}.

The limit passage in the variational inequality for $\overlinealphatau$
which is behind the inclusion \eqref{ED-4+d}, cf.\ \eqref{VI},
exploits $[\varphi_\eps]_\alpha'(\overlineEetau,\overlinealphatau)\to
[\varphi_\eps]_\alpha'(\Ee_\eps,\alpha_\eps)$ strongly in
$L^2(I{\times}\varOmega)$.

The limit passage the variational inequality,
which is behind the volume-fraction equation \eqref{ED-7+d}, towards the
variational inequality \eqref{ED-7+weak} written for the $\eps$-solution is
based on the a-priori estimates \eqref{est-disc-chi}. From them, as we already
used for $\underline\chi_\etau$, using the Aubin-Lions theorem (again
generalized for functions with measure time derivatives as in
\cite[Cor.7.9]{Roub13NPDE}) we obtain that
$\overlinechitau\to\chi_\varepsilon$ in $L^q(I{\times}\varOmega)$ to be used for
the term $\frac12({\rm div}\,\overlinevvtau)\,
\overline{\hspace*{-.1em}\chi}_{\etau}^2$ with any $q<+\infty$. Also we use the
weak convergence of $\overline{\hspace*{-.1em}\chi}_{\etau}(T)\to\chi_\eps(T)$
in $L^q(\varOmega)$ with any $q<+\infty$. 

Using also (\ref{est-strong}a,b), we can see that the Korteweg-like stress
converges even strongly, namely
\begin{align}
\overlineSstrtau\to{\bm K}_\eps^{}\ \ \ \ 
  \text{ in }\ L^s(I;L^1(\varOmega;\R^{d\times d})),\ \ 1<s<\infty\,.
\end{align}

\medskip{\it Step 4: Convergence of the dissipation rate
and of the heat-transfer equation for $\tau\to0$.}
Here we use the already proved convergence in the mechanical part
together with the mechanical-energy conservation.
We have the chain of estimates:
\begin{align}\nonumber
&
 \INT{0}{T}\!\!\INT{\varOmega}{}\GM(\W_\varepsilon)|\RR_\varepsilon|^2+
\bbD\EE(\vv_\varepsilon){:}\EE(\vv_\varepsilon)
+\nu|\nabla\EE(\vv_\varepsilon)|^p+\xi\big(\alpha_\varepsilon,\W_\varepsilon;
\DT\alpha_\varepsilon\big)
+\varkappa|\nabla\RR_\varepsilon|^2\,\d x\d t
\\[-.4em]\nonumber
&\quad\le\liminf_{\tau\to0}\INT{0}{T}\!\!\INT{\varOmega}{}\!\!
\GM(\underline\W_\etau)|\overlineRRtau|^2+
\bbD\EE(\overlinevvtau){:}\EE(\overlinevvtau)
+\nu|\nabla\EE(\overlinevvtau)|^p
\\[-.6em]\nonumber&\hspace{21.5em}
+\xi\Big(\underline\alpha_\etau,\underline\W_\etau;
 \frac{\partial\alpha_\etau}{\partial t}{+}\overlinevvtau{\cdot}\nabla\overlinealphatau\Big)
 +\varkappa|\nabla\overlineRRtau|^2\,\d x\d t
 \\[-.3em]
 \nonumber&\ \,\stackrel{\eqref{mech-energy-disc}}{\le}
 \!\INT{\varOmega}{}\frac\varrho2|\vv_0|^2
+\varphi_\varepsilon(\Ee_0,\alpha_0)+\frac{\kappa}2|\nabla\alpha_0|^2\,\d x
- \liminf_{\tau\to0}\INT{\varOmega}{}\frac\varrho2|\vv_\etau(T)|^2
+\varphi_\varepsilon(\Ee_\etau(T),\alpha_\etau(T))+\frac{\kappa}2|\nabla\alpha_\etau(T)|^2\,\d x
 \\[-.4em]\nonumber&\hspace{13.7em}+\lim_{\tau\to0}
 \INT{0}{T}\!\!\INT{\varOmega}{}\!\overline\ff_\tau\!\cdot\overlinevvtau
 -\phi(\overline\theta_\etau){\rm div}\,\overlinevvtau
 +\frac{\sqrt\tau}2\Big|\frac{\partial\alpha_\etau}{\partial t}
 {+}\overlinevvtau{\cdot}\nabla\overlinealphatau\Big|\,\d x \d t
 \\ \nonumber&\quad\le\INT{\varOmega}{}\frac\varrho2|\vv_0|^2
+\varphi_\varepsilon(\Ee_0,\alpha_0)+\frac{\kappa}2|\nabla\alpha_0|^2
 -\frac\varrho2|\vv_\varepsilon(T)|^2
-\varphi_\varepsilon(\Ee_\varepsilon(T),\alpha_\varepsilon(T))-\frac{\kappa}2|\nabla\alpha_\varepsilon(T)|^2\,\d x
\\[-.4em]\nonumber&\hspace{26.7em}
+\INT{0}{T}\!\!\INT{\varOmega}{}\!\ff{\cdot}\vv_\varepsilon
-\phi(\theta_\varepsilon){\rm div}\vv_\varepsilon\,\d x\d t
 \\[-.2em]
 &\quad=
 \INT{0}{T}\!\!\INT{\varOmega}{}\!\!\GM(\W_\varepsilon)|\RR_\varepsilon|^2+
\bbD\EE(\vv_\varepsilon){:}\EE(\vv_\varepsilon)
+\nu|\nabla\EE(\vv_\varepsilon)|^p+\xi\big(\alpha_\varepsilon,\W_\varepsilon;
\DT\alpha_\varepsilon\big)
 +\varkappa|\nabla\RR_\varepsilon|^2\,\d x\d t\,.
 \label{limsup-trick}\end{align}
The first inequality is due to the weak lower semicontinuity.
The last equality in \eqref{limsup-trick} is just the 
mechanical-energy balance \eqref{energy+} written for the
$\varepsilon$-solution.

This mechanical-energy balance follows from the tests as in the (formal)
calculations \eqref{test-velocity}--\eqref{calculus-toward-pressure}
written for the $\varepsilon$-solution. The validity of this balance 
is not automatic and the rigorous prove needs to have granted that
testing the particular equations by $\vv_\eps$, ${\bm S}_\eps$, and $\RR_\eps$ is indeed
legitimate. Here, in particular it is important that
$\bm{K}_\eps{:}\nabla\vv_\eps\in L^1(I{\times}\varOmega)$
because $\bm{K}_\eps\in L^\infty(I;L^1(\varOmega;\R^{d\times d}))$
and $\nabla\vv_\eps\in L_{\rm w*}^p(I;L^\infty(\varOmega;\R^{d\times d}))$.
Also, $\kappa\Delta\alpha_\eps\in
\partial_{\DT\alpha}\zeta(\alpha_\eps,\W_\eps;\DT{\alpha}_\eps)
+[\varphi_\eps]_\alpha'(\Ee_\eps,\alpha_\eps)$ holds pointwise a.e.\ in the sense
of $L^2(I{\times}\varOmega)$, 
we can legitimately test it by $\DT\alpha_\eps\in L^2(I{\times}\varOmega)$. Since 
$\partial_{\DT\alpha}\zeta(\alpha_\eps,\W_\eps;\cdot)$ is single-valued
except 0, cf.\ \eqref{ass:2},
$\partial_{\DT\alpha}\zeta(\alpha_\eps,\W_\eps;\DT\alpha_\eps)\DT\alpha_\eps$ is
single-valued, being equal to
$\xi(\alpha_\eps,\W_\eps;\DT\alpha_\eps)\in L^1(I{\times}\varOmega)$.

Altogether, this reveals that
there are actually equalities in \eqref{limsup-trick}.
Since the dissipation rate is uniformly convex in terms of rates on the
uniformly convex $L^2$-spaces, these rates converge not only weakly but
even strongly in these  $L^2$-spaces. Thus the dissipation rate itself
converges strongly in $L^1(I{\times}\varOmega)$.

The limit passage in the resting semilinear terms in \eqref{ED-6+d}
is then easy.

\medskip{\it Step 5: convergence for $\eps\to0$.} This final convergence
towards the weak solution due to Definition~\ref{def} copies the
arguments in the Steps~2-4 above. Actually,
derivation of the estimates \eqref{est-disc} and (\ref{est-disc+}a,b) in
Step~2 is even simplified because the time-continuous problem can rely
directly on the calculus from  Section~\ref{sec-fuller}. The manipulation
exploiting the semi-convexity \eqref{semi-convex} of $\phi_\eps$ is no
longer needed.

By \eqref{phi-eps0}, we have $\varphi_\eps(\Ee_\eps,\alpha_\eps)
=\varphi(\Ee_\eps,\alpha_\eps)+o(\eps)(1{+}
\|\Ee_\eps\|_{L^2(I{\times}\varOmega;\R^{d\times d})})\to\varphi(\Ee,\alpha)$
in $L^1(I{\times}\varOmega)$ to be used for the convergence in the
Korteweg-like stress, cf.\ \eqref{calc-energy-pressure} written for
$\eps$-solution. Also we have
$|\varphi_\eps(\Ee_\eps(T),\alpha_\eps(T))-\varphi(\Ee(T),\alpha(T))|
\le o(\eps)(1{+}|\Ee(T)|^2)$ so that
\begin{align}\nonumber
  \liminf_{\eps\to0}\!\INT{\varOmega}{}\!\varphi_\eps(\Ee_\eps(T),\alpha_\eps(T))\,\d x
\ge\liminf_{\eps\to0}&\!\INT{\varOmega}{}\!\varphi(\Ee_\eps(T),\alpha_\eps(T))\,\d x
\\[-.7em]&\ \ \ \ -o(\eps)\big(1{+}\|\Ee(T)\|_{L^2(\varOmega;\R^{d\times d})}^2\big)
\ge\INT{\varOmega}{}\!\varphi(\Ee(T),\alpha(T))\,\d x
\end{align}
to be used for \eqref{limsup-trick} written between the $\eps$-solutions
and their limit.

Furthermore, from \eqref{phi-eps1} we have
$[\varphi_\eps]_\Ee'(\Ee_\eps,\alpha_\eps)\to\varphi_\Ee'(\Ee,\alpha)$ in
$L^{1/\delta}(I;L^{2^*-\delta}(\varOmega;\R^{d\times d}))$ for any $\delta>0$
to be used e.g.\ in \eqref{est-strong-E} written
for $\eps$-solutions. Eventually, from \eqref{phi-eps2} we have
$[\varphi_\eps]_\alpha'(\Ee_\eps,\alpha_\eps)\to\varphi_\alpha'(\Ee,\alpha)$ in
$L^{1/\delta}(I;L^{2^*/2-\delta}(\varOmega))$, which is to be used in
\eqref{strong-nabla-alpha} written between $\eps$-solution and the
limit. 

The weak-continuity \eqref{temp-w-cont} follows from the 
estimate \eqref{est-disc-w} inherited for $\theta$ first as the estimate in
$L_{\rm w*}^\infty(I;{\rm Meas}(\barOmega))$. Further, due to the estimate for
$\frac{\partial}{\partial t}\theta$,
we have also $C_{\rm w*}(I;{\rm Meas}(\barOmega))$. Yet, since
$\theta(t)\in W^{1,r}(\varOmega)\subset L^1(\varOmega)$ for a.a.\ $t\in I$, the
weak* continuity of $\theta:I\to {\rm Meas}(\barOmega)$ turns into
the weak continuity of $\theta:I\to L^1(\varOmega)$.
%\end{proof}
$\hfill\Box$

\begin{remark}[{\sl Inelastic strain}]\upshape
Given an additional initial condition $\bm{P}_0^{}$, the inelastic strain
$\bm{P}$ can be reconstructed from $\ZJ{\bm{P}}=\ZJEp$, cf.\ 
Remark~\ref{rem-plast-strain}. Since 
$\ZJEp\in L^2(I;H^1(\varOmega;\R^{d\times d}))$, assuming $\bm{P}_0^{}\in L^q(\varOmega;\R^{d\times d})$,
using the estimation as in \eqref{transport+}--\eqref{Gronwall}, 
one can see that $\bm{P}\in L^\infty(I;L^q(\varOmega;\R^{d\times d}))$. If even
$\bm{P}_0^{}\in H^1(\varOmega;\R^{d\times d})$,
the estimation as in \eqref{transport+++}--\eqref{Gronwall+}
gives $\bm{P}\in L^\infty(I;H^1(\varOmega;\R^{d\times d}))$.
\end{remark}

\begin{remark}[{\sl Regularity}]\upshape
  On smooth or convex domains $\varOmega$, one can use the $H^2$-regularity
  for damage to be seen by comparison from $\Delta\alpha\in
  (\partial_{\DT\alpha}\zeta(\alpha,\W;\DT{\alpha})
  +[\varphi_\eps]_\alpha'(\Ee,\alpha))/\kappa\subset L^2(I{\times}\varOmega)$
  bounded. Knowing $\alpha\in L^2(I;H^2(\varOmega))$, we have also
  a better integrability of the Korteweg-like stress, namely 
  $\bm{K}\in L^\infty(I;L^1(\varOmega;\R^{d\times d}))\,\cap\,
  L^{1+2/d}(I{\times}\varOmega;\R^{d\times d})$. To show \eqref{nabla-dam-strong},
  one can then use the Aubin-Lions theorem instead of
  \eqref{strong-nabla-alpha}.
\end{remark}

\begin{remark}[{\sl Regularization of the stored energy}]\label{rem-example}\upshape
In the case of $\varphi$ from \eqref{psi}, the regularization 
satisfying the condition \eqref{phi-eps} can be done as
\begin{align}
\varphi_\eps(\bm E,\alpha)
=\frac d2K_\text{\sc e}^{}|{\rm sph}\,\Ee|^2+G_\text{\sc e}^{}(\alpha)
\,\frac{|{\rm dev}\,\Ee|^2}{\!\sqrt{1{+}\eps|{\rm dev}\,\Ee|^2}}
  +\frac{1}{2\kappa}G_\text{\sc d}^{}(1{-}\alpha)^2  \,.
\label{psi+}\end{align}
This ensures semiconvexity \eqref{semi-convex} of $\varphi_\eps$.
Simultaneously, the conditions (\ref{phi-eps}a,c-e) hold.
\end{remark}

\section{Concluding remarks}\label{sec-rem}
%        ~~~~~~~~~~~~~~~~~~

We close this article by several remarks to the possible modifications,
generalization, and expansion of the model and its applications.
Anyhow, several important aspects still remain open - in particular
a pressure-dependent phase-transformation temperature.

\begin{remark}[{\sl Omitting damage in solid part  --- stress formulation}]\label{rem-no-damage}\upshape
When damage (or the phase-field fracture) would 
not be considered so that the elastic moduli are fixed, one can come back
to the stress/velocity formulation like in Sect.~\ref{sec-preliminary}.
Confining ourselves to an isotropic material like in \eqref{C},
the ``undamageable'' stored energy in terms of stresses can be considered again
as \eqref{phi+} but now together with the additive decomposition of the total
 stress  $\bm{S}+{\bm B}$ to the actual elastic stress $\bm{S}$
and the  deviatoric  ``back stress'' ${\bm B}$
 and the stored energy $\varphi^*$ from \eqref{phi+}. 
In terms of a back stress rate ${\bm R}=\ZJ{\bm B}$, the resulting system then
reads as
\begin{subequations}\label{ED+thermo}\begin{align}
&\varrho\DT\vv={\rm div}\big({\bm S}+\varphi^*(\bm{S})\bbI+{\bm D}
              \big)-\frac\varrho2({\rm div}\,\vv)\,\vv+\ff\, \ \text{ with }\
               {\bm D}=\bbD\EE(\vv)-{\rm div}\big( 
\nu|\nabla\EE(\vv)|^{p-2}\nabla\EE(\vv)\big)\,,
\label{ED-1+thermo}
\\[-.5em]&
\ZJ{\bm S}=\big[[\varphi^*]'\big]^{-1}\EE(\vv)-{\bm R}\,,
 \label{ED-2+thermo}
\\\label{ED-3+thermo}
&\GM(\W){\bm R}=2G_\text{\sc e}^{}{\rm dev}{\bm S}+\varkappa\Delta{\bm R}\,,
\\\nonumber
&\pdt\W+{\rm div}\big(\vv\,\W
-\bbK(\alpha,\W)\nabla\theta\big)
=\frac{\GM(\W)}{4G_\text{\sc e}^{2}}|{\bm R}|^2
+\frac{\varkappa}{4G_\text{\sc e}^{2}}|\nabla{\bm R}|^2+\bbD\EE(\vv){:}\EE(\vv)
\\[-.1em]&\label{ED-6+thermo+}
\hspace*{13.5em}
+\nu|\nabla\EE(\vv)|^p
+\phi(\theta){\rm div}\,\vv\ \ \ \ \text{ with }\ \
\W=\widetilde\gamma(\theta)+\ell\chi\,,
\\[-.5em]&
\NU\DT\chi+H^{-1}(\chi)\ni\frac{\theta_\etau^k}{\theta_\text{\sc pt}}-1
\,;
\end{align}\end{subequations}
in \eqref{ED-2+thermo}, $[[\varphi^*]']^{-1}=\varphi'$ with the convex
conjugate $\varphi$ to $\varphi^*$. 
The energetics in the mechanical part is now revealed by testing
\eqref{ED-1+thermo} by $\vv$, \eqref{ED-2+thermo} by  the strain  
$[\varphi^*]'(\bm{S})$, and \eqref{ED-3+thermo} by
${\bm R}$. When added still \eqref{ED-6+thermo+} tested by 1, one
obtain the total energy conservation. Altogether, we obtain
\eqref{energy+}--\eqref{energy+++} but with
$\varphi^*(\bm{S})$ instead of $\varphi(\Ee,\alpha)$, with
$\zeta$ in \eqref{energy+} omitted, and with
$\GM(\W)|\ZJEp|^2+\varkappa|\nabla\ZJEp|^2$ replaced by
$\GM(\W)|{\bm R}|^2/(2G_\text{\sc e})+
\varkappa|\nabla{\bm R}|^2/(2G_\text{\sc e})$.
 For an isothermal incompressible variant with nonlinear creep (involving
possibly also plasticity) and with ${\bm B}$ eliminated
see \cite{EiHoMi21LHSV}.
\end{remark}

\begin{remark}[{\sl Nonlinear or even activated creep}]\label{rem-nonlin-creep}\upshape
  The (materially) linear creep law \eqref{ED-3+} is an idealization and
  in some materials such linear Maxwellian rheology may not be realistic.
  Then a non-quadratic (or even nonsmooth at 0) dissipation potential $\zeta_{\rm crp}^{}:\R_{\rm dev}^{d\times d}\to R$ should be considered and 
  \eqref{ED-3+} is to be generalized as
\begin{align}
  \GM(\W)\partial\zeta_{\rm crp}^{}\big(\ZJEp)
    \ni{\rm dev}\,\bm{S}+\varkappa\Delta\ZJEp\,,
\label{inclusion}\end{align}
where $\partial\zeta_{\rm crp}^{}$ denotes the convex subdifferential of
$\zeta_{\rm crp}^{}$. By the definition of the convex subdifferential,
\eqref{inclusion} integrated over $I{\times}\varOmega$ means the variational
inequality
\begin{align}
  \INT{0}{T}\!\!\INT{\varOmega}{}\zeta_{\rm crp}^{}\big(\widetildeEp\big)
  -\bm{S}{:}\big(\widetildeEp{-}\ZJEp\big)+\varkappa\nabla\ZJEp\Vdots
  \nabla\widetildeEp\,\d x\d t
\ge\INT{0}{T}\!\!\INT{\varOmega}{}\zeta_{\rm crp}^{}\big(\ZJEp\big)+
\varkappa|\nabla\ZJEp|^2\,\d x\d t
\label{VI+}\end{align}
to hold for any $\widetildeEp\in L^q(I{\times}\varOmega;\R_{\rm dev}^{d\times d})\cap
L^2(I;H^1(\varOmega;\R_{\rm dev}^{d\times d}))$ if we assume
coercivity $\zeta_{\rm crp}^{}(\cdot)\ge\epsilon|\cdot|^q$ with some $q>1$.
Then $\RR\in L^{q}(I{\times}\varOmega;\R_{\rm dev}^{d\times d})\cap
L^2(I;H^1(\varOmega;\R_{\rm dev}^{d\times d}))$. 
To converge with the approximate
solution in \eqref{VI+}, we use the strong convergence \eqref{est-strong-E} of
$\bm{S}=[\varphi_\eps]_\Ee'(\Ee,\alpha)$ and weak convergence of
$\ZJEp$ together with convexity of $\zeta_{\rm crp}^{}+\varkappa|\cdot|^2$.
\end{remark}

\begin{remark}[{\sl A phenomenological model of ice creep}]\upshape
  An example of a nonlinear creep from Remark~\ref{rem-nonlin-creep}
  is a phenomenological Glen's law
  \cite{Glen55CPI} for polycrystalline ice: strain rate is proportional to
  the $n$-power of stress (with $n$ around 3).
This fits with the potential $\zeta_{\rm crp}^{}(\cdot)=G_0|\cdot|^q$ with
$q=1+1/n$. When $\zeta_{\rm crp}^{}$ is nonsmooth at 0, we obtain an activated
creep (plasticity), which is also a phenomenon sometimes considered for ice.
\end{remark}

 \begin{remark}[{\sl Rate-and-state friction model}]\upshape
  Letting the dissipation potential $\zeta_{\rm crp}^{}$ state dependent, i.e.\
  dependent on $\alpha$ and $\W$, allows for mimicking the popular
  Dieterich-Ruina's \cite{Diet79MRF,Ruin83SISV} rate-and-state friction model
  for sliding of rocks along tectonic faults, cf.\ \cite[Sect.\,6]{Roub14NRSD}.
\end{remark}

\begin{remark}[{\sl Adiabatic effects: buoyancy}]\upshape
  Another example of the adiabatic effects
  is the Oberbeck-Boussinesq' simplified buoyancy model used
  for incompressible media in an external (typically gravitational) field.
  Although these
  media are modelled incompressible, they anyhow exhibit a slight thermal
  expansibility, which gives rise to the extra force by replacing
  $\ff$ with $\ff(1{-}b(\theta))$ in \eqref{ED-1+} 
  with $b(\cdot)$ continuous. Then \eqref{ED-6+thermo} expands 
  by the adiabatic heat source/sink $b(\theta)\vv{\cdot}\ff$.
  Like for the term $\phi(\theta){\rm div}\,\vv$, the a-priori estimation
  strategy uses several-times the Gagliardo-Nirenberg interpolation, 
 and essentially relies on that, if $\ff$ assumed bounded
  and $b(\cdot)$ has a certain at most polynomial (even superlinear) growth,
  $\theta\in L^{1+2/d-\epsilon}(I{\times}\varOmega;\R^d)$ is well in duality with
  $\vv\in L^\infty(I;L^2(\varOmega;\R^d))\,\cap L^p(I;W^{2,p}(\varOmega;\R^d))$,
cf.\ the estimates \eqref{est-disc1} and \eqref{est-disc+1},
which actually can be proved for this
adiabatic coupling, too. Here, during solid-fluid phase transformation,
the density (approximated as constant) may in fact jump up, which leads to
discontinuity of $b$ at $\theta_\text{\sc pt}$. E.g.\ water
is about 8\%  heavier than ice at $\theta_\text{\sc pt}=273^\circ$C. Such a
discontinuity may easily be treated by making the buoyancy $b$ dependent
continuously on $\W$ rather than on $\theta$,
similarly as we did for $\GM$, cf.\ Figure~\ref{fig3}.
Such jump is e.g.\ during ice/water transition and
  makes ice floats (icebergs) floating rather on the sea surface
  because ice is lighter than water even if both are on the same
  melting/freezing temperature, which comes back to the original motivation
  of Josef Stefan in melting of arctic floats (that time without any
  mechanical context, of course).
\end{remark}

\begin{remark}[{\sl Traction load}]\label{rem-g}\upshape
Considering nonhomogeneous mechanical boundary condition with the load
${\bm g}$ as \eqref{ED-BC} also in \eqref{BC} would bring technical
difficulties in estimation strategy. Note that the term
$\int_{\varGamma}{}{\bm g}_\tau^k{\cdot}\vv_\etau^k\d S$ arising in
\eqref{energy-disc} would not be well controlled. One would need to add heat
equation tested by 1/2 (instead of 1) to \eqref{mech-energy-disc} to see some
dissipation term to be able to control the traces of $\vv_\etau^k$ on
$\varGamma$ but, on the other hand, an estimation of the resting adiabatic term
  $\phi(\theta_\etau^k){\rm div}\vv_\etau^k/2$ is to be performed, for which
  we would need to impose some polynomial growth
  of heat capacity for $\theta\to\infty$.
\end{remark}

\section*{Acknowledgments}
The author is thankful for many interesting valuable discussions with
Zden\v ek Fiala, Sebastian Schwarzacher, and Giuseppe Tomassetti.


\baselineskip=10pt

\begin{thebibliography}{10}

\vspace*{-.4em}\bibitem{AmbTor90AFDJ}
L.~Ambrosio and V.M. Tortorelli.
\newblock Approximation of functional depending on jumps via by elliptic
  functionals via {$\Gamma$}-convergence.
\newblock {\em Comm. Pure Appl. Math.}, 43:999--1036, 1990.

\vspace*{-.4em}\bibitem{BabSob08HRNM}
A.Y. Babeyko and S.V. Sobolev.
\newblock High-resolution numerical modeling of stress distribution
  in visco-elasto-plastic subducting slabs.
\newblock {\em Lithos}, 103:205--216, 2008.

\vspace*{-.4em}\bibitem{BeBlNe92PBMV}
H.~Bellout, F.~Bloom, and J.~Ne\v{c}as.
\newblock Phenomenological behavior of multipolar viscous fluids.
\newblock {\em Qarterly Appl. Math.}, 1:559--583, 1992.

\vspace*{-.4em}\bibitem{Biot65MID}
M.A. Biot.
\newblock {\em Mechanics of Incremental Deformation}.
\newblock J. Wiley, New York, 1965.

\vspace*{-.4em}\bibitem{BocGal89NEPE}
L.~Boccardo and T.~Gallou{\"e}t.
\newblock Non-linear elliptic and parabolic equations involving measure data.
\newblock {\em J. Funct. Anal.}, 87:149--169, 1989.

\vspace*{-.4em}\bibitem{BoMiRo09CDPS}
G.~Bouchitt\'e, A.~Mielke, and T.~Roub\'{\i}\v{c}ek.
\newblock A complete damage problem at small strains.
\newblock {\em Zeitschrift angew. Math. Phys.}, 60:205--236, 2009.

\vspace*{-.4em}\bibitem{Bruh09EEBI}
O.T. Bruhns and A.~Meyers.
\newblock Eulerian elastoplasticity: basic issues and recent results.
\newblock {\em Theoret. Appl. Mech.}, 36:167--205, 2009.

\vspace*{-.4em}\bibitem{BuMaMi17SDRN}
J.~Burczak, J.~M\'alek, and P.~Minakowski.
\newblock Stress-diffusive regularization of non-dissipative rate-type
  materials.
\newblock {\em Disc. Cont. Dynam. Systems - S}, 10:1233--1256, 2017.

\vspace*{-.4em}\bibitem{Chal64PS}
B.~Chalmers.
\newblock {\em Principles of Solidification}.
\newblock Wiley/ Krieger, New York, 1964/1977.

\vspace*{-.4em}\bibitem{ColGra93HPCP}
P.~Colli and M.~Grasselli.
\newblock Hyperbolic phase change problems in heat conduction with memory.
\newblock {\em Proc. Royal Soc. Edinburgh}, 123A:571--592, 1993.

\vspace*{-.4em}\bibitem{ColRec02CSPP}
P.~Colli and V.~Recupero.
\newblock Convergence to the {S}tefan problem of the phase relaxation problem
  with {C}attaneo heat flux law.
\newblock {\em J. Evol. Equ.}, 2:177--195, 2002.

\vspace*{-.4em}\bibitem{ColSpr95PFMZ}
P.~Colli and J.~Sprekels.
\newblock On a {P}enrose-{F}ife model with zero interfacial energy leading to a
  phase-field system of relaxed {S}tefan type.
\newblock {\em Ann. Mat. Pura Appl.}, 169:269 -- 289, 1995.

\vspace*{-.4em}\bibitem{DaRoSt??NHFV}
E.~Davoli, T.~Roub{\'{\i}}{\v{c}}ek, and U.~Stefanelli.
\newblock A note about hardening-free viscoelastic models in {M}axwellian-type
  rheologies.
\newblock {\em Math. \& Mech. of Solids}, 26:1483--1497, 2021.

\vspace*{-.4em}\bibitem{Diet79MRF}
J.H. Dieterich.
\newblock Modelling of rock friction. {P}art 1: {E}xperimental results and
  constitutive equations.
\newblock {\em J. Geo-phys. Res.}, 84(B5):2161--2168, 1979.

\vspace*{-.4em}\bibitem{EiHoMi21LHSV}
T.~Eiter, K.~Hopf, and A.~Mielke.
\newblock Leray-{H}opf solutions to a viscoelastoplastic fluid model with
  nonsmooth stress-strain relation.
\newblock {\em Nonlin. Anal. Real world Appl. - to appear}, 2021.

\vspace*{-.4em}\bibitem{FeiMal06NSET}
E.~Feireisl and J.~M\'alek.
\newblock On the {N}avier-{S}tokes equations with temperature-dependent
  transport coefficients.
\newblock {\em Diff. Equations Nonlin. Mech.}, pages 14pp.(electronic), Art.ID
  90616, 2006.

\vspace*{-.4em}\bibitem{Fial11GSSM}
Z.~Fiala.
\newblock Geometrical setting of solid mechanics.
\newblock {\em Annals of Physics}, 326:1983--1997, 2011.

\vspace*{-.4em}\bibitem{Fial20OTDR}
Z.~Fiala.
\newblock Objective time derivatives revised.
\newblock {\em Z. angew. Math. Phys.}, 71:Art.no.4, 18 pages, 2020.

\vspace*{-.4em}\bibitem{FriGur06TBBC}
E.~Fried and M.E. Gurtin.
\newblock Tractions, balances, and boundary conditions for nonsimple materials
  with application to liquid flow at small-length scales.
\newblock {\em Archive Ration. Mech. Anal.}, 182:513--554, 2006.

\vspace*{-.4em}\bibitem{FukKen05SPCG}
T.~Fukao and N.~Kenmochi.
\newblock Stefan problems with convection governed by {N}avier-{S}tokes
  equations.
\newblock {\em Advances Math. Sci. Appl.}, 15:29--48, 2005.

\vspace*{-.4em}\bibitem{Gery19INGM}
T.V. Gerya.
\newblock {\em Introduction to Numerial Geodynamic Modelling}.
\newblock 2nd. ed., Cambridge Univ. Press, New York, 2019.

\vspace*{-.4em}\bibitem{Glen55CPI}
J.W. Glen.
\newblock The creep of polycrystalline ice.
\newblock {\em Proc. Royal Soc. London, Ser. A}, 228:519--538, 1955.

\vspace*{-.4em}\bibitem{GreNag65GTEP}
A.~Green and P.~Naghdi.
\newblock A general theory of an elastic-plastic continuum.
\newblock {\em Arch. Ration. Mech. Anal.}, 18:251--281, 1965.

\vspace*{-.4em}\bibitem{Haup02CMTM}
P.~Haupt.
\newblock {\em Continuum Mechanics and Theory of Materials}.
\newblock Springer, Berlin, second edition, 2002.

\vspace*{-.4em}\bibitem{HinZie07OCFB}
M.~Hinze and S.~Ziegenbalg.
\newblock Optimal control of the free boundary in a two-phase {S}tefan problem
  with flow driven by convection.
\newblock {\em Z. Angew. Math. Mech.}, 87:430--448, 2007.

\vspace*{-.4em}\bibitem{Jaum11GSPC}
G.~Jaumann.
\newblock Geschlossenes {S}ystem physikalischer und chemischer
  {D}ifferentialgesetze.
\newblock {\em Sitzungsber. der kaiserliche Akad. Wiss. Wien (IIa)},
  120:385--530, 1911.

\vspace*{-.4em}\bibitem{JiaFis17ADRD}
Y.~Jiao and J.~Fish.
\newblock Is an additive decomposition of a rate of deformation and objective
  stress rates pass\'e?
\newblock {\em Comput. Methods Appl. Mech. Engrg.}, 327:196--225, 2017.

\vspace*{-.4em}\bibitem{KruRou19MMCM}
M.~Kru\v{z}\'{\i}k and T.~Roub{\'{\i}}{\v{c}}ek.
\newblock {\em Mathematical Methods in Continuum Mechanics of Solids}.
\newblock Springer, Cham/Switzerland, 2019.

\vspace*{-.4em}\bibitem{LiuHon99NCHM}
C.-S. Liu and H.-K. Hong.
\newblock Non-oscillation criteria for hypoelastic models under simple shear
  deformation.
\newblock {\em J. Elast.}, 57:201--241, 1999.

\vspace*{-.4em}\bibitem{MelTes69ETCI}
M.~Mellor and R.~Testa.
\newblock Effect of temperature on the creep of ice.
\newblock {\em J. Glaciology}, 8:131--145, 1969.

\vspace*{-.4em}\bibitem{MeXiBrMe03ESRC}
A.~Meyers, H.~Xiao, and O.T. Bruhns.
\newblock Elastic stress ratchetting and corotational stress rates.
\newblock {\em Technische Mechanik}, 23:92--102, 2003.

\vspace*{-.4em}\bibitem{MieRou15RIST}
A.~Mielke and T.~Roub{\'\i}{\v{c}}ek.
\newblock {\em Rate-Independent Systems -- Theory and Application}.
\newblock Springer, New York, 2015.

\vspace*{-.4em}\bibitem{MoDuMu03LIPF}
L.~Moresi, F.~Dufour, and H.~Muhlhaus.
\newblock A {L}agrangian integration point finite element method for large
  deformation modeling of viscoelastic geomaterials.
\newblock {\em J. Comput. Phys.}, 184:476--497, 2003.

\vspace*{-.4em}\bibitem{NeNoSi89GSIC}
J.~Ne\v{c}as, A.~Novotn\'y, and M.~\v{S}ilhav\'y.
\newblock Global solution to the ideal compressible heat conductive multipolar
  fluid.
\newblock {\em Comment. Math. Univ. Carolinae}, 30:551--564, 1989.

\vspace*{-.4em}\bibitem{NecRuz92GSIV}
J.~Ne\v{c}as and M.~R{\accent23u}\v{z}i\v{c}ka.
\newblock Global solution to the incompressible viscous-multipolar material
  problem.
\newblock {\em J. Elasticity}, 29:175--202, 1992.

\vspace*{-.4em}\bibitem{PGYD20SHMM}
C.~{Petrini et al.}
\newblock Seismo-hydro-mechanical modelling of the seismic cycle: Methodology
  and implications for subduction zone seismicity.
\newblock {\em Tectonophysics}, 791:Art.no.228504, 2020.

\vspace*{-.4em}\bibitem{Prag61EDDS}
W.~Prager.
\newblock An elementary discussion of definitions of stress rate.
\newblock {\em Quarterly Appl. Math.}, 18:403--407, 1961.

\vspace*{-.4em}\bibitem{RodUrb99DSPA}
J.F. Rodrigues and J.M. Urbano.
\newblock On a {D}arcy-{S}tefan problem arising in freezing and thawing of
  saturated porous media.
\newblock {\em Continuum Mech. Thermodyn.}, 11:181--191, 1999.

\vspace*{-.4em}\bibitem{RodUrb02TDCS}
J.F. Rodrigues and J.M. Urbano.
\newblock On a three-dimensional convective {S}tefan problem for a
  non-{N}ewtonian fluid.
\newblock In Videman~J.H. Sequeira~A., da Veiga~H.B., editor, {\em Applied
  Nonlinear Analysis}, pages 457--468, Boston, 2002.

\vspace*{-.4em}\bibitem{Roub13NPDE}
T.~Roub{\'\i}{\v{c}}ek.
\newblock {\em Nonlinear Partial Differential Equations with Applications}.
\newblock Birkh\"auser, Basel, 2nd edition, 2013.

\vspace*{-.4em}\bibitem{Roub14NRSD}
T.~Roub{\'{\i}}{\v{c}}ek.
\newblock A note about the rate-and-state-dependent friction model in a
  thermodynamical framework of the {B}iot-type equation.
\newblock {\em Geophysical J. Intl.}, 199:286--295, 2014.

\vspace*{-.4em}\bibitem{Roub19MDDP}
T.~Roub{\'{\i}}{\v{c}}ek.
\newblock Models of dynamic damage and phase-field fracture, and their various
  time discretisations.
\newblock In M.~Hintterm\"uller and J.-F. Rodrigues, editors, {\em Topics in
  Applied Analysis and Optimisation}, pages 363--396. Springer, 2019.

\vspace*{-.4em}\bibitem{Roub??QISC}
T.~Roub{\'{\i}}{\v{c}}ek.
\newblock From quasi-incompressible to semi-compressible fluids.
\newblock {\em Disc. Cont. Dynam. Syst. S}, 14:4069--4092, 2021.

\vspace*{-.4em}\bibitem{Roub21TCMP}
T.~Roub{\'{\i}}{\v{c}}ek.
\newblock Thermodynamically consistent model for poroelastic rocks towards
  tectonic and volcanic processes and earthquakes.
\newblock {\em Geophysical J. Intl.}, 227:1893--1904, 2021.

\vspace*{-.4em}\bibitem{RouTom21CMPE}
T.~Roub{\'{\i}}{\v{c}}ek and G.~Tomassetti.
\newblock A convective model for poro-elastodynamics with damage and fluid flow
  towards {E}arth lithosphere modelling.
\newblock {\em Cont. Mech. Thermodynam.}, 33:2345--2361, 2021.

\vspace*{-.4em}\bibitem{Ruin83SISV}
A.L. Ruina.
\newblock Slip instability and state variable friction laws.
\newblock {\em J. Geophys. Res.}, 88:10,359--10,370, 1983.

\vspace*{-.4em}\bibitem{Schi00SCRC}
G.~Schimperna.
\newblock Some convergence results for a class of nonlinear phase-field
  evolution equations.
\newblock {\em J. Math. Anal. Appl.}, 250:406--434, 2000.

\vspace*{-.4em}\bibitem{SchDuv09CFI}
E.M. Schulson and P.~Duval.
\newblock {\em Creep and Fracture of Ice}.
\newblock Cambridge Univ. Press, Leiden, 2009.

\vspace*{-.4em}\bibitem{Stef89PTW}
J.~Stefan.
\newblock \"{U}ber einige {P}robleme der {T}heorie der {W}\"armeleitung.
\newblock {\em S.B. Wien Akad. Mat. Natur.}, 98:473--484, 1889.

\vspace*{-.4em}\bibitem{Tema69ASEN}
R.~Temam.
\newblock Sur l'approximation de la solution des \'equations de
  {N}avier-{S}tokes par la m\'ethode des pas fractionnaires ({I}).
\newblock {\em Archive Ration. Mech. Anal.}, 32:135--153, 1969.

\vspace*{-.4em}\bibitem{Toma21ITST}
G.~Tomassetti.
\newblock An interpretation of {T}emam's stabilization term in the
  quasi-incompressible {N}avier-{S}tokes system.
\newblock {\em Applications in Engr. Sci.}, 5:Art.no.\,100028, 2021.

\vspace*{-.4em}\bibitem{Visi85SPPH}
A.~Visintin.
\newblock Stefan problem with phase relaxation.
\newblock {\em IMA J. Appl. Math.}, 34:225--245, 1985.

\vspace*{-.4em}\bibitem{Visi85SSEP}
A.~Visintin.
\newblock Supercooling and superheating effects in phase transitions.
\newblock {\em IMA J. Appl. Math.}, 35:233--256, 1985.

\vspace*{-.4em}\bibitem{Visi08ISTP}
A.~Visintin.
\newblock Introduction to {S}tefan-type problems. {C}hap.8 in {H}andbook of
  {D}ifferential {E}quations {IV} ({C.M.\,Dafermos, M.\,Pokorn\'y}, eds.).
\newblock pages 377--484. North-Holland, Amsterdam, 2008.

\vspace*{-.4em}\bibitem{XiBrMe00CORF}
H.~Xiao, O.T. Bruhns, and A.~Meyers.
\newblock The choice of objective rates in finite elastoplasticity: general
  results on the uniqueness of the logarithmic rate.
\newblock {\em Proc. R. Soc. Lond. A}, 456:1865--1882, 2000.

\vspace*{-.4em}\bibitem{XiBrMe06EPSD}
H.~Xiao, O.T. Bruhns, and A.~Meyers.
\newblock Elastoplasticity beyond small deformations.
\newblock {\em Acta Mech.}, 182:31--111, 2006.

\vspace*{-.4em}\bibitem{XuShi97SPCJ}
X.~Xu and M.~Shillor.
\newblock The {S}tefan problem with convection and {J}oule's heating.
\newblock {\em Advances in Differential Equations}, 2:667--691, 1997.

\vspace*{-.4em}\bibitem{Zare03FPTR}
S.~Zaremba.
\newblock Sur une forme perfection\'ees de la th\'eorie de la relaxation.
\newblock {\em Bull. Int. Acad. Sci. Cracovie}, pages 594--614, 1903.

\end{thebibliography}
\end{document}